\documentclass[11pt,a4paper]{article}

\usepackage[T1]{fontenc}
\usepackage{lmodern}
\usepackage[margin=28mm]{geometry}
\usepackage{amsmath,amssymb,amsthm}
\usepackage{microtype}
\usepackage[hidelinks]{hyperref}

\numberwithin{equation}{section}
\newtheorem{theorem}{Theorem}[section]
\newtheorem{proposition}[theorem]{Proposition}
\newtheorem{lemma}[theorem]{Lemma}
\theoremstyle{definition}
\newtheorem{definition}[theorem]{Definition}
\theoremstyle{remark}
\newtheorem{remark}[theorem]{Remark}
\newcommand{\R}{\mathbb R}
\newcommand{\C}{\mathbb C}
\newcommand{\F}{\mathbb F}
\newcommand{\Z}{\mathbb Z}
\newcommand{\sVec}{\mathsf{sVec}_{k}}
\newcommand{\sVecG}{\mathsf{sVec}_{G}}

\newcommand{\Mor}{\mathsf{Mor}_{k}}
\newcommand{\Aut}{\operatorname{Aut}}
\newcommand{\Sec}{\operatorname{Sec}}
\newcommand{\one}{\mathbf 1}
\newcommand{\End}{\operatorname{End}}
\newcommand{\Hom}{\operatorname{Hom}}
\newcommand{\Nat}{\operatorname{Nat}}
\newcommand{\Cl}{\operatorname{Cl}}
\newcommand{\BW}{\operatorname{BW}}

\newcommand{\SH}{\mathrm{SH}}
\newcommand{\Tcal}{\mathcal T}
\newcommand{\W}{\mathcal W}
\newcommand{\Tw}{\mathrm{Tw}}
\newcommand{\Cmp}{N}
\newcommand{\ev}{\operatorname{ev}}
\newcommand{\coev}{\operatorname{coev}}
\newcommand{\id}{\mathrm{id}}
\newcommand{\gt}{\mathbin{\widehat\otimes}}
\newcommand{\simeqm}{\simeq_{\otimes}}
\newcommand{\simeqG}{\simeq_{G}}
\newcommand{\EndCat}{\mathsf{End}}
\newcommand{\Ad}{\operatorname{Ad}}
\title{Weakly Pointed Fusion $\Pi$-supercategories\\
over fields of characteristic not two}
\author{Kotaro Kawasumi}
\date{\today}
\hypersetup{
  pdftitle={Weakly Pointed Fusion Pi-supercategories over fields of characteristic not two},
  pdfauthor={Kotaro Kawasumi}
}

\begin{document}

\maketitle

\begin{abstract}
We study a class of \emph{weakly pointed} fusion $\Pi$-supercategories over a field $k$ of characteristic not equal to $2$, with unit $\one$ satisfying $\End(\one) = k$.
Assuming the existence of a \emph{Brauer--Wall reference category} over $k$, we classify these categories up to monoidal equivalence.
We show that Brauer--Wall reference categories exist for $k = \mathbb C, \mathbb R$ but not for $k = \mathbb F_5$.
For each finite group $G$, we also study a product of these categories equipped with a $G$-grading and define a group $\mathcal T_k(G)$, giving an explicit formula for the multiplication in terms of group cochains for $k = \mathbb C$ and $\mathbb R$.
\end{abstract}

\tableofcontents

\section{Introduction}

Let $k$ be a field of characteristic not equal to $2$.
All vector spaces, algebras, and categories are over $k$ unless stated
otherwise.
All $k$-vector spaces are assumed to be finite-dimensional.
Elements of groups are written $g,h,\ell,m$; the letter $k$ is reserved
for the field.

The general theory of monoidal categories enriched in a braided monoidal
category was developed in
\cite{MorrisonPenneys,MorrisonPenneysPlavnik}.
Fusion categories over non-algebraically closed fields were studied
in \cite{Sanford}.
The extension theory for non-enriched fusion categories over algebraically closed fields in \cite{EtingofNikshychOstrik}
was further developed for enriched fusion categories over algebraically closed fields in
\cite{JonesMorrisonPenneysPlavnik, DavydovNikshych} and for non-enriched fusion categories over non-algebraically closed
fields in \cite{SanfordBrauerPicard}.

The present paper studies a narrow classification problem of enriched fusion categories over non-algebraically closed fields.
Let $\mathsf{sVec}_k$ be the symmetric monoidal category of super vector spaces over $k$.
A monoidal $\Pi$-supercategory is a monoidal $\mathsf{sVec}_k$-enriched category equipped
with a parity object $\pi$ and an odd isomorphism $\pi\to\one$ as defined in \cite{BrundanEllis}.
We call such a category a fusion $\Pi$-supercategory if the category of even morphisms is a fusion category over $k$ in the sense of \cite{Sanford}.
It is \emph{weakly pointed} if every simple object $X$ satisfies
\[
 X\otimes X^*\cong\one^{\oplus p}\oplus\pi^{\oplus q}
 \qquad\text{for some }p,q\geq0~\text{with }p+q > 0.
\]
We classify weakly pointed fusion $\Pi$-supercategories $\mathcal A$
with split unit $\End_{\mathcal A}(\one)\cong k$ up to monoidal equivalence.
Over $\C$, this classification problem is covered by the general
theory of \cite{JonesMorrisonPenneysPlavnik, DavydovNikshych}.
In the present paper, we give an explicit treatment over fields that need not
be algebraically closed.

A physical motivation for studying this class comes from anomalies of finite group symmetries in
$1+1$-dimensional systems.
Let $G$ be a finite group.
In physics, bosonic anomalies in $1+1$-dimensions are described
by $H^3(BG;U(1))$
\cite{ChenGuLiuWen,ChenGuLiuWenScience}.
On the mathematical side, this agrees with the classification
of pointed fusion categories with Picard group $G$.
For fermionic anomalies in $1+1$-dimensions, the approach
initiated in \cite{GuWen} was followed by the development of
a general connection with algebraic topology
\cite{Kapustin,FreedHopkins}.
Based on this connection, a general description was obtained
in \cite{BrumfielMorgan}.
On the other hand, the classification of pointed fusion
supercategories over $\C$ gives a group that differs in general from
$H^3(BG;U(1))$, but can still be a proper subgroup of the group
described in \cite{BrumfielMorgan}.
This motivates the study of weakly pointed fusion
$\Pi$-supercategories over $\mathbb C$.
Time-reversal symmetry also motivates considering the real case
\cite{GaiottoJohnsonFreyd}.

\paragraph{Lean formalization.}
A Lean~4 formalization of the real case is available
in~\cite{RealRModLean}.
It includes the construction of a Brauer--Wall reference category,
classification up to monoidal superequivalence, and stacking
for a fixed finite grading group, together with a splitting
of the degree map.
The proof library was successfully rebuilt from a fresh clone
using Lean~4.31.0.
Build instructions and verification records are included
in the repository.

\subsection{Summary of definitions and results}
\label{subsec:summary}

We list the main definitions and results of this paper. 
See Section~\ref{subsec:weakly-pointed-supercategories} for the basic definitions of supercategories.
For the proofs of the following propositions, we need some basic knowledge of central simple superalgebras and their bimodules, collected in Section~\ref{sec:superalgebras}. 

\begin{definition}[Definition~\ref{def:weakly-pointed}]
\label{intro:def:weakly-pointed}
A \emph{semisimple tensor $\Pi$-supercategory} is a monoidal
$\Pi$-supercategory $\mathcal C$ whose underlying category
$\underline{\mathcal C}$ is an essentially small semisimple $k$-linear
abelian category with finite-dimensional Hom spaces, is rigid, and has
simple unit $\one$.
A semisimple tensor $\Pi$-supercategory $\mathcal C$ is \emph{fusion} if
it has finitely many isomorphism classes of simple objects.
An object $X$ of a semisimple tensor $\Pi$-supercategory $\mathcal C$ is
\emph{weakly invertible} if
\begin{equation}
  X\otimes X^*\cong \one^{\oplus p}\oplus\pi^{\oplus q}
  \qquad\text{for some }p,q\geq 0~\text{with }p+q > 0,
\end{equation}
where $X^*$ is a left dual of $X$ in $\underline{\mathcal C}$.
A semisimple tensor $\Pi$-supercategory is \emph{weakly pointed} if
every simple object is weakly invertible, and has \emph{split unit} if
$\End_{\mathcal C}(\one)\cong k$.
We write $\W_k$ for the class of weakly pointed fusion
$\Pi$-supercategories over $k$ with split unit.
\end{definition}

\begin{definition}[Definitions~\ref{def:brauer-wall-group}, \ref{def:components}, and~\ref{def:brauer-wall-reference}]
\label{intro:def:brauer-wall-reference}
Let $\mathcal C$ be a weakly pointed semisimple tensor
$\Pi$-supercategory with split unit.
A simple object $X\in\mathcal C$ has a \emph{Brauer--Wall type}
$r\in\BW(k)$ if the Morita equivalence class of $\End_\mathcal{C}(X)$ is $r$.
We call $\mathcal C$ a \emph{Brauer--Wall reference category} if every
$r\in\BW(k)$ occurs as the type of a simple object, and two simple
objects are isomorphic if and only if their
types agree.
\end{definition}

\begin{theorem}[Theorem~\ref{thm:three-fields}]
\label{intro:thm:reference-existence}
Up to monoidal superequivalence, the Brauer--Wall reference categories
over the three fields are as follows.

(1) Over $\C$ there are four classes of reference categories.

(2) Over $\R$ there are two classes of reference categories.

(3) Over $\F_5$ no reference category exists.
\end{theorem}

\begin{definition}[Definition~\ref{def:twist}]
\label{intro:def:twist}
A \emph{twist} is a quadruple $t=(G,\rho,\nu,\omega)$, where $G$ is a
finite group, $\rho\in Z^1(G,\BW(k))$,
$\nu\in Z^2(G,\Z/2\Z)$, and $\omega\in C^3(G,k^\times)$, such that
\begin{equation}
  \delta\omega=(-1)^{\nu\smile\nu}.
\end{equation}
We write $\Tw_k$ for the class of all twists.
\end{definition}

\begin{definition}[Definition~\ref{def:gauge}]
\label{intro:def:gauge}
Let $s=(G,\rho,\nu,\omega)$ and $t=(H,\rho',\nu',\omega')$ be twists.
A \emph{gauge} from $s$ to $t$ is a triple $(\phi,\lambda,\beta)$,
where $\phi:G\to H$ is a group isomorphism,
$\lambda\in C^1(G,\Z/2\Z)$, and $\beta\in C^2(G,k^\times)$, such that
\begin{equation}
  \phi^*\rho'=\rho,\qquad
  \phi^*\nu'=\nu+\delta\lambda,\qquad
  \phi^*\omega'
  =\omega(\delta\beta)^{-1}
    (-1)^{(\phi^*\nu')\smile\lambda+\lambda\smile\nu}.
\end{equation}
The twists $s$ and $t$ are \emph{gauge equivalent}, written $s\sim t$,
if there is a gauge from $s$ to $t$.
\end{definition}

\begin{theorem}[Definition~\ref{def:pullback-construction}, Proposition~\ref{prop:pullback-laws}, and Theorem~\ref{thm:classification}]
\label{intro:thm:classification}
For every Brauer--Wall reference category $\mathcal C$ over $k$, there exists a map $\mathcal P_{\mathcal C}: \Tw_k \to \mathcal W_k$ with the following properties:

(1) for every $\mathcal A \in \mathcal W_k$, there exists $t \in \Tw_k$ such that $\mathcal A \simeqm \mathcal P_\mathcal{C}(t)$,

(2) $\mathcal P_\mathcal{C}(t) \simeqm \mathcal P_\mathcal{C}(s)$ if and only if $s \sim t$. 
\end{theorem}

\begin{definition}[Definitions~\ref{def:stacking} and~\ref{def:labeled-pullbacks}]
\label{intro:def:stacking}
Let $G$ be a finite group.
A category $\mathcal A\in\W_k$ is \emph{$G$-labeled} if it is equipped
with a group isomorphism $\lambda_{\mathcal A}:G\to G_{\mathcal A}$, where $G_\mathcal{A}$ is the group of isomorphism classes of simple objects of $\mathcal A$.
We write $\W_k(G)$ for the class of $G$-labeled categories.
For $\mathcal A,\mathcal B\in\W_k(G)$, we write
$\mathcal A\simeqG\mathcal B$ if there is a label preserving monoidal superequivalence. 

Let $\Tw_k(G)$ be the set of twists with group $G$.
For $s,t\in\Tw_k(G)$, we write $s\sim_Gt$ if there exists a gauge
$(\id_G,\lambda,\beta)$ from $s$ to $t$.
For a Brauer--Wall reference category $\mathcal C$ and $t\in\Tw_k(G)$,
the category $\mathcal P_{\mathcal C}(t)$ is $G$-labeled.
\end{definition}

\begin{proposition}[Definition~\ref{def:stacking}, Proposition~\ref{thm:stack-existence}, and Theorem~\ref{thm:labeled-stacking-group}]
  For $\mathcal A, \mathcal B \in \mathcal W_k(G)$, we can define the \emph{stack} $\mathcal S(\mathcal A,\mathcal B) \in \mathcal W_k(G)$.
  Let $[\mathcal A]$ denote the equivalence class under $\simeqG$.
  If we define $[\mathcal A] \cdot [\mathcal B] = [\mathcal S(\mathcal A,\mathcal B)]$, then $\mathcal W_k(G)/\simeqG$ forms an abelian group. 
\end{proposition}

\begin{proposition}[Definition~\ref{def:cochain-operations}, Theorem~\ref{thm:stacking-group-categorical}, and Propositions~\ref{prop:standard-exchange-data} and~\ref{prop:standard-products}]
  Let $\mathcal C$ be a reference category over $k$. Choosing an \emph{exchange data} $(R,\beta_1,\beta_2)$ of $\mathcal C$ (Definition~\ref{def:exchange-datum}), we can define $s\star t \in \Tw_k(G)$ for $s,t\in \Tw_k(G)$ satisfying
  \begin{equation}
    \mathcal S(\mathcal P_\mathcal{C}(s),\mathcal P_\mathcal{C}(t)) \simeqG \mathcal P_\mathcal{C}(s\star t).
  \end{equation}
  If we define $[s]\cdot [t] = [s\star t]$, then $\mathcal T_k(G) = \Tw_k(G)/\sim_G$ forms an abelian group isomorphic to $\mathcal W_k(G)/\simeqG$.

  (1) For $k = \mathbb C$, there exist a reference category $\mathcal C$ over $\mathbb C$ and an exchange data of $\mathcal C$ such that
  \begin{equation}
\begin{aligned}
  s\star t&=\bigl(G,\ \rho_s+\rho_t,\ \nu_s+\nu_t+y\smile x,\
    \omega_s\omega_t\,i^{\,y\smile y\smile x}\,(-1)^{\Xi^{\C}_{s,t}}\bigr),\\
  \Xi^{\C}_{s,t}&=\nu_s\smile_1\nu_t+(y\smile x)\smile_1(\nu_s+\nu_t)
    +y\smile(yx)\smile x,
\end{aligned}
\end{equation}
where $x = \rho_s$, $y = \rho_t$, and $(yx)(h)=y(h)x(h)$.

  (2) For $k = \mathbb R$, identifying $\BW(\mathbb R)\cong\mathbb Z/8\mathbb Z$ as in Definition~\ref{def:cyclic-cochains}, there exist a reference category $\mathcal C$ over $\mathbb R$ and an exchange data of $\mathcal C$ such that
  \begin{equation}
\begin{aligned}
  s\star t&=\bigl(G,\ \rho_s+\rho_t,\ \nu_s+\nu_t+y\smile x,\
    \omega_s\omega_t(-1)^{\Xi^{\R}_{s,t}}\bigr),\\
  \Xi^{\R}_{s,t}
  &= \nu_s\smile_1\nu_t
    +(y\smile x)\smile_1(\nu_s+\nu_t)\\
  &\qquad+y\smile (yx)\smile x + y\smile y\smile (q\circ\rho_s),
\end{aligned}
\end{equation}
where $x = \rho_s\bmod 2$, $y = \rho_t\bmod 2$, $(yx)(h)=y(h)x(h)$, and $q: \mathbb Z/8 \mathbb Z \to \mathbb Z/2\mathbb Z$ is defined by
\begin{equation}
  q(a) = \frac{\tilde a(\tilde a+1)}{2} \bmod 2
\end{equation}
for the representative $\tilde a\in\{0,\ldots,7\}$ of $a$.
\end{proposition}

\begin{theorem}[Proposition~\ref{prop:supercohomology} and Theorems~\ref{thm:stacking-group}, \ref{thm:complex-stacking-group}, and~\ref{thm:real-stacking-group}]
  Assume that $k$ admits a Brauer--Wall reference category. The map
\begin{equation}
  d:\Tcal_k(G)\longrightarrow Z^1(G,\BW(k)),\qquad [s]\longmapsto\rho_s,
\end{equation}
is a surjective group homomorphism with kernel $\SH_k(G)$, and
there is a short exact sequence
\begin{equation}
\label{intro:eq:degree-sequence}
  0\longrightarrow\SH_k(G)\longrightarrow\Tcal_k(G)
  \longrightarrow Z^1(G,\BW(k))\longrightarrow0.
\end{equation}

(1) For $k=\mathbb C$, the short exact sequence does not split if there exists $\rho \in Z^1(G,\mathbb Z/2\mathbb Z)$ not the reduction modulo $2$ of a homomorphism $G\to \mathbb Z/4\mathbb Z$.
For $G=\Z/2\Z$, we have $\SH_\C(\Z/2\Z)\cong\Z/4\Z$ and
$\Tcal_\C(\Z/2\Z)\cong\Z/8\Z$.

(2) For $k = \mathbb R$, the short exact sequence splits.
\end{theorem}

\section{Central Simple Superalgebras}
\label{sec:superalgebras}

In this section, we collect basic knowledge of central simple superalgebras and their bimodules needed for the rest of the paper. 

\begin{definition}
\label{def:superalgebras}
A \emph{super vector space} is a $\Z/2\Z$-graded $k$-vector space
$X=X_0\oplus X_1$. A vector $x\in X_0\cup X_1$ is \emph{homogeneous},
of \emph{degree} $|x|\in\Z/2\Z$; the vectors in $X_0$ are \emph{even}
and those in $X_1$ are \emph{odd}.
A $k$-linear map $f:X\to Y$ between super vector spaces has
\emph{degree} $d\in\Z/2\Z$ if $f(X_i)\subseteq Y_{i+d}$ for $i\in\Z/2\Z$.
The tensor product $X\otimes_{k}Y$ of super vector spaces is graded by
$(X\otimes_{k}Y)_a=\bigoplus_{i+j=a}X_i\otimes_{k}Y_j$.
For a super vector space $X$, $\Pi X$ denotes $X$ with the reversed
grading, $(\Pi X)_a=X_{a+1}$. 
The identity map $\Pi: X \to \Pi X, x \mapsto \Pi x = x$ is an odd isomorphism. 

A \emph{superalgebra} is a super vector space $A$ with an associative
$k$-bilinear multiplication $A\times A\to A$ and a unit $1\in A_0$
such that $A_iA_j\subseteq A_{i+j}$ for $i,j\in\Z/2\Z$.
A \emph{homomorphism} of superalgebras is a $k$-linear map of degree
$0$ which preserves the multiplication and the unit.
We write $A \cong B$ if superalgebras $A$ and $B$ are isomorphic.
The \emph{supercenter} of $A$ is $Z_{\mathrm s}(A)=
Z_{\mathrm s}(A)_0\oplus Z_{\mathrm s}(A)_1$, where
\begin{equation}
\label{eq:supercenter}
 Z_{\mathrm s}(A)_d
 =\{z\in A_d:za=(-1)^{d|a|}az
       \text{ for every homogeneous }a\in A\}.
\end{equation}
A superalgebra $A$ is \emph{central} if $Z_\mathrm{s}(A) = k$,
\emph{simple} if it is nonzero and has no nonzero proper graded
two-sided ideals, and a \emph{division superalgebra} if it is nonzero
and every nonzero homogeneous element is invertible. 
The \emph{opposite superalgebra} $A^{\mathrm{sop}}$ of $A$ is the super vector space $A$ with the same
grading and the multiplication
\begin{equation}
\label{eq:superopposite}
 a\cdot_{\mathrm{sop}}b=(-1)^{|a||b|}ba.
\end{equation}
For a superalgebra $A$ and $p,q\geq0$ with $p+q>0$, the
\emph{matrix superalgebra} $M_{p|q}(A)$ is
$\operatorname{Mat}_{p+q}(A)$ with ordinary matrix multiplication and grading
\begin{equation}
\label{eq:matrix-superalgebra}
 M_{p|q}(A)_d
 =\{(a_{ij}):a_{ij}\in A_{d+\epsilon_i+\epsilon_j}\},
 \qquad
 \epsilon_i=\begin{cases}0&i\leq p,\\1&i>p.\end{cases}
\end{equation}
The \emph{graded tensor product} $A\gt B$ of superalgebras $A$ and $B$
is the super vector space $A\otimes_{k}B$ with the unit $1\gt1$ and the
multiplication
\begin{equation}
\label{eq:algebra-tensor-convention}
 (a\gt b)(a'\gt b')=(-1)^{|b||a'|}aa'\gt bb'
\end{equation}
for homogeneous $a,a'\in A$ and $b,b'\in B$.
\end{definition}

\begin{lemma}
\label{lem:graded-division-idempotents}
Let $D$ be a division superalgebra. Then $D_0$ is a division algebra,
and the only even idempotents of $D$ are $0$ and $1$.
\end{lemma}
\begin{proof}
Every nonzero element of $D_0$ is invertible, and an invertible
idempotent is $1$.
\end{proof}

\begin{lemma}
\label{lem:matrix-corners}
Let $D$ be a division superalgebra and let $e$ be a nonzero even
idempotent of $M_{p|q}(D)$. If $eM_{p|q}(D)e$ is a division
superalgebra, then it is isomorphic to $D$.
\end{lemma}
\begin{proof}
The corner superalgebra $eM_{p|q}(D)e$ is the endomorphism superalgebra
of the graded image of $e$. Graded freeness over $D$ and the division
assumption force this image to have graded rank one, so
$eM_{p|q}(D)e\cong D$;
see \cite[Section~1.3.2 and Proposition~1.4.1]{Hazrat}.
\end{proof}

\begin{definition}
\label{def:matrix-superalgebras}
Let $A$ be a superalgebra.
A \emph{graded right $A$-module} is a super vector space $V$ with an
even $k$-linear right $A$-action.
For graded right $A$-modules $V$ and $W$, a $k$-linear map $f:V\to W$
of degree $d$ is \emph{right $A$-linear} if $f(vb)=f(v)b$ for $v\in V$
and $b\in A$.
We write $\End_A(V)$ for the superalgebra whose degree-$d$ part
consists of the right $A$-linear maps $V\to V$ of degree $d$, with the composition as the multiplication.
\end{definition}

\begin{definition}
\label{def:graded-bimodules}
Let $A$, $B$, and $C$ be finite-dimensional superalgebras.
A \emph{graded $(A,B)$-bimodule} is a finite-dimensional super vector
space $M$ with even $k$-linear left $A$-action and right $B$-action satisfying $(am)b=a(mb)$.
It is \emph{simple} if it is nonzero and has no nonzero proper
graded subbimodules.
We write $M:A\to B$ for a graded $(A,B)$-bimodule $M$.
For graded $(A,B)$-bimodules $M,M'$, let $\Hom_{A-B}(M,M')$ be the
super vector space whose degree-$d$ part consists of the $k$-linear maps
$u:M\to M'$ of degree $d$ such that
\begin{equation}
\label{eq:morita-homogeneous-maps}
 u(amb)=(-1)^{|u||a|}a\,u(m)b
\end{equation}
for homogeneous $a\in A$, $m\in M$, and $b\in B$.
For $M:A\to B$ and $N:B\to C$, the balanced tensor product
$M\otimes_BN$ is a graded $(A,C)$-bimodule, and for homogeneous
$u\in\Hom_{A-B}(M,M')$ and $v\in\Hom_{B-C}(N,N')$,
\begin{equation}
\label{eq:morita-composition}
 u\otimes_Bv\in\Hom_{A-C}(M\otimes_BN,M'\otimes_BN'),\qquad
 (u\otimes_Bv)(m\otimes_Bn)=(-1)^{|v||m|}u(m)\otimes_Bv(n).
\end{equation}
For $M:A\to B$, $N:B\to C$, and $P:C\to C'$, we use the even bimodule
isomorphisms
\begin{gather}
\label{eq:morita-associator}
 \mathsf a_{M,N,P}:(M\otimes_BN)\otimes_CP\longrightarrow M\otimes_B(N\otimes_CP),\qquad
 (m\otimes_Bn)\otimes_Cp\longmapsto m\otimes_B(n\otimes_Cp),\\
\label{eq:morita-left-unitor}
 \mathsf l_M:A\otimes_AM\longrightarrow M,\qquad
 a\otimes_Am\longmapsto am,\\
\label{eq:morita-right-unitor}
 \mathsf r_M:M\otimes_BB\longrightarrow M,\qquad
 m\otimes_Bb\longmapsto mb.
\end{gather}
For $M:A\to B$ and $N:A'\to B'$, we define the graded
$(A\gt A',B\gt B')$-bimodule $M\gt N$ to be the super vector space
$M\otimes_kN$ with the actions
\begin{equation}
\label{eq:morita-external-actions}
 (a\gt a')(m\gt n)(b\gt b')=(-1)^{|a'||m|+(|a'|+|n|)|b|}amb\gt a'nb'
\end{equation}
for homogeneous $a\in A$, $a'\in A'$, $m\in M$, $n\in N$, $b\in B$, and
$b'\in B'$.
For homogeneous $u\in\Hom_{A-B}(M,M')$ and $v\in\Hom_{A'-B'}(N,N')$, we
define $u\gt v\in\Hom_{A\gt A'-B\gt B'}(M\gt N,M'\gt N')$ by
\begin{equation}
\label{eq:morita-external-maps}
 (u\gt v)(m\gt n)=(-1)^{|v||m|}u(m)\gt v(n).
\end{equation}
For $M:A\to B$, $M':B\to C$, $N:A'\to B'$, and $N':B'\to C'$, the
\emph{interchange isomorphism} is the even isomorphism of graded
$(A\gt A',C\gt C')$-bimodules
\begin{equation}
\label{eq:morita-interchange}
\begin{gathered}
 \chi_{M,N,M',N'}:(M\gt N)\otimes_{B\gt B'}(M'\gt N')\longrightarrow
 (M\otimes_BM')\gt(N\otimes_{B'}N'),\\
 (m\gt n)\otimes(m'\gt n')\longmapsto
 (-1)^{|n||m'|}(m\otimes_Bm')\gt(n\otimes_{B'}n').
\end{gathered}
\end{equation}
An isomorphism of superalgebras $f:A\to B$ is regarded as the graded
$(A,B)$-bimodule ${}_fB_B$, the bimodule $B$ with left $A$-action
through $f$; in particular the swap $\tau_{A,B}:A\gt B\to B\gt A$,
$a\gt b\mapsto(-1)^{|a||b|}b\gt a$, is a graded
$(A\gt B,B\gt A)$-bimodule.
A graded $(A,B)$-bimodule $M$ is a \emph{Morita equivalence} if there
is a graded $(B,A)$-bimodule $N$ with even isomorphisms
$M\otimes_BN\cong A$ and $N\otimes_AM\cong B$.
Superalgebras $A$ and $B$ are \emph{Morita equivalent} if there is a Morita equivalence.
\end{definition}

\begin{proposition}
\label{prop:central-simple-superalgebras}
Let $A$ and $B$ be  central simple superalgebras.

(1) $A\cong M_{p|q}(D)$ for a central division
superalgebra $D$, and $D$ is unique up to isomorphism. 

(2) $A$ and $B$ are Morita equivalent if and only if $A \cong M_{p|q}(D)$ and $B \cong M_{u|v}(D)$ for a central division superalgebra $D$.

(3) Morita equivalence is an equivalence relation.

(4) $A^{\mathrm{sop}}$ is central simple. 

(5) $A\gt B$ is central simple, and
\begin{equation}
\label{eq:enveloping-matrix}
 A\gt A^{\mathrm{sop}}\xrightarrow{\ \sim\ }\End_k(A),
 \qquad a\gt b\longmapsto\bigl(x\mapsto(-1)^{|b||x|}axb\bigr),
\end{equation}
is an isomorphism of superalgebras. 

(6) The even part $A_0$ is a semisimple $k$-algebra.

(7) $Z(A) = k$ or $Z(A) = k \oplus kz$ for an odd element $z$ with $z^2\in k^\times$, where $Z(A)$ is the ordinary center. 
If we define $\epsilon(A) \in \mathbb Z/2\mathbb Z$ by 
\begin{equation}
\label{eq:bw-parity}
 \epsilon(A)=
 \begin{cases}
  0&\dim_kZ(A)=1,\\
  1&\dim_kZ(A)=2,
 \end{cases}
\end{equation}
then $\epsilon(A)$ is invariant under  Morita equivalence, and $\epsilon(A\gt B) = \epsilon(A) + \epsilon(B)$.
\end{proposition}
\begin{proof}
Parts (1)--(2) are graded Artin--Wedderburn and Morita theory
\cite[Proposition~1.4.4, Remark~1.4.8, and Theorems~2.3.7--2.3.8]{Hazrat};
matrix superalgebras preserve the supercenter. Identity, inverse,
and composite bimodules give (3).
Parts (4)--(5), including the sign in \eqref{eq:enveloping-matrix},
are \cite[Chapter~IV, Section~4, Theorem~2.3(3), and Proposition~4.1]{LamQF}.
For (6), use ordinary central simplicity when $A_1=0$ and
\cite[Chapter~IV, Theorems~3.4(2) and~3.8]{LamQF} otherwise.
Finally, (7) follows from the center classification and product
theorem in \cite[Chapter~IV, Theorems~3.4, 3.6, 3.12, and~4.3]{LamQF}.
\end{proof}

\begin{definition}
\label{def:brauer-wall-group}
The \emph{Brauer--Wall group} $\BW(k)$ is the abelian group of graded
Morita equivalence classes of central simple superalgebras given by
\begin{equation}
\label{eq:brauer-wall-group}
 [A]+[B]=[A\gt B],\qquad 0=[k],\qquad -[A]=[A^{\mathrm{sop}}].
\end{equation}
For each $r\in\BW(k)$ fix a central division superalgebra $D_r$ of class
$r$, with $D_0=k$.
The \emph{parity} $\epsilon:\BW(k)\to\Z/2\Z$ is the homomorphism given by \eqref{eq:bw-parity}.
\end{definition}

\begin{proposition}
\label{prop:morita-equivalences}
Let $A$ and $B$ be Morita equivalent central simple superalgebras.
A graded $(A,B)$-bimodule is a Morita equivalence if and only if it is simple.
For Morita equivalences $M,N:A \xrightarrow{\sim} B$, the super vector space
$\Hom_{A-B}(M,N)$ is one-dimensional, and each of its nonzero elements is invertible.
Consequently, $N$ is evenly isomorphic to one of $M$ and
$\Pi M$, where $\Pi M$ has the actions
\begin{equation}
\label{eq:morita-parity-shift}
 a\cdot(\Pi m)\cdot b=(-1)^{|a|}\Pi(amb)
\end{equation}
for homogeneous $a\in A$, $m\in M$, and $b\in B$.
\end{proposition}
\begin{proof}
Choose a Morita equivalence $P:A\to B$ with inverse $Q$.
Tensoring with $Q$ identifies graded $(A,B)$-bimodules with graded
$(B,B)$-bimodules and sends $P$ to $B$.
By \eqref{eq:enveloping-matrix}, the latter category is equivalent
to super vector spaces; see also \cite[Theorem~4.3.1]{Kwok}.
Its simple objects correspond to the two super lines. Hence every
simple $(A,B)$-bimodule is $P$ or $\Pi P$ up to even isomorphism,
and the homogeneous intertwiner space between any two is a super
line with invertible nonzero elements.
Both $P$ and $\Pi P$ are Morita equivalences: the inverse of
$\Pi P$ is $\Pi Q$, using
\[
 \Pi P\otimes_B\Pi Q\longrightarrow P\otimes_BQ,
 \qquad \Pi x\otimes\Pi y\longmapsto(-1)^{|x|}x\otimes y.
\]
Conversely, transport along any Morita equivalence sends it to
the simple regular bimodule $B$.
\end{proof}

\begin{proposition}
\label{prop:swap-parity}
Let $A$ be a central simple superalgebra and let
$\tau:A\gt A\to A\gt A$, $\tau(x\gt y)=(-1)^{|x||y|}y\gt x$, be the swap.
The homogeneous elements $v\in A\gt A$ such that
\begin{equation}
\label{eq:signed-swap-implementation}
 vz=(-1)^{|v||z|}\tau(z)v
 \qquad(z\in A\gt A\text{ homogeneous})
\end{equation}
span a one-dimensional super vector space concentrated in degree
$\epsilon([A])$, and its nonzero elements are invertible.
\end{proposition}
\begin{proof}
For $B=A\gt A$, the solution space is
$\Hom_{B-B}(B,\tau^*B)$ via $v\mapsto(m\mapsto vm)$.
Proposition~\ref{prop:morita-equivalences} makes this a super line
with invertible nonzero elements.
Its parity can be checked after separable scalar extension, where
$A$ is $E=\End(V)$ or $E\gt C$, with $C=k^s[e]/(e^2-1)$ and $e$ odd
\cite[Theorem~3]{Wall}.
For $E$, the super flip of $V\otimes V$ implements the swap and is even.
For $C$, put $e_1=e\gt1$, $e_2=1\gt e$, and $v=e_1-e_2$. Then
\[
 ve_1=-e_2v,\qquad ve_2=-e_1v,\qquad v^2=2,
\]
so $v$ is odd, invertible, and satisfies
\eqref{eq:signed-swap-implementation}. Under the canonical
reordering of $(E\gt C)^{\gt2}$, tensoring the two elements just
constructed gives an odd invertible element satisfying the same
relation for $E\gt C$. Thus the parity is $\epsilon([A])$.
\end{proof}

\section{Supercategories}
\label{sec:supercategories}

In Section~\ref{subsec:weakly-pointed-supercategories}, we collect some basic notions on supercategories and define the class $\mathcal W_k$ of weakly pointed fusion $\Pi$-supercategories over $k$ with split unit.

In Section~\ref{subsec:components-brauer-wall-types}, we show that each category $\mathcal A$ in $\mathcal W_k$ determines a group homomorphism $\rho_\mathcal{A}$ from the group $G_\mathcal{A}$ of isomorphism classes of simple objects to the Brauer--Wall group $\BW(k)$.

In Section~\ref{subsec:twists-pullback}, we introduce the classifying parameters of $\mathcal W_k$ and construct a category from each parameter, assuming the existence of a category $\mathcal C$ such that $\rho_\mathcal{C}$ is an isomorphism.
We call such categories reference categories.
They play a role akin to that of an origin of the space $\mathcal W_k$.
Their existence is not necessary for the classification, and we can use other $2$-categorical data as the origin, as in Remark~\ref{rem:other-systems}, but the resulting classifying parameters will change.
Using reference categories as the origin makes the description of the parameters simple, compatible with the literature, and independent of the data of the origin.
Also, the definition of a reference category is simpler than that of other objects that could serve as the origin.
The existence of reference categories over $k$ is studied in Sections~\ref{sec:reference-category} and~\ref{sec:fields}.

Finally, in Section~\ref{subsec:reconstruction-classification}, we extract a twist from any category in $\mathcal W_k$ and prove that
$t\mapsto\mathcal P_{\mathcal C}(t)$ induces a bijection between gauge
equivalence classes of twists and monoidal superequivalence classes in
$\mathcal W_k$.

\subsection{Weakly pointed fusion \texorpdfstring{$\Pi$}{Pi}-supercategories}
\label{subsec:weakly-pointed-supercategories}

In this section, we give a precise definition of weakly pointed fusion $\Pi$-supercategories over $k$.

\begin{definition}
\label{def:super-vector-spaces}
The category $\sVec$ is the symmetric monoidal category of
super vector spaces (Definition~\ref{def:superalgebras}) and
degree-preserving linear maps,
with the tensor product
\begin{equation}
\label{eq:graded-tensor}
\begin{aligned}
  (X\otimes Y)_a
  &= \bigoplus_{i+j=a} X_i\otimes_{k}Y_j,\\
  (f\otimes g)(x\otimes y)
  &= f(x)\otimes g(y)
  \qquad(f:X\to X',\ g:Y\to Y'),
\end{aligned}
\end{equation}
the unit $k$ concentrated in degree zero, and the usual associativity
and unit constraints. Its braiding is given by
\begin{equation}
\label{eq:koszul}
  \beta_{X,Y}(x\otimes y)
  = (-1)^{|x||y|}y\otimes x,
\end{equation}
for homogeneous vectors $x\in X$ and $y\in Y$.
\end{definition}

\begin{definition}
\label{def:supercategories}
A \emph{supercategory} is an $\sVec$-enriched category $\mathcal C$.
Composition has the order of \cite[\S1.2]{Kelly}:
\begin{equation}
\label{eq:enriched-composition}
  \mathcal C(Y,Z)\otimes_{k}\mathcal C(X,Y)
  \longrightarrow\mathcal C(X,Z),
  \qquad g\otimes f\longmapsto g\circ f.
\end{equation}
The \emph{endomorphism superalgebra} $\End_{\mathcal C}(X)$ is
$\mathcal C(X,X)$, with its given grading and multiplication
\begin{equation}
\label{eq:endomorphism-algebra}
  fg=f\circ g,\qquad 1=\id_X.
\end{equation}
A \emph{superfunctor} is an $\sVec$-enriched functor.
A \emph{supernatural transformation} of degree $d\in\Z/2\Z$
between superfunctors $F,F':\mathcal C\to\mathcal D$ is a family of
morphisms $\theta_X:F(X)\to F'(X)$ of degree $d$ satisfying
\begin{equation}
\label{eq:homogeneous-supernaturality}
  \theta_Y\circ F(f)=(-1)^{d|f|}F'(f)\circ\theta_X
\end{equation}
for every homogeneous morphism $f:X\to Y$.
We write $\Nat(F,F')$ for the super vector space whose
degree-$d$ part consists of supernatural
transformations $F\Rightarrow F'$ of degree $d$.
A superfunctor $F:\mathcal C\to\mathcal D$ is a \emph{superequivalence}
if we have a superfunctor $G:\mathcal D\to\mathcal C$ and even supernatural
isomorphisms $GF\cong\mathrm{Id}_{\mathcal C}$ and
$FG\cong\mathrm{Id}_{\mathcal D}$.
The \emph{underlying category} $\underline{\mathcal C}$ of a supercategory
$\mathcal C$ is obtained by retaining all objects and only the even morphisms:
\begin{equation}
\label{eq:even-category}
  \underline{\mathcal C}(X,Y)=\mathcal C(X,Y)_0.
\end{equation}
A supercategory $\mathcal C$ is \emph{additive} if its underlying
category $\underline{\mathcal C}$ is additive.
An object of $\mathcal C$ is \emph{simple} if it is simple in $\underline{\mathcal C}$.
\emph{Direct sums} in $\mathcal C$ are taken in $\underline{\mathcal C}$.
An object $X$ is an \emph{even retract} of $Y$ if there exist even morphisms
$i:X\to Y$ and $p:Y\to X$ such that $p\circ i=\id_X$.
\end{definition}

\begin{definition}
\label{def:monoidal-supercategories}
For supercategories $\mathcal C$ and $\mathcal D$,
the \emph{product} $\mathcal C\boxtimes\mathcal D$ is a supercategory with
objects $X\boxtimes Y=(X,Y)$, morphism spaces
\begin{equation}
\label{eq:supercategory-product}
  (\mathcal C\boxtimes\mathcal D)(X\boxtimes Y,X'\boxtimes Y')
  =\mathcal C(X,X')\otimes_{k}\mathcal D(Y,Y'),
\end{equation}
and composition
\begin{equation}
\label{eq:super-interchange}
  (f'\boxtimes g')\circ(f\boxtimes g)
  =(-1)^{|g'||f|}(f'\circ f)\boxtimes(g'\circ g)
\end{equation}
for homogeneous morphisms, as in \cite[(1.1), Definition~1.4]{BrundanEllis}.
Here $f\boxtimes g$ denotes $f\otimes_{k}g$.
A \emph{monoidal supercategory} is a supercategory $\mathcal C$ equipped with a
superfunctor $\otimes:\mathcal C\boxtimes\mathcal C\to\mathcal C$, a unit
object $\one$, and even supernatural associativity and unit isomorphisms
satisfying the pentagon and triangle axioms. 
Its underlying category $\underline{\mathcal C}$ inherits the monoidal structure of $\mathcal C$.
A \emph{monoidal $\Pi$-supercategory} is a monoidal supercategory equipped with
an object $\pi$ and an odd isomorphism $\zeta: \pi\xrightarrow{\sim}\one$.
\end{definition}

\begin{definition}
\label{def:parity-conventions}
Let $\mathcal C$ be a monoidal $\Pi$-supercategory.
For $a\in\Z/2\Z$, put
\begin{equation}
\label{eq:parity-lines}
  P_a=
  \begin{cases}
    \one & a=0,\\
    \pi  & a=1.
  \end{cases}
\end{equation}
Define a superfunctor $\Pi:\mathcal C\to\mathcal C$ by
\begin{equation}
\label{eq:parity-shift-functor}
  \Pi X=\pi\otimes X,\qquad
  \Pi(f)=\id_\pi\otimes f.
\end{equation}
The morphisms
\begin{equation}
\label{eq:parity-shift-iso}
  z_X=l_X\circ(\zeta\otimes\id_X):\Pi X\longrightarrow X
\end{equation}
form an odd supernatural isomorphism
$z:\Pi\Rightarrow\id_{\mathcal C}$.
\end{definition}

\begin{definition}
\label{def:monoidal-superequivalence}
A \emph{monoidal superfunctor} between monoidal supercategories
$\mathcal C$ and $\mathcal D$ is a triple $(F,J,\varphi)$, where
$F:\mathcal C\to\mathcal D$ is a superfunctor, the \emph{tensorator}
\begin{equation}
\label{eq:monoidal-superfunctor-data}
 J_{X,Y}:F(X)\otimes F(Y)\xrightarrow{\ \sim\ }F(X\otimes Y)
\end{equation}
is an even supernatural isomorphism, and the \emph{unit map}
$\varphi:\one_{\mathcal D}\xrightarrow{\sim}F(\one_{\mathcal C})$ is an even
isomorphism, satisfying the usual monoidal functor axioms.
It is \emph{strict monoidal} if $\varphi$ and all $J_{X,Y}$ are identity
morphisms, and a \emph{monoidal superequivalence} if $F$ is a
superequivalence; we write $\mathcal C\simeqm\mathcal D$ when such a
superequivalence exists.
\end{definition}

\begin{proposition}
\label{prop:monoidal-parity}
Let $(F,J,\varphi):\mathcal C\to\mathcal D$ be a monoidal superfunctor
between monoidal $\Pi$-supercategories with parity data
$(\pi_{\mathcal C},\zeta_{\mathcal C})$ and $(\pi_{\mathcal D},\zeta_{\mathcal D})$.
Then there is a unique even isomorphism
$c_F:F(\pi_{\mathcal C})\to\pi_{\mathcal D}$ such that
\begin{equation}
\label{eq:parity-object-comparison}
 \zeta_{\mathcal D}\circ c_F=\varphi^{-1}\circ F(\zeta_{\mathcal C}),
\end{equation}
namely $c_F=\zeta_{\mathcal D}^{-1}\circ\varphi^{-1}\circ F(\zeta_{\mathcal C})$.
\end{proposition}
\begin{proof}
The displayed formula composes two odd isomorphisms and an even one,
so $c_F$ is even. Invertibility of $\zeta_{\mathcal D}$ gives uniqueness.
\end{proof}

\begin{definition}
\label{def:weakly-pointed}
A \emph{semisimple tensor $\Pi$-supercategory} is a monoidal
$\Pi$-supercategory $\mathcal C$ whose underlying category
$\underline{\mathcal C}$ is an essentially small semisimple $k$-linear
abelian category with finite-dimensional Hom spaces, is rigid, and has
simple unit $\one$.
A semisimple tensor $\Pi$-supercategory $\mathcal C$ is \emph{fusion} if
it has finitely many isomorphism classes of simple objects.
An object $X$ of a semisimple tensor $\Pi$-supercategory $\mathcal C$ is
\emph{weakly invertible} if
\begin{equation}
\label{eq:weak-invertibility}
  X\otimes X^*\cong \one^{\oplus p}\oplus\pi^{\oplus q}
  \qquad\text{for some }p,q\geq 0~\text{with }p+q > 0,
\end{equation}
where $X^*$ is a left dual of $X$ in $\underline{\mathcal C}$.
A semisimple tensor $\Pi$-supercategory is \emph{weakly pointed} if
every simple object is weakly invertible, and has \emph{split unit} if
$\End_{\mathcal C}(\one)\cong k$.
We write $\W_k$ for the class of weakly pointed fusion
$\Pi$-supercategories with split unit.
\end{definition}

\begin{remark}
\label{rem:fusion-sanford}
If $k$ is perfect, the underlying category of a fusion
$\Pi$-supercategory is a fusion category over $k$ in the sense of
\cite[Definition~2.9]{Sanford}, since a finite-dimensional semisimple
algebra over a perfect field is separable.
\end{remark}

\subsection{Components and Brauer--Wall types}
\label{subsec:components-brauer-wall-types}

For a weakly pointed semisimple tensor $\Pi$-supercategory $\mathcal A$
with split unit, we show that the isomorphism classes of simple objects
in $\mathcal A$ form a group $G_{\mathcal A}$, and that the Morita
equivalence classes of their endomorphism superalgebras define a group
homomorphism $\rho_{\mathcal A}:G_{\mathcal A}\to\BW(k)$.
Brauer--Wall reference categories are precisely those for which
$\rho_{\mathcal A}$ is an isomorphism.

\begin{lemma}
\label{lem:homogeneous-schur}
Let $\mathcal A$ be a supercategory such that the underlying category
$\underline{\mathcal A}$ is abelian, let $\Pi:\mathcal A\to\mathcal A$ be
a superfunctor, and let $z:\Pi\Rightarrow\id_{\mathcal A}$ be an odd
supernatural isomorphism.

(1) Every nonzero homogeneous morphism between simple objects of
$\mathcal A$ is invertible. If two simple objects of $\mathcal A$ are
isomorphic then there is a homogeneous isomorphism between them.

(2) Suppose moreover that $\underline{\mathcal A}$ is semisimple with
finite-dimensional Hom spaces, let $S$ be a simple object of
$\mathcal A$, let $S'$ be an object of $\mathcal A$, and let
$u:S'\to S$ be an odd isomorphism. An object isomorphic in
$\mathcal A$ to a finite direct sum of copies of $S$ and $S'$ is evenly
isomorphic to $S^{\oplus p}\oplus S'^{\oplus q}$ for some $p,q\geq0$.
\end{lemma}
\begin{proof}
The even isomorphism $z_Xz_{\Pi X}:\Pi^2X\to X$ makes $\Pi$ an
autoequivalence of $\underline{\mathcal A}$. For odd $f:X\to Y$,
$fz_X:\Pi X\to Y$ is even, so ordinary Schur's lemma applies.
Any isomorphism between simples has a nonzero homogeneous component.
In particular,
\begin{equation}
\label{eq:simple-hom-orthogonality}
 X\not\cong Y\ \Longrightarrow\ \mathcal A(X,Y)=0,
 \qquad \End(X)\text{ is a division superalgebra}.
\end{equation}
For (2), decompose the object evenly into simples. Each summand
has a nonzero homogeneous map to $S$ or $S'$, hence is homogeneously
isomorphic to $S$. Composing with $u$ or $u^{-1}$ when necessary
makes this an even isomorphism to $S$ or $S'$.
\end{proof}

\begin{proposition}
\label{prop:simple-tensor}
Let $\mathcal A$ be a weakly pointed semisimple tensor $\Pi$-supercategory with
split unit. For simple objects $X,Y$, the map
\begin{equation}
\label{eq:tensor-end-map}
 \End(X)\gt\End(Y)\xrightarrow{\sim}\End(X\otimes Y),
 \qquad f\gt g\longmapsto f\otimes g
\end{equation}
is an isomorphism of superalgebras.
\end{proposition}
\begin{proof}
Semisimplicity and the split unit make composition
$\mathcal A(U,\one)\times\mathcal A(\one,U)\to k$ a perfect pairing for
every object $U$. Using the left-dual isomorphism
$\End(V)\xrightarrow{\sim}\mathcal A(V^*\otimes V,\one)$,
$f\mapsto\ev_V(\id_{V^*}\otimes f)$, each homogeneous functional
$\varphi\in\End(V)^*$ of degree $d$ is represented by a unique
$\tilde\psi\in\mathcal A(\one,V^*\otimes V)_d$ with
$\varphi(f)=\ev_V(\id_{V^*}\otimes f)\tilde\psi$.
Define $\Sigma_\varphi:\End(V\otimes W)\to\End(W)$ by
\[
 \Sigma_\varphi(h)
 =l_W(\ev_V\otimes\id_W)\alpha_{V^*,V,W}^{-1}
  (\id_{V^*}\otimes h)\alpha_{V^*,V,W}(\tilde\psi\otimes\id_W)l_W^{-1}.
\]
By super interchange, $\Sigma_\varphi(f\otimes g)=(-1)^{d|g|}\varphi(f)g$
for homogeneous $f\in\End(V)$ and $g\in\End(W)$.
Choose a homogeneous basis $(f_i)$ of $\End(V)$ and its dual basis
$(\varphi_i)$. If $\sum_i f_i\otimes g_i=0$ in $\End(V\otimes W)$,
applying $\Sigma_{\varphi_j}$ annihilates both parity components of
$g_j$, so $g_j=0$ for every $j$. Hence
$\End(V)\gt\End(W)\to\End(V\otimes W)$ is injective.

For simple $Y$, choose an even decomposition
$Y\otimes Y^*\cong\one^{\oplus p}\oplus\pi^{\oplus q}$.
Left duality gives
$\mathcal A(X\otimes Y,X\otimes Y)\cong\mathcal A(X,X\otimes(Y\otimes Y^*))$,
whence
\[
 \dim_k\End(Y)=p+q,\qquad
 \dim_k\End(X\otimes Y)=(p+q)\dim_k\End(X).
\]
Thus the injective map is surjective; super interchange makes it
an algebra homomorphism.
\end{proof}

\begin{definition}
\label{def:components}
Let $\mathcal A$ be a weakly pointed semisimple tensor $\Pi$-supercategory with
split unit, and let $G_{\mathcal A}$ be the set of classes of simple
objects of $\mathcal A$ under homogeneous isomorphism in $\mathcal A$;
by Lemma~\ref{lem:homogeneous-schur}, two simple objects have the same
class if and only if they are isomorphic in $\mathcal A$.
For $g\in G_{\mathcal A}$, the \emph{component} $\mathcal A_g$
is the full subsupercategory of objects whose simple summands have class
$g$, including the zero object.
\end{definition}

\begin{definition}
\label{def:brauer-wall-reference}
Let $\mathcal C$ be a weakly pointed semisimple tensor $\Pi$-supercategory with
split unit. A simple object $X\in\mathcal C$ has \emph{Brauer--Wall type}
$r\in\BW(k)$ if there is an isomorphism of superalgebras
\begin{equation}
\label{eq:brauer-wall-type}
  D_r\cong\End_{\mathcal C}(X).
\end{equation}
We call $\mathcal C$ a \emph{Brauer--Wall reference category} if every
$r\in\BW(k)$ occurs as the type of a simple object, and two simple objects
are isomorphic in $\mathcal C$, that is, have the same class in
$G_{\mathcal C}$, if and only if their types agree.
Section~\ref{sec:reference-category} decides when such a category exists.
\end{definition}

\begin{proposition}
\label{prop:component-types}
Every component of a weakly pointed semisimple tensor $\Pi$-supercategory
with split unit has a unique
Brauer--Wall type, given by the type of any of its simple objects.
\end{proposition}
\begin{proof}
For simple $X$, Lemma~\ref{lem:homogeneous-schur} makes $D=\End(X)$ a
division superalgebra. Weak invertibility and
Proposition~\ref{prop:simple-tensor} give
\[
 D\gt\End(X^*)\cong\End(X\otimes X^*)\cong M_{p|q}(k).
\]
The inclusion $z\mapsto z\gt1$ sends $Z_{\mathrm s}(D)$ into
$Z_{\mathrm s}(M_{p|q}(k))=k$, hence $D$ is central.
The unique representative $D_r\cong D$ determines the type.
Conjugation by a homogeneous isomorphism preserves the endomorphism
superalgebra, so this type is constant on a component.
\end{proof}

\begin{proposition}
\label{prop:intrinsic-group}
Let $\mathcal A$ be a weakly pointed semisimple tensor $\Pi$-supercategory with
split unit. The rule
\begin{equation}
\label{eq:intrinsic-product}
 [X][Y]=[S],\qquad S\text{ any simple summand of }X\otimes Y,
\end{equation}
makes $G_{\mathcal A}$ a group, finite when $\mathcal A\in\W_k$, with identity $[\one]$ and
inverse $[X]^{-1}=[X^*]$. Brauer--Wall type defines a homomorphism
\begin{equation}
\label{eq:intrinsic-rho}
 \rho_{\mathcal A}:G_{\mathcal A}\longrightarrow\BW(k),
 \qquad \rho_{\mathcal A}([X])=r\quad\text{if }\End(X)\cong D_r.
\end{equation}
The category $\mathcal A$ decomposes as
\begin{equation}
\label{eq:component-decomposition}
 \mathcal A=\bigoplus_{g\in G_{\mathcal A}}\mathcal A_g,
 \qquad \mathcal A_g\otimes\mathcal A_h\subset\mathcal A_{gh}.
\end{equation}
For a Brauer--Wall reference category $\mathcal C$,
$\rho_{\mathcal C}$ is an isomorphism onto $\BW(k)$; we index its
components by $r\in\BW(k)$ and write them as $\mathcal C_r$.
\end{proposition}
\begin{proof}
Homogeneous Schur's lemma and semisimplicity give
\begin{equation}
\label{eq:component-orthogonality}
 \mathcal A(X,Y)=0
 \quad(X\in\mathcal A_g,\ Y\in\mathcal A_h,\ g\ne h)
\end{equation}
and the direct-sum decomposition into components.
For simple $X,Y$ of types $r,s$,
\[
 \End(X\otimes Y)\cong D_r\gt D_s\cong M_{u|v}(D_{r+s}).
\]
The projections onto different components are supercentral
idempotents, so exactly one is nonzero. Every simple summand of
$X\otimes Y$ has endomorphism superalgebra isomorphic to $D_{r+s}$
by Lemma~\ref{lem:matrix-corners}. This proves
well-definedness of the product and
\begin{equation}
\label{eq:intrinsic-type-additivity}
 \rho_{\mathcal A}(gh)=\rho_{\mathcal A}(g)+\rho_{\mathcal A}(h).
\end{equation}
Associativity and the identity follow from the associator and
unitors. The nonzero evaluation and coevaluation split, so $\one$
is a summand of $X^*\otimes X$ and $X\otimes X^*$; this gives inverses.
Finiteness follows in the fusion case. For a reference category,
the defining existence and uniqueness of types make
$\rho_{\mathcal C}$ bijective.
\end{proof}

\begin{proposition}
\label{prop:comparison-lines}
(1) Components of weakly pointed semisimple tensor $\Pi$-supercategories with
split unit and the same Brauer--Wall type are superequivalent.

(2) Let $n\geq1$ and let $\mathcal B_1,\ldots,\mathcal B_n,\mathcal D$ be
components of weakly pointed semisimple tensor $\Pi$-supercategories with split
unit.
Choose a simple object $S_i\in\mathcal B_i$ for each $i$, and put
$\mathcal B=\mathcal B_1\boxtimes\cdots\boxtimes\mathcal B_n$ and
$U=S_1\boxtimes\cdots\boxtimes S_n\in\mathcal B$,
with both products parenthesized to the left.
Suppose that the superfunctors
\begin{equation}
\label{eq:comparison-functors}
 F,F':\mathcal B\longrightarrow\mathcal D
\end{equation}
induce superalgebra isomorphisms
\begin{equation}
\label{eq:comparison-endomorphisms}
\begin{aligned}
 F_{U,U}:\End_{\mathcal B}(U)
   &\xrightarrow{\sim}\End_{\mathcal D}(F(U)),
   &a&\longmapsto F(a),\\
 F'_{U,U}:\End_{\mathcal B}(U)
   &\xrightarrow{\sim}\End_{\mathcal D}(F'(U)),
   &a&\longmapsto F'(a).
\end{aligned}
\end{equation}
Then $\Nat(F,F')$ is one-dimensional, every nonzero element is a
supernatural isomorphism, and a transformation is determined by its
value at $U$.
\end{proposition}

\begin{proof}
(1) For a simple $S$ in a component $\mathcal E$, put $D=\End(S)$.
Every object is an even sum of $S$ and $\Pi S$, and every graded
right $D$-module is graded free. Matrix entries therefore give a
superequivalence
\begin{equation}
\label{eq:component-module-equivalence}
 \Hom(S,-):\mathcal E\longrightarrow
 \mathsf{Mod}^{\mathrm{gr}}_{\mathrm{fd}}(D),\qquad
 X\longmapsto\mathcal E(S,X),\quad f\longmapsto(u\mapsto fu).
\end{equation}
Components of the same type have isomorphic $D$ and hence are
superequivalent.

(2) Put $A=\End(U)$, choose simple $T\in\mathcal D$, and set
$D=\End(T)$, $V_0=\mathcal D(T,F(U))$, $W_0=\mathcal D(T,F'(U))$.
The hypotheses identify $A$ with $\End_D(V_0)$ and $\End_D(W_0)$,
so $V_0,W_0$ are Morita equivalences from $A$ to $D$.
A degree-$d$ intertwiner is represented by $u:F(U)\to F'(U)$ with
$uF(a)=(-1)^{d|a|}F'(a)u$.
For homogeneous splittings $\sum_\alpha i_\alpha p_\alpha=\id_V$
through $U$, with $|i_\alpha|=|p_\alpha|=e_\alpha$, extend it by
\begin{equation}
\label{eq:comparison-extension}
 \theta_V=\sum_\alpha(-1)^{de_\alpha}
 F'(i_\alpha)uF(p_\alpha).
\end{equation}
The intertwining equation gives supernaturality and independence
of the splittings. Conversely, supernaturality forces this formula.
Thus evaluation gives
\begin{equation}
\label{eq:comparison-evaluation}
 \Nat(F,F')\xrightarrow{\sim}\Hom_{A-D}(V_0,W_0),\qquad
 \theta\longmapsto(v\mapsto\theta_Uv).
\end{equation}
Proposition~\ref{prop:morita-equivalences} proves the assertions;
extending $u^{-1}$ gives the inverse transformation.
\end{proof}

\begin{lemma}
\label{lem:component-tensor-homs}
Let $\mathcal A$ be a weakly pointed semisimple tensor $\Pi$-supercategory with split
unit, and let $X,X'\in\mathcal A_g$ and $Y,Y'\in\mathcal A_h$. Then
tensoring gives an isomorphism of super vector spaces
\begin{equation}
\label{eq:component-tensor-homs}
 \mathcal A(X,X')\otimes_{k}\mathcal A(Y,Y')
 \xrightarrow{\sim}\mathcal A(X\otimes Y,X'\otimes Y').
\end{equation}
\end{lemma}
\begin{proof}
For chosen simple representatives in the two components this is
Proposition~\ref{prop:simple-tensor}. Homogeneous isomorphisms
reduce arbitrary simples to these representatives, with the
invertible signs of super interchange. Even direct-sum
decompositions then reduce the general case to these matrix entries.
\end{proof}

\subsection{Classifying parameters and construction}
\label{subsec:twists-pullback}

We introduce group cochain parameters $t$. Assuming the existence of a reference category $\mathcal C$ over $k$, we construct a weakly pointed fusion $\Pi$-supercategory $\mathcal P_\mathcal{C}(t)$, which will be a representative of each monoidal superequivalence class.

Definitions of cochains and cup products are collected in Appendix~\ref{app:cochains}.

\begin{definition}
\label{def:twist}
A \emph{twist} is a quadruple $t=(G,\rho,\nu,\omega)$, where $G$ is a
finite group, $\rho\in Z^1(G,\BW(k))$,
$\nu\in Z^2(G,\Z/2\Z)$, and $\omega\in C^3(G,k^\times)$, such that
\begin{equation}
  \label{eq:twist-cocycle}
  \delta\omega=(-1)^{\nu\smile\nu}.
\end{equation}
We write $\Tw_k$ for the class of all twists.
\end{definition}

\begin{definition}
\label{def:gauge}
Let $s=(G,\rho,\nu,\omega)$ and $t=(H,\rho',\nu',\omega')$ be twists.
A \emph{gauge} from $s$ to $t$ is a triple $(\phi,\lambda,\beta)$,
where $\phi:G\to H$ is a group isomorphism,
$\lambda\in C^1(G,\Z/2\Z)$, and $\beta\in C^2(G,k^\times)$, satisfying
\begin{equation}
\label{eq:gauge-degree}
  \phi^*\rho'=\rho,\qquad
  \phi^*\nu'=\nu+\delta\lambda,
\end{equation}
\begin{equation}
\label{eq:gauge-associator}
  \phi^*\omega'
  =\omega(\delta\beta)^{-1}
    (-1)^{(\phi^*\nu')\smile\lambda+\lambda\smile\nu}.
\end{equation}
The twists $s$ and $t$ are \emph{gauge equivalent}, written $s\sim t$,
if there exists a gauge from $s$ to $t$.
\end{definition}

\begin{proposition}
\label{prop:gauge-equivalence}
Gauge equivalence is an equivalence relation on $\Tw_k$: $(\id_G,0,1)$ is a
gauge from $s$ to $s$, the composite of gauges $(\phi,\lambda,\beta)$ from
$s$ to $t$ and $(\psi,\mu,\eta)$ from $t$ to $u$ is a gauge from $s$ to $u$
with group isomorphism $\psi\phi$, and the inverse of $(\phi,\lambda,\beta)$
is a gauge from $t$ to $s$ with group isomorphism $\phi^{-1}$.
\end{proposition}

\begin{proof}
The identity gauge is $(\id_G,0,1)$. Composition and inversion are
\begin{gather*}
 (\psi,\mu,\eta)\circ(\phi,\lambda,\beta)
 =\bigl(\psi\phi,\lambda+\phi^*\mu,
   \beta\,\phi^*\eta\,(-1)^{(\phi^*\mu)\smile\lambda}\bigr),\\
 (\phi,\lambda,\beta)^{-1}
 =\bigl(\phi^{-1},(\phi^{-1})^*\lambda,
   (\phi^{-1})^*(\beta^{-1}(-1)^{\lambda\smile\lambda})\bigr).
\end{gather*}
The scalar correction in the composite follows from
$\delta(b\smile\lambda)=\delta b\smile\lambda+b\smile\delta\lambda$
with $b=\phi^*\mu$; substitution gives the inverse formula.
Lemma~\ref{lem:coboundary-identities} ensures preservation of the
twist equations and normalization.
\end{proof}

\begin{definition}
\label{def:parity-maps}
Let $\mathcal C$ be a monoidal $\Pi$-supercategory.
Put $\zeta_0=\id_{\one}$ and $\zeta_1=\zeta$.
For $a,b\in\Z/2\Z$ and an object $X$, define
even isomorphisms $c_{a,X}:P_a\otimes X\to X\otimes P_a$ and
$m_{a,b}:P_a\otimes P_b\to P_{a+b}$ by
\begin{equation}
\label{eq:parity-maps}
\begin{aligned}
 c_{a,X}
 &=(\id_X\otimes\zeta_a^{-1})r_X^{-1}l_X
   (\zeta_a\otimes\id_X),\\
 m_{a,b}
 &=\zeta_{a+b}^{-1}l_{\one}(\zeta_a\otimes\zeta_b).
\end{aligned}
\end{equation}
Let $G$ be a group and $\nu\in C^2(G,\Z/2\Z)$ with $\delta\nu=0$.
For $g,h,\ell\in G$ and objects $X,Y,Z$, define an even isomorphism
\begin{equation}
\label{eq:parity-associator-components}
\begin{aligned}
 A^\nu_{g,h,\ell}(X,Y,Z):{}\;
 &(((X\otimes Y)\otimes P_{\nu(g,h)})\otimes Z)\otimes P_{\nu(gh,\ell)}\\
 &\longrightarrow
   (X\otimes((Y\otimes Z)\otimes P_{\nu(h,\ell)}))\otimes P_{\nu(g,h\ell)}
\end{aligned}
\end{equation}
by
\begin{equation}
\label{eq:parity-associator}
\begin{aligned}
 A^\nu_{g,h,\ell}(X,Y,Z)
  ={}&\bigl(\id_{X\otimes(Y\otimes Z)}\otimes
      m_{\nu(h,\ell),\nu(g,h\ell)}^{-1}m_{\nu(g,h),\nu(gh,\ell)}\bigr)\\
     &\circ\bigl(((\id_X\otimes\id_Y)\otimes c_{\nu(g,h),Z})
      \otimes\id_{P_{\nu(gh,\ell)}}\bigr),
\end{aligned}
\end{equation}
with associators suppressed.
Finally, define even isomorphisms
\begin{equation}
\label{eq:parity-unitors}
\begin{aligned}
 \tilde l_X&=l_Xr_{\one\otimes X}:
   (\one\otimes X)\otimes\one\longrightarrow X,\\
 \tilde r_X&=r_Xr_{X\otimes\one}:
   (X\otimes\one)\otimes\one\longrightarrow X.
\end{aligned}
\end{equation}
\end{definition}

\begin{definition}
\label{def:pullback-construction}
Let $\mathcal C$ be a Brauer--Wall reference category. 
For each finite family $(X_i)_{i\in I}$ of objects of $\mathcal C$, 
fix a direct sum
$\bigoplus_{i\in I}X_i$ with even inclusions $\iota_i$ and
projections $p_i$.
If $X_i=0$ for all $i\ne j$, we choose $X_j$ itself as
the direct sum, with $\iota_j=p_j=\id_{X_j}$.
In the formulas below, we identify
$\bigoplus_i(\bigoplus_j X_{i,j})$ with
$\bigoplus_{(i,j)}X_{i,j}$, and
$(\bigoplus_iX_i)\otimes(\bigoplus_jY_j)$ with
$\bigoplus_{(i,j)}(X_i\otimes Y_j)$,
using the isomorphisms induced by the inclusions and projections.

For $t=(G,\rho,\nu,\omega)\in\Tw_k$, we define a monoidal $\Pi$-supercategory
$\mathcal P_{\mathcal C}(t)$, a \emph{pullback}, as follows.
Its objects are families
$X=(X_g)_{g\in G}$ with $X_g\in\mathcal C_{\rho(g)}$.
Its morphism spaces are
\begin{equation}
\label{eq:pullback-hom}
 \mathcal P_{\mathcal C}(t)(X,Y)
  =\bigoplus_{g\in G}\mathcal C(X_g,Y_g).
\end{equation}
Identities and composition are defined componentwise.
The tensor product is
\begin{equation}
\label{eq:pullback-tensor}
\begin{aligned}
 (X\otimes_tY)_\ell
  &=\bigoplus_{gh=\ell}
    \bigl((X_g\otimes Y_h)\otimes P_{\nu(g,h)}\bigr),\\
 (f\otimes_t f')_\ell
  &=\bigoplus_{gh=\ell}
    \bigl((f_g\otimes f'_h)\otimes\id_{P_{\nu(g,h)}}\bigr).
\end{aligned}
\end{equation}
The unit object $\one_t$, parity object $\pi_t$, and odd isomorphism
$\zeta_t:\pi_t\to\one_t$ are supported at $1\in G$, with 
components $\one,\pi,\zeta$, respectively.
The associator
$\alpha^t_{X,Y,Z}:(X\otimes_tY)\otimes_tZ
\to X\otimes_t(Y\otimes_tZ)$
is
\begin{equation}
\label{eq:pullback-associator}
 (\alpha^t_{X,Y,Z})_r
  =\bigoplus_{gh\ell=r}
    \omega(g,h,\ell)\,A^\nu_{g,h,\ell}(X_g,Y_h,Z_\ell).
\end{equation}
The unitors $l^t_X:\one_t\otimes_tX\to X$ and
$r^t_X:X\otimes_t\one_t\to X$ are
\begin{equation}
\label{eq:pullback-unitors}
 (l^t_X)_g=\tilde l_{X_g},\qquad
 (r^t_X)_g=\tilde r_{X_g}.
\end{equation}
The pentagon and triangle axioms follow from
Proposition~\ref{prop:parity-pentagon} below.
\end{definition}

\begin{proposition}
\label{prop:parity-pentagon}
Let $\mathcal C$ be a Brauer--Wall reference category and let $t\in\Tw_k$.
The associator $\alpha^t$ and unitors $l^t,r^t$ defined in
Definition~\ref{def:pullback-construction} are even supernatural
isomorphisms satisfying the pentagon and triangle identities.
\end{proposition}

\begin{proof}
Substitution in \eqref{eq:parity-maps}, using super interchange,
gives
\begin{equation}
\label{eq:parity-unit-cases}
 c_{0,X}=r_X^{-1}l_X,\qquad m_{0,b}=l_{P_b},\qquad m_{b,0}=r_{P_b},
\end{equation}
and, with constraints suppressed,
\begin{equation}
\label{eq:parity-crossing-contraction}
 c_{a+b,X}(m_{a,b}\otimes\id_X)
 =(\id_X\otimes m_{a,b})(c_{a,X}\otimes\id_{P_b})
 (\id_{P_a}\otimes c_{b,X}).
\end{equation}
The maps $m_{a,b}$ combine parity factors associatively, and the
maps $c_{a,X}$ interchange $P_a$ with $X$ compatibly with tensor
products. Interchanging $P_a$ and $P_b$ gives the sign
\[
 c_{a,P_b}=(-1)^{ab}m_{b,a}^{-1}m_{a,b}.
\]
For the pentagon of $A^\nu$ at $(g,h,\ell,m)$, let $T_3$ and $T_2$
denote the three-edge and two-edge composites. Using $m_{a,b}$ to
combine the parity factors in both composites gives the same map,
up to the sign from interchanging $P_{\nu(g,h)}$ and
$P_{\nu(\ell,m)}$ using $c$. Thus
\begin{equation}
\label{eq:parity-pentagon-defect}
 T_3=(-1)^{\nu(g,h)\nu(\ell,m)}T_2.
\end{equation}
The scalar ratio for $\alpha^t=\omega A^\nu$ is
$\delta\omega=(-1)^{\nu\smile\nu}$, which cancels this sign.
Normalization and \eqref{eq:parity-unit-cases} give
\begin{equation}
\label{eq:parity-triangle}
 \bigl((\id_X\otimes\tilde l_Y)\otimes\id_{P_{\nu(g,h)}}\bigr)
 A^\nu_{g,1,h}(X,\one,Y)
 =(\tilde r_X\otimes\id_Y)\otimes\id_{P_{\nu(g,h)}}.
\end{equation}
All maps are even supernatural by their defining formulas.
The identities extend to families by direct sums.
\end{proof}

\begin{proposition}
\label{prop:pullback-laws}
Let $\mathcal C$ be a Brauer--Wall reference category and
$t=(G,\rho,\nu,\omega)$ a twist. Then $\mathcal P_{\mathcal C}(t)\in\W_k$.
Choose simple $S_r\in\mathcal C_r$ with $S_0=\one$, and let $X_g$ be
$S_{\rho(g)}$ supported at $g$. Every simple object is isomorphic to a
unique $X_g$. The map $g\mapsto[X_g]$ is a group isomorphism
\begin{equation}
\label{eq:pullback-group}
 G\xrightarrow{\sim}G_{\mathcal P_{\mathcal C}(t)},
 \qquad \rho_{\mathcal P_{\mathcal C}(t)}([X_g])=\rho(g).
\end{equation}
\end{proposition}

\begin{proof}
Write $\mathcal P=\mathcal P_{\mathcal C}(t)$ and $S[g]$ for a
family supported at $g$. Its even category is
$\prod_{g\in G}\underline{\mathcal C_{\rho(g)}}$, so it is finite
semisimple with split unit. Its homogeneous simple classes are
the $X_g$, and $\End(X_g)\cong D_{\rho(g)}$.

For rigidity, let $X=S[g]$ be simple and put
\[
 a=\nu(g,g^{-1})=\nu(g^{-1},g),\qquad
 Y=(S^*\otimes P_a)[g^{-1}].
\]
Suppressing constraints and the supported-summand inclusions,
take
\begin{align*}
 \eta&=(\id_{S\otimes S^*}\otimes m_{a,a}^{-1})\coev_S,\\
 \vartheta&=\omega(g,g^{-1},g)^{-1}
 (\ev_S\otimes m_{a,a})
 ((\id_{S^*}\otimes c_{a,S})\otimes\id_{P_a}).
\end{align*}
Equation~\eqref{eq:parity-crossing-contraction} cancels the parity
factors in the first snake, while the scalar cancels the associator
factor. The second snake is then a nonzero even idempotent on the
simple object $Y$, hence is $\id_Y$.
Left duality is injective on the finite set of even simple classes,
so right duals also exist. For even direct sums, with inclusions
$i_i,j_i$ and projections $p_i,q_i$ on an object and its dual, use
\begin{equation}
\label{eq:direct-sum-duals}
 \ev_V=\sum_i\vartheta_i(q_i\otimes_t p_i),\qquad
 \coev_V=\sum_i(i_i\otimes_t j_i)\eta_i.
\end{equation}
The product of a simple with its dual lies in the neutral component,
which consists of sums of $\one_t$ and $\pi_t$, proving weak
invertibility. Finally, $X_g\otimes_tX_h$ is nonzero and supported
at $gh$, giving the group and type assertions.
\end{proof}

\subsection{Classification}
\label{subsec:reconstruction-classification}

Comparing the components of $\mathcal A\in\W_k$ with those of a fixed
reference category, we extract a twist. We prove that gauge classes of
twists classify the categories in $\W_k$ up to monoidal superequivalence.

\begin{definition}
\label{def:reconstruction-comparisons}
Let $\mathcal A\in\W_k$ and let $\mathcal C$ be a Brauer--Wall reference
category. Put $G=G_{\mathcal A}$ and $\rho=\rho_{\mathcal A}$.
\emph{Normalized comparison data} for $\mathcal A$ consist of a pair
$\Cmp=(F,\theta)$ of a family $F=(F_g)_{g\in G}$ of superequivalences
and a family $\theta=(\theta_{g,h})_{g,h\in G}$ of homogeneous
supernatural isomorphisms,
\begin{equation}
\label{eq:binary-comparison}
\begin{aligned}
 F_g&:\mathcal A_g\longrightarrow\mathcal C_{\rho(g)}
 &&(g\in G),\\
 \theta_{g,h}(X,Y)&:F_{gh}(X\otimes Y)\longrightarrow F_gX\otimes F_hY
 &&(g,h\in G,\ X\in\mathcal A_g,\ Y\in\mathcal A_h),
\end{aligned}
\end{equation}
where $\theta_{g,h}$ is supernatural between the superfunctors
$(X,Y)\mapsto F_{gh}(X\otimes Y)$ and $(X,Y)\mapsto F_gX\otimes F_hY$
from $\mathcal A_g\boxtimes\mathcal A_h$ to $\mathcal C_{\rho(gh)}$,
such that $F_1\one=\one$, $F_1\pi=\pi$, $F_1\zeta=\zeta$, and
\begin{equation}
\label{eq:unit-comparison}
 \theta_{1,h}(\one,Y)=l_{F_hY}^{-1}F_h(l_Y),\qquad
 \theta_{g,1}(X,\one)=r_{F_gX}^{-1}F_g(r_X).
\end{equation}
\end{definition}

\begin{proposition}
\label{prop:reconstruction}
Normalized comparison data exist for every $\mathcal A\in\W_k$ and every
Brauer--Wall reference category $\mathcal C$.
Let $\Cmp=(F,\theta)$ be normalized comparison data for $\mathcal A$.
Then $\nu(g,h)=|\theta_{g,h}|$ defines a normalized cocycle
$\nu\in C^2(G,\Z/2\Z)$.
For $V\in\mathcal C$, put $e_{a,V}=r_V(\id_V\otimes\zeta_a)$,
a homogeneous isomorphism of degree $a$.
For $X\in\mathcal A_g$ and $Y\in\mathcal A_h$, define the even isomorphism
\begin{equation}
\label{eq:even-binary-comparison}
\begin{aligned}
 K_{g,h}(X,Y)
 ={}&e_{\nu(g,h),\,F_gX\otimes F_hY}^{-1}\theta_{g,h}(X,Y)\\
 &\colon F_{gh}(X\otimes Y)
 \longrightarrow(F_gX\otimes F_hY)\otimes P_{\nu(g,h)}.
\end{aligned}
\end{equation}
For $X\in\mathcal A_g$, $Y\in\mathcal A_h$, $Z\in\mathcal A_\ell$,
define even isomorphisms
\begin{equation}
\label{eq:triple-comparisons}
\begin{aligned}
 \mathcal L_{X,Y,Z}
 ={}&A^\nu_{g,h,\ell}(F_gX,F_hY,F_\ell Z)\\
 &\circ\bigl((K_{g,h}(X,Y)\otimes\id_{F_\ell Z})\otimes\id_{P_{\nu(gh,\ell)}}\bigr)
   K_{gh,\ell}(X\otimes Y,Z),\\
 \mathcal R_{X,Y,Z}
 ={}&\bigl((\id_{F_gX}\otimes K_{h,\ell}(Y,Z))\otimes\id_{P_{\nu(g,h\ell)}}\bigr)\\
 &\circ K_{g,h\ell}(X,Y\otimes Z)F_{gh\ell}(\alpha^{\mathcal A}_{X,Y,Z})\\
 &\colon F_{gh\ell}((X\otimes Y)\otimes Z)\\
 &\longrightarrow
 (F_gX\otimes((F_hY\otimes F_\ell Z)\otimes P_{\nu(h,\ell)}))\otimes P_{\nu(g,h\ell)}.
\end{aligned}
\end{equation}
There is a unique map $\omega:G^3\to k^\times$ such that
\begin{equation}
\label{eq:extracted-omega}
 \mathcal R_{X,Y,Z}=\omega(g,h,\ell)\mathcal L_{X,Y,Z}
\end{equation}
for all $X,Y,Z$ in these components.
Then $t=(G,\rho,\nu,\omega)$ is a twist.
\end{proposition}

\begin{proof}
Choose component equivalences using
\eqref{eq:component-module-equivalence}, normalized at $\one,\pi,\zeta$.
Proposition~\ref{prop:comparison-lines} supplies binary comparison
lines; their unit values can be fixed by \eqref{eq:unit-comparison}.
The two triple composites of $\theta$ lie in the same super line,
so their degrees agree:
\[
 \nu(g,h)+\nu(gh,\ell)=\nu(h,\ell)+\nu(g,h\ell).
\]
After inserting the parity factors, the comparisons $K$ are even.
The same line argument gives the unique scalar
$\mathcal R=\omega(g,h,\ell)\mathcal L$.
Transporting the associator of $\mathcal A$ through $K$ gives
$\omega A^\nu$; its pentagon and
\eqref{eq:parity-pentagon-defect} imply
\[
 \delta\omega(g,h,\ell,m)=(-1)^{\nu(g,h)\nu(\ell,m)}.
\]
The unit comparisons give $\nu(1,g)=\nu(g,1)=0$ and
$\omega(g,1,h)=1$. The pentagon at $(1,1,g,h)$ and $(g,h,1,1)$
then gives $\omega(1,g,h)=\omega(g,h,1)=1$.
\end{proof}

\begin{proposition}
\label{prop:reconstruction-functor}
Let $\mathcal A\in\W_k$ and let $\mathcal C$ be a Brauer--Wall reference category.
Let $\Cmp=(F,\theta)$ be normalized comparison data for $\mathcal A$, and
let $K_{g,h}$ and $t$ be as in Proposition~\ref{prop:reconstruction}.
The following data define a monoidal superequivalence
$(F,J,\varphi):\mathcal A\to\mathcal P_{\mathcal C}(t)$.
For every object $X$ of $\mathcal A$ choose a decomposition
$X\cong\bigoplus_g X_g$ into components, with even inclusions $i_g^X$
and projections $p_g^X$.
The superfunctor $F$ is
\begin{equation}
\label{eq:reconstruction-functor}
 (FX)_g=F_gX_g,\qquad
 (Ff)_g=F_g(p_g^Yfi_g^X)\quad(f:X\to Y),
\end{equation}
and the tensorator $J_{X,Y}:FX\otimes_{t}FY\to F(X\otimes Y)$ and the
unit map $\varphi:\one_{t}\to F\one$ are
\begin{equation}
\label{eq:reconstruction-tensorator}
\begin{aligned}
 (J_{X,Y})_\ell
 &=\sum_{gh=\ell}
   F_\ell\bigl(p_\ell^{X\otimes Y}(i_g^X\otimes i_h^Y)\bigr)\,
   K_{g,h}(X_g,Y_h)^{-1}\,p_{(g,h)},\\
 \varphi_1&=F_1(p_1^{\one}),\qquad
 \varphi_g=0\quad(g\neq1),
\end{aligned}
\end{equation}
where $p_{(g,h)}$ is the projection of the direct sum
$(FX\otimes_{t}FY)_\ell=\bigoplus_{gh=\ell}(F_gX_g\otimes F_hY_h)\otimes P_{\nu(g,h)}$
onto its $(g,h)$ summand, as fixed in Definition~\ref{def:pullback-construction}.
\end{proposition}

\begin{proof}
Component orthogonality makes $F$ a superfunctor, fully faithful
and essentially surjective because each $F_g$ is so.
The tensorator is the direct sum of $K_{g,h}^{-1}$ under the
chosen component decompositions, hence is even, supernatural,
and invertible. The required equations are
\begin{equation}
\label{eq:reconstruction-monoidal-coherence}
\begin{aligned}
 F(\alpha^{\mathcal A}_{X,Y,Z})J_{X\otimes Y,Z}
 (J_{X,Y}\otimes_t\id_{FZ})
 &=J_{X,Y\otimes Z}(\id_{FX}\otimes_tJ_{Y,Z})\alpha^t_{FX,FY,FZ},\\
 F(l_X)J_{\one,X}(\varphi\otimes_t\id_{FX})&=l^t_{FX},\\
 F(r_X)J_{X,\one}(\id_{FX}\otimes_t\varphi)&=r^t_{FX}.
\end{aligned}
\end{equation}
On homogeneous summands these are respectively
\eqref{eq:extracted-omega} and \eqref{eq:unit-comparison}; even
direct sums give the general case.
\end{proof}

\begin{proposition}
\label{prop:gauge-functors}
Let $s=(G,\rho,\nu,\omega)$ and $t=(H,\rho',\nu',\omega')$ be twists,
and let $\mathcal C$ be a Brauer--Wall reference category. Then
\begin{equation}
\label{eq:gauge-monoidal-equivalence}
 \mathcal P_{\mathcal C}(s)\simeqm\mathcal P_{\mathcal C}(t)
 \quad\Longleftrightarrow\quad s\sim t.
\end{equation}
More precisely, a gauge $(\phi,\lambda,\beta)$ from $s$ to $t$ defines a
monoidal superequivalence
$(F,J,\varphi):\mathcal P_{\mathcal C}(s)\to\mathcal P_{\mathcal C}(t)$ by
\begin{equation}
\label{eq:gauge-functor}
 (FX)_{\phi(g)}=X_g\otimes P_{\lambda(g)},\qquad
 (Ff)_{\phi(g)}=f_g\otimes\id_{P_{\lambda(g)}},
\end{equation}
\begin{equation}
\label{eq:gauge-tensorator}
\begin{aligned}
 (J_{X,Y})_{\phi(\ell)}
 &=\sum_{gh=\ell}\beta(g,h)\,
   (\iota_{(g,h)}\otimes\id_{P_{\lambda(\ell)}})\,
   \Gamma_{g,h}(X_g,Y_h)\,p_{(\phi(g),\phi(h))},\\
 \varphi_1&=r_{\one}^{-1},\qquad
 \varphi_h=0\quad(h\neq1).
\end{aligned}
\end{equation}
Here $p_{(\phi(g),\phi(h))}$ is the projection of the fixed direct sum
$(FX\otimes_tFY)_{\phi(\ell)}$ onto its $(\phi(g),\phi(h))$ summand,
$\iota_{(g,h)}$ is the inclusion of the $(g,h)$ summand of
$(X\otimes_sY)_\ell$, and $\Gamma_{g,h}(X_g,Y_h)$ is the even isomorphism
\begin{equation}
\label{eq:gauge-tensorator-summand}
\begin{aligned}
 \Gamma_{g,h}(X_g,Y_h)
 ={}&\bigl(\id_{X_g\otimes Y_h}\otimes
   m_{\nu(g,h),\lambda(gh)}^{-1}\,
   m_{\lambda(g)+\lambda(h),\nu'(\phi(g),\phi(h))}\,
   (m_{\lambda(g),\lambda(h)}\otimes\id_{P_{\nu'(\phi(g),\phi(h))}})\bigr)\\
 &\circ\bigl(((\id_{X_g}\otimes c_{\lambda(g),Y_h})\otimes\id_{P_{\lambda(h)}})
   \otimes\id_{P_{\nu'(\phi(g),\phi(h))}}\bigr)\\
 &\colon ((X_g\otimes P_{\lambda(g)})\otimes(Y_h\otimes P_{\lambda(h)}))
   \otimes P_{\nu'(\phi(g),\phi(h))}\\
 &\longrightarrow
 ((X_g\otimes Y_h)\otimes P_{\nu(g,h)})\otimes P_{\lambda(gh)}
\end{aligned}
\end{equation}
with associators suppressed. The maps $m_{a,b}$ and their inverses
in this formula are composable by \eqref{eq:gauge-degree}.
Conversely, every monoidal superequivalence $E$ from
$\mathcal P_{\mathcal C}(s)$ to $\mathcal P_{\mathcal C}(t)$ induces a gauge $(\phi,\lambda,\beta)$ from
$s$ to $t$, where $\phi:G\to H$ is the isomorphism induced by $E$ on the
groups \eqref{eq:pullback-group} of simple objects.
In particular, if $G=H$ and $E(X_g)\cong X_g$ for every $g\in G$, with
$X_g$ as in Proposition~\ref{prop:pullback-laws}, then
$(\id_G,\lambda,\beta)$ is a gauge from $s$ to $t$.
\end{proposition}

\begin{proof}
Write $\bar\nu=\phi^*\nu'$ and $\bar\omega=\phi^*\omega'$.
Reindexing by $\phi$ and tensoring with $P_{\lambda(g)}$ make $F$
a superequivalence. The degree equation makes $\Gamma$ well defined.
To check its scalar coherence, use $e_{a,V}=r_V(\id_V\otimes\zeta_a)$
to delete parity factors. For three objects $U_1,U_2,U_3$, set
\begin{equation}
\label{eq:gauge-ordered-deletions}
\begin{aligned}
 D_L^{U_1,U_2,U_3}(u,v)
 &=\alpha_{U_1,U_2,U_3}e_{v,(U_1\otimes U_2)\otimes U_3}
 ((e_{u,U_1\otimes U_2}\otimes\id_{U_3})\otimes\id_{P_v}),\\
 D_R^{U_1,U_2,U_3}(u,v)
 &=e_{v,U_1\otimes(U_2\otimes U_3)}
 ((\id_{U_1}\otimes e_{u,U_2\otimes U_3})\otimes\id_{P_v}).
\end{aligned}
\end{equation}
Put $(n_1,n_2,n_3,n_4)=(\nu(g,h),\nu(gh,\ell),\nu(h,\ell),\nu(g,h\ell))$.
Super interchange gives
\begin{equation}
\label{eq:gauge-associator-deletion}
 D_R(n_3,n_4)A^\nu_{g,h,\ell}
 =(-1)^{n_1n_2+n_3n_4}D_L(n_1,n_2).
\end{equation}
Applying the same deletion to the two binary tensorators makes
both coherence routes scalar multiples of a common isomorphism.
Their relative sign is
$(-1)^{\bar\nu(g,h)\lambda(\ell)+\lambda(g)\nu(h,\ell)}$.
Thus coherence is precisely
\[
 \bar\omega=\omega(\delta\beta)^{-1}
 (-1)^{\bar\nu\smile\lambda+\lambda\smile\nu}.
\]
Normalization gives the unit equations.

Conversely, a monoidal superequivalence $E$ induces a group
isomorphism $\phi$ preserving types. Comparison lines give
$\eta_g:E_g\Rightarrow\id$ of degree $\lambda(g)$, with
$\lambda(1)=0$. Inserting parity factors turns these into an even
comparison with the displayed $F$. Transport the tensorator of
$E$ to $F$. Since the endomorphism line of the binary tensor
functor is $k\id$, its degree and scalar parts give
\[
 \bar\nu=\nu+\delta\lambda,\qquad
 \widehat J_{g,h}=\beta(g,h)\Gamma_{g,h}.
\]
The preceding coherence calculation gives the scalar gauge
equation. Its unit instances imply
$\beta(1,g)=\beta(g,1)=\beta(1,1)$, so dividing by $\beta(1,1)$
normalizes $\beta$ without changing its coboundary.
The final assertion follows from the induced map on simple classes.
\end{proof}

\begin{theorem}
\label{thm:classification}
Let $\mathcal C$ be a Brauer--Wall reference category. The map
$t\mapsto\mathcal P_{\mathcal C}(t)$ induces a bijection
\begin{equation}
\label{eq:classification}
 \Tw_k/{\sim}\xrightarrow{\ \sim\ }\W_k/{\simeqm}.
\end{equation}
\end{theorem}

\begin{proof}
Proposition~\ref{prop:pullback-laws} makes the map well defined,
and Proposition~\ref{prop:gauge-functors} gives injectivity on
equivalence classes. Propositions~\ref{prop:reconstruction}
and~\ref{prop:reconstruction-functor} give surjectivity.
\end{proof}

\section{Stacking}
\label{sec:stacking}

Given a finite group $G$, we consider a product of $G$-graded categories in $\mathcal W_k$.
After defining the product, we show that the equivalence classes with $G$-grading form an abelian group.
Then, we represent the multiplication in terms of twists.

\subsection{Completions and completed products}
\label{subsec:completions-products}

We prepare additive and idempotent completion for supercategories and show that $\mathcal W_k$ is closed under the product $\boxtimes$ after the completion. 

\begin{definition}
\label{def:completion-construction}
Let $\mathcal E$ be a supercategory.
An even idempotent $e:X\to X$ is \emph{split} if there are an object $Y$
and even morphisms $i:Y\to X$ and $p:X\to Y$ such that
$i\circ p=e$ and $p\circ i=\id_Y$.
The supercategory $\mathcal E$ has \emph{split even idempotents} if every
even idempotent in $\mathcal E$ is split.
The \emph{additive envelope} $\operatorname{Add}(\mathcal E)$ is the
supercategory whose objects are finite formal direct sums of objects
of $\mathcal E$, and whose morphisms are matrices of morphisms in
$\mathcal E$. A matrix has degree $d$ if all its entries have degree $d$.
The \emph{even idempotent completion} $\operatorname{Kar}_0(\mathcal E)$
is the supercategory whose objects are the pairs $(X,e)$ such that
$e\in\End_{\mathcal E}(X)_0$ and $e^2=e$, and whose morphism spaces are
\begin{equation}
\label{eq:even-idempotent-completion}
 \operatorname{Kar}_0(\mathcal E)((X,e),(Y,f))
 =f\mathcal E(X,Y)e,
\end{equation}
with the grading inherited from $\mathcal E(X,Y)$.
\end{definition}

\begin{definition}
\label{def:completions}
The \emph{completion} of a supercategory $\mathcal E$ is
\begin{equation}
\label{eq:supercategory-completion}
 \overline{\mathcal E}
 =\operatorname{Kar}_0(\operatorname{Add}(\mathcal E)).
\end{equation}
Define a superfunctor
\begin{equation}
\label{eq:completion-inclusion}
 \iota:\mathcal E\longrightarrow\overline{\mathcal E},
 \qquad X\longmapsto((X),\id_X),\quad f\longmapsto(f).
\end{equation}
In $\overline{\mathcal E}$, we write $X$ for $\iota(X)$
and $f$ for $\iota(f)$.
\end{definition}

\begin{proposition}
\label{prop:completion-construction}
For every supercategory $\mathcal E$, its completion
$\overline{\mathcal E}$ is additive and has split even idempotents.
The superfunctor $\iota$ is fully faithful, and every object
of $\overline{\mathcal E}$ is an even retract of a finite direct sum
of objects of $\mathcal E$.
\end{proposition}
\begin{proof}
Direct sums are given by concatenation and block idempotents:
\[
 (X,e)\oplus(Y,f)=(X\oplus Y,e\oplus f).
\]
An even idempotent $u$ on $(X,e)$ splits through $(X,u)$, with
both maps equal to $u$. The morphism-space formula
\eqref{eq:even-idempotent-completion} gives full faithfulness of
$\iota$, and the two maps $e$ exhibit $(X,e)$ as a retract of
$(X,\id_X)=\bigoplus_i\iota X_i$.
\end{proof}

\begin{proposition}
\label{prop:completion-extension}
Let $\mathcal E$ be a supercategory and let $\mathcal D$ be an additive
supercategory with split even idempotents.
Every superfunctor $F:\mathcal E\to\mathcal D$ admits an extension
\begin{equation}
\label{eq:completed-product-extension}
 \overline F:\overline{\mathcal E}\longrightarrow\mathcal D,
 \qquad \overline F\circ\iota=F.
\end{equation}
For any two superfunctors $G,G':\overline{\mathcal E}\to\mathcal D$,
every supernatural transformation $\theta:G\circ\iota\Rightarrow G'\circ\iota$
extends uniquely to a supernatural transformation
$\overline\theta:G\Rightarrow G'$ with $\overline\theta_X=\theta_X$
for all $X\in\mathcal E$, yielding an isomorphism of super vector spaces
\begin{equation}
\label{eq:completion-natural-transformations}
 \Nat(G\circ\iota,G'\circ\iota)\xrightarrow{\sim}\Nat(G,G'),
 \qquad \theta\longmapsto\overline\theta.
\end{equation}
In particular, for any two extensions $\overline F,\overline F'$ of $F$,
there is a unique even supernatural isomorphism
\begin{equation}
\label{eq:completed-product-extension-uniqueness}
 \theta:\overline F\overset{\sim}{\Longrightarrow}\overline F',
 \qquad \theta_X=\id_{FX}
 \quad(X\in\mathcal E).
\end{equation}
These statements also hold for finitely many source variables, with
each source completed separately.
\end{proposition}
\begin{proof}
First extend $F$ by matrices:
\begin{equation}
\label{eq:additive-functor-extension}
 F^+(X)=\bigoplus_iFX_i,\qquad
 F^+((f_{ji}))=\sum_{i,j}\iota_j^{\mathcal D}F(f_{ji})p_i^{\mathcal D}.
\end{equation}
Choose even splittings $I_eP_e=F^+(e)$ and $P_eI_e=\id_{Q_e}$,
using identity splittings on the original objects, and define
\begin{equation}
\label{eq:completion-extension-construction}
 \overline F(X,e)=Q_e,\qquad
 \overline F(u)=P_fF^+(u)I_e
 \quad\bigl(u:(X,e)\to(Y,f)\bigr).
\end{equation}
The identities $fu=u=ue$ prove functoriality.

For arbitrary superfunctors $G,G':\overline{\mathcal E}\to\mathcal D$
and an object $(X,e)$ of $\overline{\mathcal E}$ with $X=(X_1,\ldots,X_r)$,
let $i_j:\iota X_j\to(X,e)$ and $p_j:(X,e)\to\iota X_j$ be the canonical
even morphisms obtained from the direct-sum inclusions, projections, and
the retraction $e$, so that $\sum_j i_jp_j=\id_{(X,e)}$.
For a homogeneous $\theta\in\Nat(G\circ\iota,G'\circ\iota)$ of degree $d$,
set
\[
 \overline\theta_{(X,e)}=\sum_j G'(i_j)\,\theta_{X_j}\,G(p_j).
\]
For any homogeneous $u:(X,e)\to(Y,f)$ with corresponding maps $j_k,q_k$
on $(Y,f)$, the composites $q_kui_j$ lie in $\mathcal E(X_j,Y_k)$ of
degree $|u|$; supernaturality of $\theta$ on $\mathcal E$ therefore
implies supernaturality of $\overline\theta$. Conversely, for any
$\eta\in\Nat(G,G')$ extending $\theta$, supernaturality with respect to
the even maps $i_j,p_j$ forces
$\eta_{(X,e)}=\sum_j G'(i_j)\eta_{\iota X_j}G(p_j)=\overline\theta_{(X,e)}$.
Since unique extension preserves identities and composition, applying
\eqref{eq:completion-natural-transformations} to $\id_F$ gives
\eqref{eq:completed-product-extension-uniqueness}.
For several variables, extend by matrices separately and split the even
idempotent $F^+(e_1\boxtimes\cdots\boxtimes e_n)$; the tensor products of
the even maps $i_j,p_j$ satisfy the same completeness identity, so the
same formula gives the multivariable version of
\eqref{eq:completion-natural-transformations}.
\end{proof}

\begin{proposition}
\label{prop:monoidal-completion}
Let $\mathcal E$ be a monoidal supercategory.
Choose an extension of
$\iota\circ\otimes:
\mathcal E\boxtimes\mathcal E\to\overline{\mathcal E}$
in both variables:
\begin{equation}
\label{eq:completed-product-monoidal-structure}
\begin{gathered}
 \overline{\otimes}:
 \overline{\mathcal E}\boxtimes\overline{\mathcal E}
 \longrightarrow\overline{\mathcal E},\\
\begin{aligned}
 \overline{\otimes}\circ
 (\iota\boxtimes\iota)
 &=\iota\circ\otimes,\\
 \one_{\overline{\mathcal E}}
 &=\iota(\one_{\mathcal E}).
\end{aligned}
\end{gathered}
\end{equation}
The associator and unitors are the unique extensions of $\iota(\alpha)$,
$\iota(l)$, and $\iota(r)$ by
\eqref{eq:completion-natural-transformations} and its finite-variable
version. These define a monoidal structure on $\overline{\mathcal E}$
for which $\iota$ is strict monoidal.
If $\mathcal E$ is a monoidal $\Pi$-supercategory, then
$\overline{\mathcal E}$ is a monoidal $\Pi$-supercategory with
\begin{equation}
\label{eq:completion-parity}
 \pi_{\overline{\mathcal E}}=\iota(\pi_{\mathcal E}),
 \qquad
 \zeta_{\overline{\mathcal E}}=\iota(\zeta_{\mathcal E}).
\end{equation}
Let $\mathcal D$ be an additive monoidal supercategory with split even
idempotents. For a monoidal superfunctor
$(F,J,\varphi):\mathcal E\to\mathcal D$, every extension $\overline F$
becomes monoidal with tensorator $\overline J$, the unique extension
of $J$ by \eqref{eq:completion-natural-transformations} in two variables,
and unit map $\overline\varphi=\varphi$.
Then $\overline F\circ\iota=F$ as monoidal superfunctors.
\end{proposition}
\begin{proof}
Extend the tensor product in two variables and the constraints by
Proposition~\ref{prop:completion-extension}. Uniqueness preserves
composition and inverses, so the pentagon and triangle follow from
their restrictions to the original objects. The parity data extend
through $\iota$.
Likewise, extend the tensorator to
\begin{equation}
\label{eq:completion-extended-tensorator}
 \overline J:\overline F(-)\otimes\overline F(-)
 \xrightarrow{\sim}\overline F(-\overline\otimes-),
 \qquad\overline\varphi=\varphi.
\end{equation}
Its coherence follows by restricting the equations to the original
objects and using uniqueness of extension.
\end{proof}

\begin{proposition}
\label{prop:completion-uniqueness}
Let $\mathcal E$ be a supercategory and let $\mathcal D$ be an additive
supercategory with split even idempotents.
Let $F:\mathcal E\to\mathcal D$ be a fully faithful superfunctor.
Suppose that every object of $\mathcal D$ is an even retract of a finite
direct sum of objects $FX$. Every extension $\overline F$ in
\eqref{eq:completed-product-extension} is a superequivalence.
If $\mathcal E$ and $\mathcal D$ are monoidal and $F$ is monoidal,
then $\overline F$, with the structure of
Proposition~\ref{prop:monoidal-completion}, is a monoidal superequivalence.
\end{proposition}
\begin{proof}
Full faithfulness extends to finite direct sums and even retracts, so
$\overline F$ is fully faithful. Every object of $\mathcal D$ splits
an even idempotent on some $F^+X$; full faithfulness lifts that
idempotent to $e$ on $X$, giving the object $\overline F(X,e)$.
This proves essential surjectivity. The monoidal assertion follows
from Proposition~\ref{prop:monoidal-completion}.
\end{proof}

\begin{definition}
\label{def:completed-product}
Let $\mathcal A$ and $\mathcal B$ be monoidal supercategories.
The monoidal structure on $\mathcal A\boxtimes\mathcal B$ is defined by
\begin{equation}
\label{eq:product-monoidal-structure}
\begin{aligned}
 (X\boxtimes Y)\otimes(X'\boxtimes Y')
 &=(X\otimes X')\boxtimes(Y\otimes Y'),\\
 (f\boxtimes g)\otimes(f'\boxtimes g')
 &=(-1)^{|g||f'|}(f\otimes f')\boxtimes(g\otimes g')
\end{aligned}
\end{equation}
for homogeneous morphisms. The unit object, associator, and unitors are
\begin{equation}
\label{eq:product-coherence}
\begin{aligned}
 \one_{\mathcal A\boxtimes\mathcal B}
 &=\one_{\mathcal A}\boxtimes\one_{\mathcal B},\\
 \alpha_{X\boxtimes Y,\,X'\boxtimes Y',\,X''\boxtimes Y''}
 &=\alpha^{\mathcal A}_{X,X',X''}\boxtimes
   \alpha^{\mathcal B}_{Y,Y',Y''},\\
 l_{X\boxtimes Y}&=l^{\mathcal A}_X\boxtimes l^{\mathcal B}_Y,
 \qquad r_{X\boxtimes Y}=r^{\mathcal A}_X\boxtimes r^{\mathcal B}_Y.
\end{aligned}
\end{equation}
The \emph{completed product} is
\begin{equation}
\label{eq:completed-product}
 \mathcal A\mathbin{\overline{\boxtimes}}\mathcal B
 =\overline{\mathcal A\boxtimes\mathcal B},
\end{equation}
with monoidal structure chosen as in
Proposition~\ref{prop:monoidal-completion}.
If $\mathcal A$ and $\mathcal B$ are monoidal $\Pi$-supercategories,
the $\Pi$-structure on $\mathcal A\boxtimes\mathcal B$ is defined by
\begin{equation}
\label{eq:product-parity}
\begin{aligned}
 \pi_{\mathcal A\boxtimes\mathcal B}
 &=\pi_{\mathcal A}\boxtimes\one_{\mathcal B},\\
 \zeta_{\mathcal A\boxtimes\mathcal B}
 &=\zeta_{\mathcal A}\boxtimes\id_{\one_{\mathcal B}}.
\end{aligned}
\end{equation}
The $\Pi$-structure on
$\mathcal A\mathbin{\overline{\boxtimes}}\mathcal B$
is then defined by Proposition~\ref{prop:monoidal-completion}.
\end{definition}

\begin{proposition}
\label{thm:completed-product}
If $\mathcal A,\mathcal B\in\W_k$, then
\begin{equation}
\label{eq:completed-product-closure}
 \mathcal A\mathbin{\overline{\boxtimes}}\mathcal B\in\W_k,
\end{equation}
and there is a group isomorphism
\begin{equation}
\label{eq:completed-product-group}
 G_{\mathcal A}\times G_{\mathcal B}
 \xrightarrow{\sim}G_{\mathcal A\mathbin{\overline{\boxtimes}}\mathcal B},
 \qquad ([X],[Y])\longmapsto[S],
\end{equation}
where $S$ is any simple summand of $X\boxtimes Y$.
Under this isomorphism,
\begin{equation}
\label{eq:completed-product-type}
 \rho_{\mathcal A\mathbin{\overline{\boxtimes}}\mathcal B}(g,h)
 =\rho_{\mathcal A}(g)+\rho_{\mathcal B}(h).
\end{equation}
\end{proposition}
\begin{proof}
Write $\mathcal T=\mathcal A\mathbin{\overline\boxtimes}\mathcal B$.
Choose simple $X_g\in\mathcal A_g$, $Y_h\in\mathcal B_h$ and put
$Z_{g,h}=X_g\boxtimes Y_h$. Then
\[
 E_{g,h}=\End(Z_{g,h})
 \cong D_{\rho_{\mathcal A}(g)}\gt D_{\rho_{\mathcal B}(h)}
\]
is central simple of type $\rho_{\mathcal A}(g)+\rho_{\mathcal B}(h)$.
Every object is an even retract of sums of
$P=\bigoplus_{g,h}(Z_{g,h}\oplus\Pi Z_{g,h})$.
Since $\End(P)\cong\bigoplus_{g,h}M_{1|1}(k)\gt E_{g,h}$,
Proposition~\ref{prop:central-simple-superalgebras}(6) makes
$E=\End(P)_0$ semisimple, and
\[
 \Hom_0(P,-):\underline{\mathcal T}
 \simeq\mathsf{proj}_{\mathrm{fd}}(E)=\mathsf{mod}_{\mathrm{fd}}(E).
\]
Thus $\mathcal T$ is finite semisimple.
Duals are componentwise on pairs and sums. For an even retract
$(Z,e)$, write $e^*$ for the dual idempotent and set
\[
 (Z,e)^*=(Z^*,e^*),\qquad
 \ev_{(Z,e)}=\ev_Z(e^*\otimes e),\qquad
 \coev_{(Z,e)}=(e\otimes e^*)\coev_Z.
\]
These give rigidity; the unit has endomorphism algebra $k$.
Graded simplicity of $E_{g,h}$ makes all simple retracts of
$Z_{g,h}$ and $\Pi Z_{g,h}$ homogeneously isomorphic. Their
endomorphism superalgebras are isomorphic to
$D_{\rho_{\mathcal A}(g)+\rho_{\mathcal B}(h)}$ by
Lemma~\ref{lem:matrix-corners}. Hence the blocks satisfy
\[
 \mathcal T=\bigoplus_{g,h}\mathcal T_{g,h},\qquad
 \mathcal T_{g,h}\otimes\mathcal T_{g',h'}\subset\mathcal T_{gg',hh'},
 \qquad(\mathcal T_{g,h})^*\subset\mathcal T_{g^{-1},h^{-1}}.
\]
The neutral block consists of sums of $\one$ and $\pi$, proving
weak invertibility. The same block decomposition and the types of
these simple retracts give the asserted group isomorphism and type formula.
\end{proof}

\subsection{Stacking of \texorpdfstring{$G$}{G}-labeled categories}
\label{subsec:stacking-labeled-categories}

Fixing a group $G$ with an isomorphism $G \to G_\mathcal{A}$, we remove the group isomorphism redundancy and consider a class $\mathcal W_k(G)$ of categories $\mathcal A$ in $\mathcal W_k$ with $G_\mathcal{A} \cong G$. 
We introduce a product $\mathcal S(-,-)$, called \emph{stacking}, on $\mathcal W_k(G)$ and show that the equivalence classes in $\mathcal W_k(G)$ form an abelian group.

\begin{definition}
\label{def:stacking}
Let $G$ be a finite group. A category $\mathcal A\in\W_k$ is
\emph{$G$-labeled} if it is equipped with a group isomorphism
$\lambda_{\mathcal A}:G\to G_{\mathcal A}$; write $\W_k(G)$ for the class of
$G$-labeled categories and
$\mathcal A_{(g)}=\mathcal A_{\lambda_{\mathcal A}(g)}$ for the components.
For $\mathcal A,\mathcal B\in\W_k(G)$, write $\mathcal A\simeqG\mathcal B$ if
there is a monoidal superequivalence $E:\mathcal A\to\mathcal B$ such that
$E(\mathcal A_{(g)})\subset\mathcal B_{(g)}$ for all $g\in G$; such an $E$
is called \emph{label-preserving}.
The \emph{stack} $\mathcal S(\mathcal A,\mathcal B)$ is the full
subsupercategory of $\mathcal A\mathbin{\overline{\boxtimes}}\mathcal B$
whose objects are the even retracts of finite direct sums of the objects
\begin{equation}
\label{eq:diagonal-pairs}
 X\boxtimes Y,\qquad X\in\mathcal A_{(g)},\quad Y\in\mathcal B_{(g)},\quad g\in G.
\end{equation}
Let $\mathcal D(\mathcal A,\mathcal B)$ denote the monoidal full $\Pi$-subsupercategory
on the objects \eqref{eq:diagonal-pairs}, with
inclusion $j:\mathcal D(\mathcal A,\mathcal B)\hookrightarrow\mathcal S(\mathcal A,\mathcal B)$.
\end{definition}

\begin{proposition}
\label{thm:stack-existence}
If $\mathcal A,\mathcal B\in\W_k(G)$, then $\mathcal S(\mathcal A, \mathcal B) \in \mathcal W_k(G)$ such that
\begin{equation}
\label{eq:stack-labeling}
 \lambda_{\mathcal S(\mathcal A,\mathcal B)}(g)=[S],\qquad
 S\ \text{any simple summand of}\ X\boxtimes Y,\quad
 X\in\mathcal A_{(g)},\ Y\in\mathcal B_{(g)}\ \text{simple}.
\end{equation}
\end{proposition}
\begin{proof}
Write $\mathcal T=\mathcal A\mathbin{\overline\boxtimes}\mathcal B$.
By Proposition~\ref{thm:completed-product}, the stack is the sum
of the blocks indexed by the diagonal subgroup:
\[
 \mathcal S(\mathcal A,\mathcal B)
 =\bigoplus_{g\in G}\mathcal T_{(\lambda_{\mathcal A}(g),\lambda_{\mathcal B}(g))}.
\]
Tensor products and duals follow multiplication and inversion in
this subgroup, giving the asserted labeling and closure.
\end{proof}

\begin{lemma}
\label{lem:stack-functoriality}
If $\mathcal A\simeqG\mathcal A'$ and $\mathcal B\simeqG\mathcal B'$ in
$\W_k(G)$, then
$\mathcal S(\mathcal A,\mathcal B)\simeqG\mathcal S(\mathcal A',\mathcal B')$.
\end{lemma}
\begin{proof}
The product of the given equivalences is a monoidal
superequivalence on pairs. Extend it to completions by
Propositions~\ref{prop:monoidal-completion}
and~\ref{prop:completion-uniqueness}. Since it preserves labels,
it restricts to an equivalence of the diagonal blocks, hence of
the stacks.
\end{proof}

\begin{definition}
\label{def:graded-svec}
For a finite group $G$, let $\sVecG$ be the supercategory of families
$V=(V_g)_{g\in G}$ of finite-dimensional super vector spaces, with
morphism spaces $\bigoplus_{g\in G}\Hom_k(V_g,W_g)$ graded by the parity
of linear maps. For a family indexed by $G$, $X[g]$ denotes the family
with $X$ at $g$ and $0$ elsewhere.
It is a monoidal $\Pi$-supercategory with tensor product
\begin{equation}
\label{eq:graded-svec-tensor}
 (V\otimes W)_\ell=\bigoplus_{gh=\ell}V_g\otimes_kW_h,\qquad
 (f\otimes f')(v\otimes w)=(-1)^{|f'||v|}f(v)\otimes f'(w),
\end{equation}
for homogeneous $f,f'$ and $v$, unit $k[1]$, associator and unitors
acting on each summand by
\begin{equation}
\label{eq:graded-svec-constraints}
 (v\otimes w)\otimes u\longmapsto v\otimes(w\otimes u),\qquad
 c\otimes v\longmapsto cv,\qquad v\otimes c\longmapsto cv
 \qquad(c\in k),
\end{equation}
parity object $\pi=k^{0|1}[1]$, where $k^{0|1}$ is $k$ placed in odd
degree, and $\zeta:k^{0|1}[1]\to k[1]$ the identity map of $k$, which is
odd.
\end{definition}

\begin{lemma}
\label{lem:stack-unit}
Let $G$ be a finite group. We have $\sVecG\in\W_k(G)$ with
$\lambda_{\sVecG}(g)=[k[g]]$, and $\mathcal A\simeqG\mathcal S(\mathcal A,\sVecG)$
for every $\mathcal A\in\W_k(G)$.
\end{lemma}
\begin{proof}
Definition~\ref{def:graded-svec} gives $\sVecG\in\W_k(G)$.
On objects in a single component, use the fully faithful monoidal
superfunctor
\[
 X\in\mathcal A_{(g)}\longmapsto X\boxtimes k[g],
 \qquad f\longmapsto f\boxtimes\id.
\]
Its tensorator uses $k\otimes_k k\cong k$.
Its image generates the stack under sums and even retracts, since
every super vector space is a sum of even and odd lines and
$X\boxtimes k^{0|1}[g]\cong(X\otimes\pi)\boxtimes k[g]$ evenly.
The inclusion of these homogeneous objects also generates
$\mathcal A$. Completion uniqueness therefore gives the desired
label-preserving monoidal equivalence.
\end{proof}

\begin{definition}
\label{def:opposite}
For a monoidal $\Pi$-supercategory $\mathcal A$, the \emph{opposite}
$\mathcal A^{\mathrm{sop}}$ has the same objects, morphism spaces
$\mathcal A^{\mathrm{sop}}(X,Y)=\mathcal A(Y,X)$, composition and tensor
product
\begin{equation}
\label{eq:opposite-structure}
 g\circ^{\mathrm{sop}}f=(-1)^{|f||g|}f\circ g,\qquad
 X\otimes^{\mathrm{sop}}Y=X\otimes Y,\qquad
 f\otimes^{\mathrm{sop}}g=f\otimes g,
\end{equation}
associator $\alpha^{-1}$, unitors $l^{-1}$ and $r^{-1}$, unit $\one$,
parity object $\pi$, and $\zeta^{\mathrm{sop}}=\zeta^{-1}$.
\end{definition}

\begin{lemma}
\label{lem:opposite}
$\mathcal A^{\mathrm{sop}}$ is a monoidal $\Pi$-supercategory. If
$\mathcal A\in\W_k(G)$, then $\mathcal A^{\mathrm{sop}}\in\W_k(G)$ with the same
simple objects and labeling.
\end{lemma}
\begin{proof}
The Koszul sign in opposite composition gives associativity and
super interchange; inverting the even constraints preserves their
coherence. The even category is $(\underline{\mathcal A})^{\mathrm{op}}$,
with left and right duals exchanged, and
\[
 \End_{\mathcal A^{\mathrm{sop}}}(X)=\End_{\mathcal A}(X)^{\mathrm{sop}}.
\]
Both duals of a simple $X$ represent $[X]^{-1}$, so its product
with either dual lies in the neutral component. Thus weak
invertibility and the split unit persist. Tensor products of
objects are unchanged, giving the same simple classes and labeling.
\end{proof}

\begin{lemma}
\label{lem:stack-inverse}
For $\mathcal A\in\W_k(G)$,
$\mathcal S(\mathcal A,\mathcal A^{\mathrm{sop}})\simeqG\sVecG$.
\end{lemma}
\begin{proof}
On diagonal pairs in degree $g$, define
\[
 H(X\boxtimes Y)=\mathcal A(Y,X)[g],\qquad
 H(a\boxtimes b)(f)=(-1)^{|b||f|}afb.
\]
Super interchange and opposite composition make $H$ a superfunctor.
Its tensorator is $f\otimes f'\mapsto f\otimes f'$, invertible by
Lemma~\ref{lem:component-tensor-homs}; associativity and units follow
from those of $\mathcal A$.
For a simple generator $S$ with $D=\End(S)$, full faithfulness is
the isomorphism
\[
 D\gt D^{\mathrm{sop}}\xrightarrow{\sim}\End_k(D),\qquad
 a\gt b\longmapsto(f\mapsto(-1)^{|b||f|}afb).
\]
Extend this to sums and parity shifts. The images of
$S\boxtimes S$ and $(\Pi S)\boxtimes S$ contain the even and odd
lines at $g$ as retracts, so they generate $\sVecG$.
Monoidal completion uniqueness gives the asserted equivalence.
\end{proof}

\begin{lemma}
\label{lem:stack-products}
For $\mathcal A,\mathcal B,\mathcal C\in\W_k(G)$,
\begin{equation}
\label{eq:stack-commutativity-associativity}
 \mathcal S(\mathcal A,\mathcal B)\simeqG\mathcal S(\mathcal B,\mathcal A),
 \qquad
 \mathcal S(\mathcal S(\mathcal A,\mathcal B),\mathcal C)
 \simeqG\mathcal S(\mathcal A,\mathcal S(\mathcal B,\mathcal C)).
\end{equation}
\end{lemma}
\begin{proof}
The super swap
\[
 X\boxtimes Y\longmapsto Y\boxtimes X,\qquad
 f\boxtimes u\longmapsto(-1)^{|f||u|}u\boxtimes f
\]
is strict monoidal and preserves diagonal pairs. It extends to
completions and proves commutativity.
For associativity, both iterated stacks are completions of the
same category of diagonal triples
$(X\boxtimes Y)\boxtimes Z\cong X\boxtimes(Y\boxtimes Z)$.
Indeed, expanding sums and even retracts in either order gives
the same generators. Proposition~\ref{prop:completion-uniqueness}
gives a monoidal equivalence, preserving the common label $g$.
\end{proof}

\begin{theorem}
\label{thm:labeled-stacking-group}
Let $G$ be a finite group. The set $\W_k(G)/{\simeqG}$ is an abelian group
under
\begin{equation}
\label{eq:labeled-product}
 [\mathcal A]\cdot[\mathcal B]=[\mathcal S(\mathcal A,\mathcal B)],
\end{equation}
with unit $[\sVecG]$ and inverse $[\mathcal A]^{-1}=[\mathcal A^{\mathrm{sop}}]$.
\end{theorem}
\begin{proof}
Proposition~\ref{thm:stack-existence} and
Lemma~\ref{lem:stack-functoriality} make the product well defined.
Lemma~\ref{lem:stack-products} gives associativity and commutativity,
Lemma~\ref{lem:stack-unit} gives the unit, and
Lemmas~\ref{lem:opposite}--\ref{lem:stack-inverse} give inverses.
\end{proof}

\subsection{Exchange data}
\label{subsec:labeled-pullbacks-exchange-data}

To write the multiplication law of $\mathcal W_k(G)$ in terms of twists, we prepare an \emph{exchange data} on a reference category $\mathcal C$. 
This data is akin to braiding, but with the difference that the component morphisms may be odd, and the two composites in each hexagon may differ by a scalar multiple. 
In this paper, we present this data solely as what is required for direct computation; a more category-theoretic investigation may be found in \cite{SamSnowden}, and a 2-categorical interpretation remains a topic for future discussion.

\begin{lemma}
\label{lem:reference-exchange}
Let $\mathcal C$ be a Brauer--Wall reference category and let
$a,b\in\BW(k)$. Define superfunctors
$T,T^{\mathrm{op}}:\mathcal C_a\boxtimes\mathcal C_b\to\mathcal C_{a+b}$
by
\begin{equation}
\label{eq:exchange-functors}
\begin{aligned}
 T(X\boxtimes Y)&=X\otimes Y,&
 T(f\boxtimes g)&=f\otimes g,\\
 T^{\mathrm{op}}(X\boxtimes Y)&=Y\otimes X,&
 T^{\mathrm{op}}(f\boxtimes g)&=(-1)^{|f||g|}g\otimes f
\end{aligned}
\end{equation}
for homogeneous morphisms $f,g$.
Then $\Nat(T,T^{\mathrm{op}})$ is a one-dimensional super vector space,
and its nonzero elements are supernatural isomorphisms.
\end{lemma}
\begin{proof}
Super interchange makes $T^{\mathrm{op}}$ a superfunctor.
At simple generators, both $T$ and $T^{\mathrm{op}}$ induce
endomorphism-algebra isomorphisms by
Proposition~\ref{prop:simple-tensor}, with the graded swap for the
latter. Proposition~\ref{prop:comparison-lines} therefore applies.
For a degree-$d$ comparison, supernaturality reads
\begin{equation}
\label{eq:exchange-supernaturality}
 R_{X',Y'}(f\otimes g)
 =(-1)^{d(|f|+|g|)+|f||g|}(g\otimes f)R_{X,Y}.
\end{equation}
\end{proof}

\begin{definition}
\label{def:exchange-datum}
Let $\mathcal C$ be a Brauer--Wall reference category. The
\emph{exchange degree} of $\mathcal C$ is the map
$b_{\mathcal C}:\BW(k)^2\to\Z/2\Z$ sending $(a,b)$ to the degree of
$\Nat(T,T^{\mathrm{op}})$ in Lemma~\ref{lem:reference-exchange}.
An \emph{exchange data} of $\mathcal C$ is a triple $(R,\beta_1,\beta_2)$
consisting of a nonzero element $R$ of $\Nat(T,T^{\mathrm{op}})$ for every
$a,b\in\BW(k)$, with components
\begin{equation}
\label{eq:exchange-isomorphism}
 R_{X,Y}:X\otimes Y\longrightarrow Y\otimes X
 \qquad(X\in\mathcal C_a,\ Y\in\mathcal C_b)
\end{equation}
and the unit values
\begin{equation}
\label{eq:exchange-unit}
 R_{\one,X}=r_X^{-1}l_X,\qquad R_{X,\one}=l_X^{-1}r_X,
\end{equation}
and maps $\beta_1,\beta_2:\BW(k)^3\to k^\times$, the
\emph{hexagon scalars}, such that for $X\in\mathcal C_a$,
$Y\in\mathcal C_b$, $Z\in\mathcal C_c$
\begin{equation}
\label{eq:exchange-left-hexagon}
 \alpha_{Y,Z,X}R_{X,Y\otimes Z}\alpha_{X,Y,Z}
 =\beta_1(a;b,c)\,
   (\id_Y\otimes R_{X,Z})\alpha_{Y,X,Z}(R_{X,Y}\otimes\id_Z),
\end{equation}
\begin{equation}
\label{eq:exchange-right-hexagon}
\begin{aligned}
 R_{X\otimes Y,Z}
 ={}&\beta_2(a,b;c)\,
 \alpha_{Z,X,Y}(R_{X,Z}\otimes\id_Y)\alpha_{X,Z,Y}^{-1}\\
 &\qquad\circ(\id_X\otimes R_{Y,Z})\alpha_{X,Y,Z}.
\end{aligned}
\end{equation}
\end{definition}

\begin{proposition}
\label{prop:exchange-existence}
Every Brauer--Wall reference category $\mathcal C$ has an exchange data.
For a family $R$ as in Definition~\ref{def:exchange-datum}, the hexagon
scalars exist, are unique, and equal $1$ whenever an argument is $0$.
The exchange degree satisfies
\begin{equation}
\label{eq:exchange-degree-laws}
 b_{\mathcal C}(a,b)=b_{\mathcal C}(b,a),\qquad
 b_{\mathcal C}(a+b,c)=b_{\mathcal C}(a,c)+b_{\mathcal C}(b,c),\qquad
 b_{\mathcal C}(a,a)=\epsilon(a),
\end{equation}
and $b_{\mathcal C}(a,b)=\epsilon(a)\epsilon(b)$ for all $a,b$ if $\BW(k)$ is
cyclic.
\end{proposition}
\begin{proof}
Choose nonzero exchange comparisons, using the unitors when one
argument is $\one$. The comparison-line uniqueness extends these
unit values to the neutral components.
Both routes in each hexagon lie in one super line, so their degrees
agree and their ratio defines a unique scalar. This gives
biadditivity of $b_{\mathcal C}$; evaluating at a unit gives
normalization of the scalars. Since $R_{Y,X}R_{X,Y}$ is a scalar
endomorphism of the tensor functor, its degree is zero, giving symmetry.
For a simple $S_a$, \eqref{eq:exchange-supernaturality} identifies
$R_{S_a,S_a}$ with an invertible homogeneous element of
$D_a\gt D_a$ satisfying \eqref{eq:signed-swap-implementation}.
Proposition~\ref{prop:swap-parity} gives
$b_{\mathcal C}(a,a)=\epsilon(a)$.
If $\BW(k)$ is cyclic with generator $r$, then
\[
 b_{\mathcal C}(mr,nr)=mn\epsilon(r)
 =\epsilon(mr)\epsilon(nr).
\]
\end{proof}

\begin{definition}
\label{def:exchange-rescaling}
Let $(R,\beta_1,\beta_2)$ be an exchange data of a
Brauer--Wall reference category $\mathcal C$, and let
$\lambda:\BW(k)^2\to k^\times$ satisfy $\lambda(0,b)=\lambda(a,0)=1$ for
all $a,b$. Put $(\lambda R)_{X,Y}=\lambda(a,b)R_{X,Y}$ for
$X\in\mathcal C_a$, $Y\in\mathcal C_b$, and
\begin{equation}
\label{eq:factor-rescaling}
 \beta^\lambda_1(a;b,c)
 =\beta_1(a;b,c)\frac{\lambda(a,b+c)}{\lambda(a,b)\lambda(a,c)},\qquad
 \beta^\lambda_2(a,b;c)
 =\beta_2(a,b;c)\frac{\lambda(a+b,c)}{\lambda(a,c)\lambda(b,c)}.
\end{equation}
Then $(\lambda R,\beta^\lambda_1,\beta^\lambda_2)$ is the
\emph{rescaling} of the exchange data by $\lambda$.
\end{definition}

\begin{proposition}
\label{prop:exchange-identities}
Let $(R,\beta_1,\beta_2)$ be an exchange data of a
Brauer--Wall reference category $\mathcal C$. Every rescaling of it is
an exchange data of $\mathcal C$, and every exchange data of
$\mathcal C$ is its rescaling by a unique $\lambda$.
\end{proposition}
\begin{proof}
Each exchange space is a line, so any other normalized exchange data has
$R'=\lambda R$ for a unique normalized $\lambda$.
Substitution in the two hexagons gives exactly the ratios
in \eqref{eq:factor-rescaling}.
\end{proof}

\begin{proposition}
\label{prop:hexagon-identities}
Let $(R,\beta_1,\beta_2)$ be an exchange data of a
Brauer--Wall reference category $\mathcal C$. For $a,b,c,d\in\BW(k)$,
\begin{equation}
\label{eq:exchange-datum-identities}
\begin{aligned}
 \beta_1(a;b,c)\beta_1(a;b+c,d)
 &=\beta_1(a;b,c+d)\beta_1(a;c,d),\\
 \beta_2(a+b,c;d)\beta_2(a,b;d)
 &=\beta_2(b,c;d)\beta_2(a,b+c;d),\\
 \frac{\beta_1(b;c,d)\beta_1(a;c,d)}{\beta_1(a+b;c,d)}
 &=(-1)^{b_{\mathcal C}(a,c)\,b_{\mathcal C}(b,d)}
   \frac{\beta_2(a,b;c)\beta_2(a,b;d)}{\beta_2(a,b;c+d)},\\
 \frac{\beta_1(a;b,c)}{\beta_1(a;c,b)}
 &=(-1)^{b_{\mathcal C}(a+b,c)\,b_{\mathcal C}(a,b)
         +b_{\mathcal C}(a,b+c)\,b_{\mathcal C}(b,c)}
   \frac{\beta_2(b,a;c)}{\beta_2(a,b;c)}.
\end{aligned}
\end{equation}
\end{proposition}
\begin{proof}
Use the two hexagons to express exchange maps involving tensor
products, such as $R_{X,Y\otimes Z}$ and $R_{X\otimes Y,Z}$,
in terms of exchanges between individual factors, and cancel
common isomorphisms. For the first two identities, the two expansions
of $R_{X,(Y\otimes Z)\otimes W}$ and
$R_{(X\otimes Y)\otimes Z,W}$ agree by the pentagon.
For the third, the two expansions of $R_{X\otimes Y,Z\otimes W}$
differ only by interchanging $R_{Y,W}$ and $R_{X,Z}$, giving
$(-1)^{b_{\mathcal C}(b,d)b_{\mathcal C}(a,c)}$.
For the fourth, apply \eqref{eq:exchange-supernaturality} to
$R_{X,Y}$ and to $R_{Y,Z}$, then use the hexagons to expand
the exchanges involving tensor products.
The two supernaturality signs have exponents
\[
 b_{\mathcal C}(a+b,c)b_{\mathcal C}(a,b),\qquad
 b_{\mathcal C}(a,b+c)b_{\mathcal C}(b,c),
\]
whose sum gives the stated exponent.
\end{proof}

\subsection{Labeled pullback and stacking of twists}
\label{subsec:stacking-twists-degree}

Pullbacks of twists on $G$ carry a natural $G$-labeling.
Given an exchange data on $\mathcal C$, we compute a twist $s\star t$
representing the stack of $\mathcal P_{\mathcal C}(s)$ and
$\mathcal P_{\mathcal C}(t)$.  

We use the stacking formula to describe the resulting group $\Tcal_k(G)$ by two successive
group extensions. Recording the Brauer--Wall type homomorphism $\rho$
gives a surjection onto $Z^1(G,\BW(k))$ with kernel $\SH_k(G)$.
Within this kernel, a parity class $[\nu]\in H^2(G,\Z/2\Z)$ occurs
precisely when the pentagon equation
$\delta\omega=(-1)^{\nu\smile\nu}$ admits a solution. For each such
parity class, the choices of scalar associator up to gauge form a
torsor for $H^3(G,k^\times)$.

\begin{definition}
\label{def:labeled-pullbacks}
For a finite group $G$, let $\Tw_k(G)$ be the set of twists with group $G$;
for $s,t\in\Tw_k(G)$, write $s\sim_Gt$ if there is a gauge
$(\id_G,\lambda,\beta)$ from $s$ to $t$; by
Proposition~\ref{prop:gauge-equivalence} this is an equivalence relation,
and we put $\Tcal_k(G)=\Tw_k(G)/{\sim_G}$. Put $e_G=(G,0,0,1)$ and
$p_\rho=(G,\rho,0,1)$ for $\rho\in Z^1(G,\BW(k))$.
For a Brauer--Wall reference category $\mathcal C$ and $s\in\Tw_k(G)$, we
have $\mathcal P_{\mathcal C}(s)\in\W_k(G)$ with $\lambda(g)=[X_g]$ by
\eqref{eq:pullback-group}.
\end{definition}

\begin{lemma}
\label{lem:trivial-pullback}
Let $\mathcal C$ be a Brauer--Wall reference category and let $G$ be a
finite group. Then
\begin{equation}
\label{eq:trivial-pullback}
 \mathcal P_{\mathcal C}(e_G)\simeqG\sVecG.
\end{equation}
\end{lemma}
\begin{proof}
On the neutral component, \eqref{eq:component-module-equivalence}
with $S=\one$ gives the superequivalence $\mathcal C(\one,-)$
to super vector spaces. Its tensorator is
\[
 u\otimes v\longmapsto(u\otimes v)l_{\one}^{-1},
\]
an isomorphism by Lemma~\ref{lem:component-tensor-homs}; coherence
follows from that of the unit. Apply this equivalence componentwise
to the families in $\mathcal P_{\mathcal C}(e_G)$.
For $\nu=0$, $\omega=1$, all inserted parity factors are units,
so it is monoidal and sends $\one[g]$ to $k[g]$.
\end{proof}

\begin{definition}
\label{def:cochain-operations}
Let $\mathcal C$ be a Brauer--Wall reference category with an exchange
data $(R,\beta_1,\beta_2)$. For $s,t\in\Tw_k(G)$ define quadruples $s\star t$ and $s^\vee$ by
\begin{equation}
\label{eq:twist-product}
 s\star t=\bigl(G,\ \rho_s+\rho_t,\ \nu_s+\nu_t+c_{s,t},\
   \omega_s\omega_t\,K_{s,t}\,(-1)^{\Xi_{s,t}}\bigr),
\end{equation}
\begin{equation}
\label{eq:twist-inverse}
 s^\vee=\bigl(G,\ -\rho_s,\ \nu_s+c_{s,s},\
   \bigl(\omega_s\,K^\vee_s\,(-1)^{\Xi^\vee_s}\bigr)^{-1}\bigr),
\end{equation}
where
\begin{equation}
\label{eq:stacking-factor-cochains}
\begin{aligned}
 c_{s,t}(g,h)&=b_{\mathcal C}(\rho_t(g),\rho_s(h)),\qquad
 K_{s,t}(g,h,\ell)
 =\frac{\beta_2(\rho_t(g),\rho_t(h);\rho_s(\ell))}
        {\beta_1(\rho_t(g);\rho_s(h),\rho_s(\ell))},\\
 \Xi_{s,t}(g,h,\ell)&=\bigl(\nu_s\smile_1\nu_t+c_{s,t}\smile_1(\nu_s+\nu_t)\bigr)(g,h,\ell)
   +c_{s,t}(g,h)\,c_{s,t}(h,\ell),\\
 K^\vee_s(g,h,\ell)
 &=\frac{\beta_2(\rho_s(g),\rho_s(h);-\rho_s(\ell))}
        {\beta_1(\rho_s(g);-\rho_s(h),-\rho_s(\ell))},\\
 \Xi^\vee_s(g,h,\ell)&=\bigl(\nu_s\smile_1\nu_s+c_{s,s}\smile_1(\nu_s+c_{s,s})\bigr)(g,h,\ell)
   +c_{s,s}(g,h)\,c_{s,s}(h,\ell),
\end{aligned}
\end{equation}
with $\smile_1$ as in \eqref{eq:cup-one}.
By definition, $c_{s^\vee,s}=c_{s,s}$, $K_{s^\vee,s}=K^\vee_s$, and
$\Xi_{s^\vee,s}=\Xi^\vee_s$.
\end{definition}

\begin{lemma}
\label{lem:product-gauges}
For $s,t\in\Tw_k(G)$ we have $s\star t\in\Tw_k(G)$ and $s^\vee\in\Tw_k(G)$, and
\begin{equation}
\label{eq:stacking-unit-inverse}
 s^\vee\star s=e_G.
\end{equation}
\end{lemma}
\begin{proof}
Put $c=c_{s,t}$, $T(g,h,\ell)=c(g,h)c(h,\ell)$, and
$\hat c(g,h,\ell,m)=c(g,\ell)c(h,m)$.
Biadditivity of $b_{\mathcal C}$ makes $c$ a normalized
bimultiplicative $2$-cocycle. The hexagon identities and the
cochain identities of Appendix~\ref{app:cochains} give
\begin{align*}
 \delta K_{s,t}&=(-1)^{\hat c},&
 \delta T&=c\smile c+\hat c,\\
 \delta\Xi_{s,t}
 &=\nu_{s\star t}\smile\nu_{s\star t}
   +\nu_s\smile\nu_s+\nu_t\smile\nu_t+\hat c.
\end{align*}
Consequently
\[
 \delta\bigl(K_{s,t}(-1)^{\Xi_{s,t}}\bigr)
 =(-1)^{\nu_{s\star t}\smile\nu_{s\star t}
          +\nu_s\smile\nu_s+\nu_t\smile\nu_t},
\]
which proves the twist equation for $s\star t$.
This identity uses only the degree cocycles, so applying it to
$(s^\vee,s)$ also proves the twist equation for $s^\vee$ without
assuming it in advance. All cochains are normalized, and
\eqref{eq:twist-inverse} cancels both the degree and scalar parts
of $s^\vee\star s$, giving $e_G$.
\end{proof}

\begin{proposition}
\label{prop:stacking-comparison-coherence}
Let $\mathcal C$ and $(R,\beta_1,\beta_2)$ be as in
Definition~\ref{def:cochain-operations} and let $s,t\in\Tw_k(G)$. There is
a fully faithful monoidal superfunctor
\begin{equation}
\label{eq:stacking-comparison-monoidal}
 (F_0,J,\varphi):\mathcal D(\mathcal P_{\mathcal C}(s),\mathcal P_{\mathcal C}(t))
 \longrightarrow\mathcal P_{\mathcal C}(s\star t)
\end{equation}
whose image generates $\mathcal P_{\mathcal C}(s\star t)$ under finite direct
sums and even retracts. Its data are as follows. For morphisms
$f:X\to X'$ and $f':Y\to Y'$,
\begin{equation}
\label{eq:stacking-comparison-functor}
\begin{aligned}
 F_0\bigl(X[g]\boxtimes Y[g]\bigr)&=(X\otimes Y)[g],\\
 F_0\bigl(f[g]\boxtimes f'[g]\bigr)&=(f\otimes f')[g],
\end{aligned}
\end{equation}
extended linearly, where $f[g]$ denotes the morphism supported at $g$
with component $f$. For objects $X,Y,Z,W$ of components of
$\mathcal C$, let $\chi_{X,Y,Z,W}$ be the middle exchange
\begin{equation}
\label{eq:middle-exchange}
\begin{aligned}
 \chi_{X,Y,Z,W}
 ={}&\alpha_{X,Z,Y\otimes W}^{-1}
       (\id_X\otimes\alpha_{Z,Y,W})
       (\id_X\otimes(R_{Y,Z}\otimes\id_W))\\
 &\circ(\id_X\otimes\alpha_{Y,Z,W}^{-1})
       \alpha_{X,Y,Z\otimes W}\\
 &\colon (X\otimes Y)\otimes(Z\otimes W)
 \longrightarrow (X\otimes Z)\otimes(Y\otimes W).
\end{aligned}
\end{equation}
Write $e_{d,Z}=r_Z(\id_Z\otimes\zeta_d)$ as in
Proposition~\ref{prop:reconstruction}. For $X\in\mathcal C_{\rho_s(g)}$,
$Y\in\mathcal C_{\rho_t(g)}$, $X'\in\mathcal C_{\rho_s(h)}$,
$Y'\in\mathcal C_{\rho_t(h)}$, put $a=\nu_s(g,h)$, $b=\nu_t(g,h)$,
$n=\nu_{s\star t}(g,h)$, and
\begin{equation}
\label{eq:stacking-tensorator}
\begin{aligned}
 J_{g,h}(X,Y;X',Y')
 ={}&(e_{a,X\otimes X'}\otimes e_{b,Y\otimes Y'})^{-1}\\
 &\circ\chi_{X,Y,X',Y'}\,
       e_{n,(X\otimes Y)\otimes(X'\otimes Y')}\\
 &\colon ((X\otimes Y)\otimes(X'\otimes Y'))\otimes P_n\\
 &\longrightarrow
 ((X\otimes X')\otimes P_a)\otimes((Y\otimes Y')\otimes P_b).
\end{aligned}
\end{equation}
The tensorator and unit map are
\begin{equation}
\label{eq:stacking-comparison-structure}
\begin{aligned}
 J_{X[g]\boxtimes Y[g],\,X'[h]\boxtimes Y'[h]}
   &=J_{g,h}(X,Y;X',Y')[gh],\\
 \varphi&=(l_{\one}^{-1})[1].
\end{aligned}
\end{equation}
\end{proposition}
\begin{proof}
Put $c=c_{s,t}$ and $n=\nu_s+\nu_t+c$.
The middle exchange in $J_{g,h}$ has degree $c(g,h)$, so
$|J_{g,h}|=n(g,h)+\nu_s(g,h)+\nu_t(g,h)+c(g,h)=0$.
The same cancellation, using \eqref{eq:exchange-supernaturality},
proves supernaturality. In deleting parity factors, retain the sign
\[
 (e_a\otimes e_b)^{-1}=(-1)^{ab}(e_a^{-1}\otimes e_b^{-1}).
\]
For associativity, delete all parity factors from both routes
using \eqref{eq:gauge-ordered-deletions} and
\eqref{eq:gauge-associator-deletion}. Expanding the two exchange
maps involving tensor products using the hexagons contributes
\[
 K_{s,t}(g,h,\ell)(-1)^{c(g,h)c(h,\ell)};
\]
the sign comes from interchanging the two exchange maps of degrees
$c(g,h)$ and $c(h,\ell)$ acting on disjoint tensor factors.
The remaining deletion signs have total exponent
$(\nu_s\smile_1\nu_t+c\smile_1(\nu_s+\nu_t))(g,h,\ell)$.
Thus the routes agree precisely when
\[
 \omega_{s\star t}=\omega_s\omega_tK_{s,t}(-1)^{\Xi_{s,t}},
\]
as prescribed. Unit coherence follows from normalization and
\eqref{eq:exchange-unit}.
Full faithfulness is Lemma~\ref{lem:component-tensor-homs} on each
support. A tensor product of simple generators is nonzero of type
$\rho_s(g)+\rho_t(g)$, so it and its parity shift generate that
component under sums and even retracts. This proves the generation
assertion in every degree.
\end{proof}

\begin{theorem}
\label{thm:stacking-group-categorical}
Let $\mathcal C$ be a Brauer--Wall reference category with an exchange
data and let $G$ be a finite group. For $s,t\in\Tw_k(G)$ there is a
label-preserving monoidal superequivalence
\begin{equation}
\label{eq:stacking-comparison}
 F:\mathcal S(\mathcal P_{\mathcal C}(s),\mathcal P_{\mathcal C}(t))
 \xrightarrow{\ \simeqm\ }\mathcal P_{\mathcal C}(s\star t).
\end{equation}
The set $\Tcal_k(G)$ is an abelian group under
\begin{equation}
\label{eq:class-operations}
 [s]\cdot[t]=[s\star t],\qquad 1=[e_G],\qquad [s]^{-1}=[s^\vee],
\end{equation}
and the map
\begin{equation}
\label{eq:labeled-bijection}
 \Tcal_k(G)\longrightarrow\W_k(G)/{\simeqG},\qquad
 [s]\longmapsto[\mathcal P_{\mathcal C}(s)],
\end{equation}
is a group isomorphism onto the stacking group of
Theorem~\ref{thm:labeled-stacking-group}.
\end{theorem}
\begin{proof}
The diagonal category and the image of $F_0$ generate their
respective targets under sums and even retracts.
By Propositions~\ref{prop:monoidal-completion}
and~\ref{prop:completion-uniqueness}, $F_0$ extends to the
label-preserving equivalence \eqref{eq:stacking-comparison}.
Reconstruction, followed by reindexing along the chosen labeling,
gives surjectivity of \eqref{eq:labeled-bijection};
Proposition~\ref{prop:gauge-functors} gives injectivity because a
label-preserving equivalence induces $\phi=\id_G$.
The equivalence $F$ identifies multiplication with $\star$,
Lemma~\ref{lem:trivial-pullback} identifies the unit with $e_G$,
and Lemma~\ref{lem:product-gauges} identifies the inverse with $s^\vee$.
\end{proof}

\begin{proposition}
\label{prop:supercohomology}
Let $\mathcal C$ be a Brauer--Wall reference category with an exchange
data and let $G$ be a finite group. The subset
\begin{equation}
\label{eq:supercohomology-group}
 \SH_k(G)=\{t\in\Tw_k(G):\rho_t=0\}/{\sim_G}
\end{equation}
of $\Tcal_k(G)$ is a subgroup, on which the product \eqref{eq:class-operations}
is induced by
\begin{equation}
\label{eq:degree-zero-product}
 (\nu,\omega)\star(\nu',\omega')
 =\bigl(\nu+\nu',\,\omega\omega'(-1)^{\nu\smile_1\nu'}\bigr),
\end{equation}
for every exchange data, with
$[(\nu,\omega)]^{-1}=[(\nu,\omega^{-1}(-1)^{\nu\smile_1\nu})]$.
There is an exact sequence
\begin{equation}
\label{eq:supercohomology-extension}
\begin{gathered}
 0\longrightarrow H^3(G,k^\times)\longrightarrow\SH_k(G)
 \longrightarrow H^2(G,\Z/2\Z)^+\longrightarrow0,\\
 H^2(G,\Z/2\Z)^+=\bigl\{[\nu]\in H^2(G,\Z/2\Z):
   (-1)^{\nu\smile\nu}\in B^4(G,k^\times)\bigr\},
\end{gathered}
\end{equation}
in which the first map is $[\omega]\mapsto[(0,\omega)]$ and the second is
$[(\nu,\omega)]\mapsto[\nu]$.
\end{proposition}
\begin{proof}
When $\rho_s=\rho_t=0$, normalization gives $c_{s,t}=0$, $K_{s,t}=1$,
and $\Xi_{s,t}=\nu_s\smile_1\nu_t$. The product and inverse formulas
follow, and the zero-degree subset is a subgroup.
The map to $[\nu]$ is well defined under gauge transformations and
has image exactly $H^2(G,\Z/2\Z)^+$; independence of representatives
follows from \eqref{eq:gauge-obstruction-coboundary}.
If $\nu=\delta\lambda$, the gauge $(\id_G,\lambda,1)$ replaces it by
$(0,\omega(-1)^{\lambda\smile\nu})$.
For pairs with $\nu=0$, gauge equivalence is exactly
$\omega'=\omega(\delta\beta)^{-1}$, and multiplication is ordinary
multiplication of cocycles. Thus the kernel is $H^3(G,k^\times)$.
\end{proof}

\begin{theorem}
\label{thm:stacking-group}
Let $\mathcal C$ be a Brauer--Wall reference category with an exchange
data and let $G$ be a finite group. The degree map
\begin{equation}
\label{eq:degree}
 d:\Tcal_k(G)\longrightarrow Z^1(G,\BW(k)),\qquad [s]\longmapsto\rho_s,
\end{equation}
is a surjective group homomorphism with $d[p_\rho]=\rho$ and kernel
$\SH_k(G)$, so that there is an exact sequence
\begin{equation}
\label{eq:degree-sequence}
 0\longrightarrow\SH_k(G)\longrightarrow\Tcal_k(G)
 \xrightarrow{\ d\ } Z^1(G,\BW(k))\longrightarrow0.
\end{equation}
\end{theorem}
\begin{proof}
Gauges preserve $\rho$, the product adds it, and $d[p_\rho]=\rho$.
The kernel is $\SH_k(G)$ by definition.
\end{proof}

\section{Reference categories over \texorpdfstring{$k$}{k}}
\label{sec:reference-category}

In this section we define an obstruction
$o_4(k)\in H^4(\BW(k),k^\times)/\operatorname{im}Q$, where
$Q([n])=[(-1)^{n\smile n}]:H^2(\BW(k),\Z/2\Z)\to H^4(\BW(k),k^\times)$, and prove that a Brauer--Wall reference
category over $k$ exists if and only if $o_4(k)=1$. When it exists, we also
classify reference categories up to monoidal superequivalence by
$\SH_k(\BW(k))$.

The strategy of this section is essentially the same as the extension
theory of \cite[Section~8]{EtingofNikshychOstrik}, which classifies
$G$-graded extensions of a fusion category $\mathcal D$ over an
algebraically closed field in terms of its Brauer--Picard group
$\operatorname{BrPic}(\mathcal D)$, the group
$\operatorname{Inv}(\mathcal Z(\mathcal D))$ of invertible objects in
its Drinfeld center, and $k^\times$.
Given a homomorphism $G\to\operatorname{BrPic}(\mathcal D)$, one first
encounters a degree-three obstruction in
$H^3(G,\operatorname{Inv}(\mathcal Z(\mathcal D)))$ to choosing
compatible multiplication equivalences between the graded components;
when it vanishes, the choices form a torsor under
$H^2(G,\operatorname{Inv}(\mathcal Z(\mathcal D)))$.
For each choice, there is then a degree-four pentagon obstruction in
$H^4(G,k^\times)$ to choosing associativity isomorphisms, and when it
vanishes, the choices form a torsor under $H^3(G,k^\times)$.

In our setting, a Brauer--Wall reference category is a $\BW(k)$-graded
extension of $\sVec$ associated with $\id:\BW(k)\to\BW(k)$, and the
roles of $\operatorname{BrPic}(\mathcal D)$,
$\operatorname{Inv}(\mathcal Z(\mathcal D))$, and $k^\times$ are
played by $\BW(k)$, $\Z/2\Z$ (generated by the parity shift $\Pi k$,
with Koszul braiding $-1$), and $k^\times$.
Accordingly, we first show that the degree-three obstruction
$\kappa_3(k)\in H^3(\BW(k),\Z/2\Z)$ to an \emph{even Morita system} always
vanishes, leaving an $H^2(\BW(k),\Z/2\Z)$-torsor of choices; we then
define the pentagon obstruction $O_{D,M}\in H^4(\BW(k),k^\times)$,
which shifts by $Q([n])=[(-1)^{n\smile n}]$ under this $H^2$-action
because of the Koszul sign, yielding $o_4(k)$; and finally we construct
and classify reference categories from \emph{reference systems}.

Although enriched and non-algebraically closed generalizations of
\cite{EtingofNikshychOstrik} have been developed in
\cite{JonesMorrisonPenneysPlavnik,DavydovNikshych,SanfordBrauerPicard},
they do not simultaneously cover $\sVec$-enrichment and
arbitrary non-algebraically closed fields $k$, so they cannot be applied directly here. We therefore
give a self-contained, elementary treatment in terms of central simple
superalgebras and graded bimodules.

\subsection{Sections of braided categorical groups}
\label{subsec:braided-categorical-group-sections}

We recall basic facts on braided categorical groups $\mathcal G$
from \cite[\S3]{JoyalStreet} and describe the obstruction to lifting
$\pi_0\mathcal G$ to a monoidal section. The associator class in
$H^3(\pi_0\mathcal G,\Aut(\one))$ vanishes if and only if the quadratic
function of the braiding is the diagonal of a biadditive map, in which
case the monoidal sections modulo monoidal isomorphism form a torsor
under $H^2(\pi_0\mathcal G,\Aut(\one))$.

\begin{definition}
\label{def:braided-categorical-group}
A \emph{braided categorical group} is a braided monoidal category in
which every morphism is invertible and every object is invertible with
respect to the tensor product \cite[\S3]{JoyalStreet}.
A map $q:\Gamma\to\Lambda$ between abelian groups is \emph{quadratic}
if $q(-a)=q(a)$ for all $a\in\Gamma$ and the map
$\Gamma\times\Gamma\to\Lambda$, $(a,a')\mapsto q(a+a')-q(a)-q(a')$, is
additive in each variable.
\end{definition}

\begin{proposition}
\label{prop:quadratic-function}
Let $\mathcal G$ be a braided categorical group with braiding $c$.
The set $\pi_0\mathcal G$ of isomorphism classes of objects of
$\mathcal G$ is an abelian group with the multiplication induced by
$\otimes$, and $\Aut(\one)$ is abelian.
For every object $Y$ of $\mathcal G$, the map $f\mapsto\id_Y\otimes f$
is an isomorphism $\Aut(\one)\to\Aut(Y\otimes\one)\cong\Aut(Y)$.
The map
\begin{equation}
\label{eq:quadratic-function}
\begin{gathered}
 q:\pi_0\mathcal G\longrightarrow\Aut(\one),\\
 q([X])=\text{the element of }\Aut(\one)\text{ corresponding to }
 c_{X,X}\in\Aut(X\otimes X),
\end{gathered}
\end{equation}
is well defined and quadratic; it is called the \emph{quadratic
function} of $\mathcal G$.
\end{proposition}
\begin{proof}
Tensor inverses and the braiding make $\pi_0\mathcal G$ abelian.
The two unit multiplications give the Eckmann--Hilton argument
for commutativity of $\Aut(\one)$.
Tensoring with an invertible object is an equivalence, giving
$\Aut(\one)\cong\Aut(Y)$ compatibly with conjugation.
Naturality of the braiding therefore makes $q$ well defined.
Writing the group laws additively, the element
$q(a+b)-q(a)-q(b)$ for $a=[X]$, $b=[Y]$ is represented by
$c_{Y,X}c_{X,Y}$. The resulting map in $(a,b)$ is biadditive,
and $q(-a)=q(a)$; see \cite[\S3]{JoyalStreet}.
\end{proof}

\begin{definition}
\label{def:sections}
Let $\mathcal G$ be a braided categorical group with associator $\alpha$, and unitors $l$ and $r$,
and let $(X_a)_{a\in\pi_0\mathcal G}$ be objects of $\mathcal G$ such that
$X_a$ is of class $a$ and $X_0=\one$.
Automorphism groups of objects are identified with $\Aut(\one)$ by
Proposition~\ref{prop:quadratic-function}, and cochains with values in
$\Aut(\one)$ are written multiplicatively.
A \emph{section} of $\mathcal G$ on $(X_a)_{a\in\pi_0\mathcal G}$ is a
family $t=(t_{a,b})_{a,b\in\pi_0\mathcal G}$ of isomorphisms
$t_{a,b}:X_a\otimes X_b\to X_{a+b}$ with $t_{0,a}=l_{X_a}$ and
$t_{a,0}=r_{X_a}$.
Its \emph{associativity cocycle} $h_t\in C^3(\pi_0\mathcal G,\Aut(\one))$
is defined by
\begin{equation}
\label{eq:section-cocycle}
 t_{a,b+c}\circ(\id_{X_a}\otimes t_{b,c})\circ\alpha_{X_a,X_b,X_c}
 =h_t(a,b,c)\circ t_{a+b,c}\circ(t_{a,b}\otimes\id_{X_c}),
\end{equation}
with $h_t(a,b,c)$ acting on $X_{a+b+c}$.
For a normalized $n\in C^2(\pi_0\mathcal G,\Aut(\one))$, the section
$n\cdot t$ is $(n(a,b)\circ t_{a,b})_{a,b}$, with $n(a,b)$ acting on
$X_{a+b}$.
The section $t$ is \emph{monoidal} if $h_t=1$.
\end{definition}

\begin{proposition}
\label{prop:sections-associator-class}
Let $\mathcal G$ and $(X_a)_{a\in\pi_0\mathcal G}$ be as in
Definition~\ref{def:sections}, and let $t$ be a section of
$\mathcal G$ on $(X_a)_{a\in\pi_0\mathcal G}$.

(1) The cochain $h_t$ is a normalized $3$-cocycle, every section on
$(X_a)_{a\in\pi_0\mathcal G}$ is $n\cdot t$ for a unique normalized
$n\in C^2(\pi_0\mathcal G,\Aut(\one))$, and $h_{n\cdot t}=h_t\,\delta n$.
The class
\begin{equation}
\label{eq:associator-class}
 [h_t]\in H^3(\pi_0\mathcal G,\Aut(\one))
\end{equation}
is independent of $t$ and of $(X_a)_{a\in\pi_0\mathcal G}$; it is
called the \emph{associator class} of $\mathcal G$, and braided
equivalent braided categorical groups have the same associator class
\cite[\S3]{JoyalStreet}.

(2) The section $t$ is monoidal if and only if $a\mapsto X_a$ with the
structure isomorphisms $t_{a,b}$ is a monoidal functor
$(\pi_0\mathcal G)_{\mathrm{disc}}\to\mathcal G$, where
$(\pi_0\mathcal G)_{\mathrm{disc}}$ is the strict monoidal category with
object set $\pi_0\mathcal G$ and only identity morphisms, whose tensor
product and unit are the multiplication and the unit of the group
$\pi_0\mathcal G$. Two monoidal sections $t,t'$ are monoidally
isomorphic as such functors if and only if $t'=(\delta\mu)\cdot t$ for
some normalized $\mu\in C^1(\pi_0\mathcal G,\Aut(\one))$.

(3) A monoidal section on $(X_a)_{a\in\pi_0\mathcal G}$ exists if and
only if the associator class of $\mathcal G$ is trivial. In that case
$[n]\cdot t=n\cdot t$ makes the set of monoidal sections on
$(X_a)_{a\in\pi_0\mathcal G}$ modulo monoidal isomorphism a torsor under
$H^2(\pi_0\mathcal G,\Aut(\one))$.

(4) Let $q$ be the quadratic function of $\mathcal G$. The associator
class of $\mathcal G$ is trivial if and only if there is a biadditive
map $\beta:\pi_0\mathcal G\times\pi_0\mathcal G\to\Aut(\one)$ such that 
$q(a)=\beta(a,a)$ for all $a\in\pi_0\mathcal G$.

(5) Let $(X'_a)_{a\in\pi_0\mathcal G}$ be another family of objects such that
$X'_a$ is of class $a$ and $X'_0=\one$, and let $\xi_a:X_a\to X'_a$ be
isomorphisms with $\xi_0=\id_{\one}$. Then
\begin{equation}
\label{eq:section-transport}
 t\longmapsto t^\xi,\qquad
 t^\xi_{a,b}=\xi_{a+b}\circ t_{a,b}\circ(\xi_a\otimes\xi_b)^{-1},
\end{equation}
is a bijection from the sections on $(X_a)_{a\in\pi_0\mathcal G}$ to the
sections on $(X'_a)_{a\in\pi_0\mathcal G}$ with $h_{t^\xi}=h_t$ and
$(n\cdot t)^\xi=n\cdot t^\xi$; it induces a bijection between the
torsors of (3) which does not depend on the $\xi_a$.
\end{proposition}
\begin{proof}
(1) The unit axioms normalize $h_t$, and the pentagon gives
\[
 h_t(b,c,d)h_t(a,b+c,d)h_t(a,b,c)
 =h_t(a+b,c,d)h_t(a,b,c+d).
\]
Any other section is uniquely $n\cdot t$; substitution gives
$h_{n\cdot t}=h_t\delta n$. Transport along representatives
preserves $h_t$, and invariance under braided equivalence follows
from \cite[\S3]{JoyalStreet}.

(2)--(3) The monoidal-functor equation is $h_t=1$.
A normalized natural isomorphism $\mu$ changes the section by
$t'=(\delta\mu^{-1})\cdot t$. Thus a monoidal section exists exactly
when $[h_t]=1$, and their classes form the asserted $H^2$-torsor.

(4) A monoidal section transports the braiding to a biadditive
map with diagonal $q$. Conversely, a biadditive $\beta$ with this
diagonal defines a skeletal braided categorical group with trivial
associator and quadratic function $q$. The classification in
\cite[\S3]{JoyalStreet} gives a braided equivalence with $\mathcal G$.

(5) Conjugation by $\xi$ gives the stated bijection and commutes
with the cochain action. Replacing $\xi$ by $\xi u$ changes the
transported section by $\delta u^{-1}$, so the induced torsor
bijection is independent of the chosen isomorphisms.
\end{proof}

\subsection{The Morita groupoid and even Morita systems}
\label{subsec:morita-groupoid-systems}

We apply the preceding framework to the Morita groupoid $\Mor$ of
central simple superalgebras over $k$. Identifying its quadratic
function with the parity $\epsilon:\BW(k)\to\Z/2\Z$ (defined in \eqref{eq:bw-parity}) shows
that the degree-three obstruction $\kappa_3(k)\in H^3(\BW(k),\Z/2\Z)$
always vanishes. Translating monoidal sections of $\Mor$ into bimodule
data, we conclude that even Morita systems exist and form a torsor
under $H^2(\BW(k),\Z/2\Z)$ up to isomorphism and parity transport.

\begin{definition}
\label{def:morita-groupoid}
The \emph{Morita groupoid} $\Mor$ is the category whose objects are the
finite-dimensional central simple superalgebras and whose morphisms
$A\to B$ are the even isomorphism classes $[M]$ of Morita equivalences
$M:A\to B$, with the composite of $[M]:A\to B$ and $[N]:B\to C$ equal
to $[M\otimes_BN]$ and with identities $[A]$.
It is equipped with the product $\gt$ on objects and on classes and
with the unit $k$. An isomorphism of superalgebras is regarded as a
Morita equivalence by Definition~\ref{def:graded-bimodules}.
The associator is the class of the isomorphism
$(A\gt B)\gt C\to A\gt(B\gt C)$, $(a\gt b)\gt c\mapsto a\gt(b\gt c)$,
the unitors are the classes of the isomorphisms $k\gt A\to A$,
$\lambda\gt a\mapsto\lambda a$, and $A\gt k\to A$,
$a\gt\lambda\mapsto a\lambda$, and the braiding is the class of the
swap $\tau_{A,B}:A\gt B\to B\gt A$, $a\gt b\mapsto(-1)^{|a||b|}b\gt a$.
\end{definition}

\begin{lemma}
\label{lem:truncated-picard}
The category $\Mor$ with the structure of
Definition~\ref{def:morita-groupoid} is a braided categorical group.
By Definitions~\ref{def:graded-bimodules} and
\ref{def:brauer-wall-group}, the isomorphism classes of
objects of $\Mor$ are the graded Morita classes, so $\pi_0\Mor=\BW(k)$,
the multiplication induced by $\gt$ being the addition of
Definition~\ref{def:brauer-wall-group}.
For every object $X$, the automorphism group $\Aut(X)$ of $X$ in $\Mor$
consists of the classes of the two Morita equivalences ${}_XX_X$ and
$\Pi({}_XX_X)$, where $\Pi({}_XX_X)$ has the actions
\eqref{eq:morita-parity-shift}, and the isomorphism $\Aut(k)\to\Aut(X)$
of Proposition~\ref{prop:quadratic-function} sends the class of $\Pi k$
to the class of $\Pi({}_XX_X)$.
We identify $\Aut(k)$ with $\Z/2\Z$ by sending the class of $\Pi k$ to
$1$. Its quadratic function is then
\begin{equation}
\label{eq:truncated-picard-quadratic-function}
 q(a)=\epsilon(a),
\end{equation}
where $\epsilon:\BW(k)\to\Z/2\Z$ is the parity \eqref{eq:bw-parity}.
\end{lemma}
\begin{proof}
Balanced tensor products give composition and inverses of Morita
classes. The even interchange isomorphism makes $\gt$ functorial,
and the graded swaps give a braiding; their coherence descends
from super vector spaces. Each object is invertible because
$A\gt A^{\mathrm{sop}}\cong\End_k(A)$ is Morita equivalent to $k$.
Hence $\pi_0\Mor=\BW(k)$.
Proposition~\ref{prop:morita-equivalences} gives the two
automorphisms $[A]$ and $[\Pi A]$, with
\[
 A\gt\Pi k\xrightarrow{\sim}\Pi A,
 \qquad a\gt\Pi\lambda\longmapsto(-1)^{|a|}\Pi(a\lambda).
\]
For $B=D_a\gt D_a$, Proposition~\ref{prop:swap-parity} identifies
the swap bimodule with $\Pi^{\epsilon(a)}B$ up to even isomorphism.
Thus $q(a)=\epsilon(a)$.
\end{proof}

\begin{lemma}
\label{lem:parity-obstruction-vanishes}
The associator class $\kappa_3(k)$ of $\Mor$, read in $H^3(\BW(k),\Z/2\Z)$ through
the identification $\Aut(k)\cong\Z/2\Z$ of Lemma~\ref{lem:truncated-picard}, is trivial.
\end{lemma}
\begin{proof}
The biadditive map $\beta(a,b)=\epsilon(a)\epsilon(b) \in \mathbb F_2$ has diagonal
$\beta(a,a)=\epsilon(a)=q(a)$. Apply
Proposition~\ref{prop:sections-associator-class}(4).
\end{proof}

\begin{definition}
\label{def:morita-system}
Let $D=(D_a)_{a\in\BW(k)}$ be a family of representatives as in
Definition~\ref{def:brauer-wall-group}.
A \emph{Morita system} is a pair $(D,M)$ where $M$ is a family
$M=(M_{a,b})_{a,b\in\BW(k)}$ of Morita equivalences
$M_{a,b}:D_a\gt D_b\to D_{a+b}$ such that
$M_{0,a}=D_a=M_{a,0}$, where $k\gt D_a$ and $D_a\gt k$ are identified
with $D_a$ by the unit isomorphisms.
For $a,b,c\in\BW(k)$ put
\begin{equation}
\label{eq:morita-triple-composites}
\begin{aligned}
 L_{a,b,c}
 &=(M_{a,b}\gt D_c)\otimes_{D_{a+b}\gt D_c}M_{a+b,c},\\
 U_{a,b,c}
 &=(D_a\gt M_{b,c})\otimes_{D_a\gt D_{b+c}}M_{a,b+c},
\end{aligned}
\end{equation}
both Morita equivalences $(D_a\gt D_b)\gt D_c\to D_{a+b+c}$, the left
action on $U_{a,b,c}$ being transported along the associativity
isomorphism of $\gt$.
By Proposition~\ref{prop:morita-equivalences}, we define the \emph{parity cochain} $p_{D,M}\in C^3(\BW(k),\Z/2\Z)$ of $(D,M)$ by
\begin{equation}
\label{eq:morita-parity-cochain}
 p_{D,M}(a,b,c)=|A_{a,b,c}|
\end{equation}
for any nonzero homogeneous bimodule map
$A_{a,b,c}:L_{a,b,c}\to U_{a,b,c}$.
The Morita system $(D,M)$ is \emph{even} if $p_{D,M}=0$.
When one of $a,b,c$ is $0$, the normalization
$M_{0,a}=D_a=M_{a,0}$ and the unit isomorphisms
\eqref{eq:morita-left-unitor}, \eqref{eq:morita-right-unitor},
$k\gt N\cong N$, and $N\gt k\cong N$ give canonical even
bimodule isomorphisms
\[
 \lambda_L:L_{a,b,c}\xrightarrow{\sim}N,
 \qquad
 \lambda_U:U_{a,b,c}\xrightarrow{\sim}N,
\]
where $N=M_{b,c}$, $M_{a,c}$, or $M_{a,b}$ according as
$a=0$, $b=0$, or $c=0$.
We call the composite
\[
 L_{a,b,c}\xrightarrow{\lambda_L}N
 \xrightarrow{\lambda_U^{-1}}U_{a,b,c}
\]
the \emph{unit identification}.
An \emph{even system} is a triple $(D,M,A)$ where $(D,M)$ is an even Morita system and $A$ is a family $A = (A_{a,b,c})_{a,b,c\in \BW(k)}$ of even bimodule isomorphisms $A_{a,b,c}:L_{a,b,c}\to U_{a,b,c}$ such that
$A_{a,b,c}$ is the unit identification whenever one of $a,b,c$ is $0$.
For a normalized $n\in C^2(\BW(k),\Z/2\Z)$, put $\Pi^nM=(\Pi^{n(a,b)}M_{a,b})_{a,b}$, with the actions \eqref{eq:morita-parity-shift}.
Then $(D,\Pi^nM)$ is a Morita system.
An \emph{isomorphism} of Morita systems $M\to M'$ is a family of even
bimodule isomorphisms $\phi_{a,b}:M_{a,b}\to M'_{a,b}$, $a,b\in\BW(k)$;
two Morita systems are \emph{isomorphic} if an isomorphism between them
exists, and $M'$ is a \emph{parity transport} of $M$ if
$M'$ is isomorphic to $\Pi^{\delta\mu}M$ for some normalized
$\mu\in C^1(\BW(k),\Z/2\Z)$.
\end{definition}

\begin{lemma}
\label{lem:morita-systems-sections}
Let $D=(D_a)_{a\in\BW(k)}$ be the family of representatives of
Definition~\ref{def:brauer-wall-group}, regarded as objects of $\Mor$,
identify $\Aut(k)$ with $\Z/2\Z$ as in Lemma~\ref{lem:truncated-picard},
and let $\Sec(D)$ be the set of even Morita systems with respect
to $D$ modulo isomorphism and parity transport.
For a Morita system $M$ with respect to $D$, the family
$t_M=([M_{a,b}])_{a,b\in\BW(k)}$ is a section of $\Mor$ on
$D$, and $M\mapsto t_M$ is a bijection from the set of
Morita systems modulo isomorphism onto the set of sections of $\Mor$ on
$D$, such that
\begin{equation}
\label{eq:morita-systems-sections}
 h_{t_M}=p_M,\qquad t_{\Pi^nM}=n\cdot t_M
 \qquad(n\in C^2(\BW(k),\Z/2\Z)\text{ normalized}).
\end{equation}
Consequently $M$ is even if and only if $t_M$ is monoidal, parity
transports of $M$ correspond to sections monoidally isomorphic to
$t_M$, even Morita systems exist, and $M\mapsto t_M$ identifies
$\Sec(D)$ with the set of monoidal sections of $\Mor$ on
$D$ modulo monoidal isomorphism; in particular
$\Sec(D)$ is a torsor under $H^2(\BW(k),\Z/2\Z)$, acting by
$[n]\cdot M=\Pi^nM$.
\end{lemma}
\begin{proof}
A section is exactly a choice of the even isomorphism classes
$[M_{a,b}]$. Its two triple composites are $[L_{a,b,c}]$ and
$[U_{a,b,c}]$, and
$U_{a,b,c}\cong\Pi^{p_M(a,b,c)}L_{a,b,c}$ evenly. Thus
\[
 h_{t_M}=p_M,\qquad t_{\Pi^nM}=n\cdot t_M.
\]
Even systems therefore correspond to monoidal sections, and parity
transports to monoidal isomorphisms. Existence and the $H^2$-torsor
follow from Lemma~\ref{lem:parity-obstruction-vanishes} and
Proposition~\ref{prop:sections-associator-class}.
\end{proof}

\subsection{The pentagon obstruction}
\label{subsec:pentagon-obstruction}

For an even system $(D,M,A)$, the two composites around the pentagon
for four Brauer--Wall classes differ by a scalar, defining a normalized
$4$-cocycle $O_{D,M,A}\in Z^4(\BW(k),k^\times)$. Its cohomology class
$O_{D,M}\in H^4(\BW(k),k^\times)$ is independent of $A$ and shifts by
$Q([n])=[(-1)^{n\smile n}]$ under the $H^2(\BW(k),\Z/2\Z)$-action on
$M$; the resulting coset $o_4(k)=O_{D,M}\operatorname{im}Q$ is the
exact obstruction to the existence of a reference system.

\begin{definition}
\label{def:pentagon-scalar}
Let $(M,A)$ be an even system.
A \emph{word} is either a letter $a\in\BW(k)$ or $(w_1w_2)$ for words
$w_1,w_2$.
Define $|w|\in\BW(k)$, a superalgebra $D_w$, and a Morita equivalence
$M_w:D_w\to D_{|w|}$ by recursion on $w$, starting from $D$:
\begin{equation}
\label{eq:iterated-bimodule}
\begin{gathered}
 |a|=a,\qquad |(w_1w_2)|=|w_1|+|w_2|,\qquad
 D_{(w_1w_2)}=D_{w_1}\gt D_{w_2},\\
 M_a=D_a,\qquad
 M_{(w_1w_2)}=(M_{w_1}\gt M_{w_2})\otimes_{D_{|w_1|}\gt D_{|w_2|}}M_{|w_1|,|w_2|}.
\end{gathered}
\end{equation}
For words $w_1,w_2,w_3$, let
\begin{equation}
\label{eq:iterated-identifications}
\begin{aligned}
 \iota^L_{w_1,w_2,w_3}&=\mathsf a_{(M_{w_1}\gt M_{w_2})\gt M_{w_3},\,
   M_{|w_1|,|w_2|}\gt D_{|w_3|},\,M_{|w_1|+|w_2|,|w_3|}}\\
 &\qquad\circ\Bigl(\bigl(
  \chi^{-1}_{M_{w_1}\gt M_{w_2},\,M_{w_3},\,M_{|w_1|,|w_2|},\,D_{|w_3|}}
  \circ(\id\gt\mathsf r_{M_{w_3}}^{-1})\bigr)\otimes\id\Bigr)\\
 &:M_{((w_1w_2)w_3)}\longrightarrow
 \bigl((M_{w_1}\gt M_{w_2})\gt M_{w_3}\bigr)
 \otimes_{D_{|w_1|}\gt D_{|w_2|}\gt D_{|w_3|}}L_{|w_1|,|w_2|,|w_3|},\\
 \iota^U_{w_1,w_2,w_3}&=\mathsf a_{M_{w_1}\gt(M_{w_2}\gt M_{w_3}),\,
   D_{|w_1|}\gt M_{|w_2|,|w_3|},\,M_{|w_1|,|w_2|+|w_3|}}\\
 &\qquad\circ\Bigl(\bigl(
  \chi^{-1}_{M_{w_1},\,M_{w_2}\gt M_{w_3},\,D_{|w_1|},\,M_{|w_2|,|w_3|}}
  \circ(\mathsf r_{M_{w_1}}^{-1}\gt\id)\bigr)\otimes\id\Bigr)\\
 &:M_{(w_1(w_2w_3))}\longrightarrow
 \bigl(M_{w_1}\gt(M_{w_2}\gt M_{w_3})\bigr)
 \otimes_{D_{|w_1|}\gt D_{|w_2|}\gt D_{|w_3|}}U_{|w_1|,|w_2|,|w_3|},
\end{aligned}
\end{equation}
with $\mathsf a$, $\mathsf r$, and $\chi$ the isomorphisms
\eqref{eq:morita-associator}, \eqref{eq:morita-right-unitor}, and
\eqref{eq:morita-interchange}, and with
$M_{w_1}\gt(M_{w_2}\gt M_{w_3})$ identified with
$(M_{w_1}\gt M_{w_2})\gt M_{w_3}$ along the associativity isomorphism of
$\gt$.
Words $w,w'$ \emph{differ by one reassociation} if $w'$ is obtained
from $w$ by replacing one subword of the form $((w_1w_2)w_3)$ by
$(w_1(w_2w_3))$. For such $w,w'$, define the even isomorphism
$c_{w,w'}:M_w\to M_{w'}$ by recursion on $w$, according to whether
the replaced subword is $w$ itself, lies in the left factor of $w$, or
lies in the right factor of $w$:
\begin{equation}
\label{eq:iterated-comparison}
\begin{aligned}
 c_{((w_1w_2)w_3),(w_1(w_2w_3))}
 &=(\iota^U_{w_1,w_2,w_3})^{-1}\circ(\id\otimes A_{|w_1|,|w_2|,|w_3|})
   \circ\iota^L_{w_1,w_2,w_3},\\
 c_{(w_1w_2),(w_1'w_2)}&=(c_{w_1,w_1'}\gt\id_{M_{w_2}})\otimes\id,\\
 c_{(w_1w_2),(w_1w_2')}&=(\id_{M_{w_1}}\gt c_{w_2,w_2'})\otimes\id.
\end{aligned}
\end{equation}
For $a,b,c,d\in\BW(k)$, the even isomorphisms
\begin{equation}
\label{eq:reference-pentagon-scalar}
\begin{gathered}
 P_L=c_{(a((bc)d)),(a(b(cd)))}\circ c_{((a(bc))d),(a((bc)d))}\circ
     c_{(((ab)c)d),((a(bc))d)},\\
 P_R=c_{((ab)(cd)),(a(b(cd)))}\circ c_{(((ab)c)d),((ab)(cd))}
\end{gathered}
\end{equation}
from $M_{(((ab)c)d)}$ to $M_{(a(b(cd)))}$ differ by a unique scalar
(Proposition~\ref{prop:morita-equivalences}),
\begin{equation}
\label{eq:pentagon-scalar}
 P_L=O_{M,A}(a,b,c,d)\,P_R,\qquad O_{M,A}(a,b,c,d)\in k^\times .
\end{equation}
The function $O_{M,A}:\BW(k)^4\to k^\times$ is the \emph{pentagon
scalar} of $(M,A)$, and an even system with $O_{M,A}=1$ is called a
\emph{reference system}.
\end{definition}

\begin{lemma}
\label{lem:pentagon-cocycle}
The pentagon scalar $O_{D,M,A}$ of an even system is a normalized
$4$-cocycle.
\end{lemma}
\begin{proof}
If an argument is $0$, the unit identifications give $P_L=P_R$, so
$O=O_{D,M,A}$ is normalized. On the associahedron for five letters,
the square faces commute by naturality of the even comparisons, and the
six pentagonal faces give the factors of $\delta O$. Each edge occurs
twice with opposite orientations, so $\delta O=1$.
\end{proof}

\begin{lemma}
\label{lem:representatives-independence}
Let $(D,M,A)$ be an even system and let $D'$ be another family of
representatives. There exists an even system $(D',M',A')$ such that
\begin{equation}
\label{eq:representatives-pentagon}
 O_{D,M,A}=O_{D',M',A'}.
\end{equation}
\end{lemma}
\begin{proof}
Choose isomorphisms $\phi_a:D_a\to D'_a$ with $\phi_0=\id_k$ by
Proposition~\ref{prop:central-simple-superalgebras}(1)--(2).
Transport the bimodule actions and comparisons along the $\phi_a$, using
$M'_{0,a}=M'_{a,0}=D'_a$ at the unit. This gives an even system, and
the induced isomorphisms on the iterated bimodules conjugate both
pentagon composites, preserving their ratio.
\end{proof}

\begin{lemma}
\label{lem:comparison-change}
Let $(D,M,A)$ be an even system.
For every even system $(D,M,A')$, there exists
$\theta\in C^3(\BW(k),k^\times)$ such that
\begin{equation}
\label{eq:pentagon-cohomology}
 O_{D,M,A'}=O_{D,M,A}\,\delta\theta.
\end{equation}
Conversely, for every $\theta\in C^3(\BW(k),k^\times)$, there exists an
even system $(D,M,A')$ satisfying \eqref{eq:pentagon-cohomology}.
\end{lemma}
\begin{proof}
By Proposition~\ref{prop:morita-equivalences}, each comparison space is
one-dimensional, so $A'=\theta A$ for a unique $3$-cochain $\theta$.
The unit condition is precisely the normalization of $\theta$.
In \eqref{eq:reference-pentagon-scalar}, the three factors of $\theta$
on $P_L$ divided by the two on $P_R$ give $(\delta\theta)(a,b,c,d)$.
\end{proof}

\begin{definition}
  Let $(D,M)$ be an even Morita system. 
  Choose an even system $(D,M,A)$ and define
\begin{equation}
\label{eq:pentagon-cohomology-class}
 O_{D,M}:=[O_{D,M,A}]\in H^4(\BW(k),k^\times),
\end{equation}
which does not depend on the choice of $A$.
\end{definition}

\begin{lemma}
\label{lem:pentagon-shift}
Let $(D,M)$ be an even Morita system.
For every even Morita system $(D,M')$, there exists
$n\in Z^2(\BW(k),\Z/2\Z)$ such that
\begin{equation}
\label{eq:pentagon-sign}
 O_{D,M'}=O_{D,M}[(-1)^{n\smile n}].
\end{equation}
Conversely, for every $n\in Z^2(\BW(k),\Z/2\Z)$, there exists an even
Morita system $(D,M')$ satisfying \eqref{eq:pentagon-sign}.
\end{lemma}
\begin{proof}
By Proposition~\ref{prop:morita-equivalences}, every $M'$ is evenly
isomorphic to $\Pi^nM$ for a normalized $2$-cochain $n$. The shifts on
$L_{a,b,c}$ and $U_{a,b,c}$ differ in degree by $(\delta n)(a,b,c)$,
so $\Pi^nM$ is even exactly when $n$ is a cocycle.
For such $n$, transport a choice of $A$ to $A^n$ by tensoring the shift
maps in the order of the factors in \eqref{eq:morita-triple-composites}.
After cancellation along the two pentagon paths, the only interchange
is between the shifts on $(a,b)$ and $(c,d)$ at $((ab)(cd))$. Its
Koszul sign gives
\[
 O_{D,\Pi^nM,A^n}(a,b,c,d)
 =O_{D,M,A}(a,b,c,d)\,(-1)^{n(a,b)n(c,d)}.
\]
Even isomorphisms preserve the pentagon ratio, so passing to cohomology
gives \eqref{eq:pentagon-sign} in both directions.
\end{proof}

\begin{lemma}
\label{lem:pentagon-gauge-classes}
The map
\begin{equation}
\label{eq:reference-secondary-obstruction}
 Q:H^2(\BW(k),\Z/2\Z)\longrightarrow H^4(\BW(k),k^\times),
 \qquad Q([n])=[(-1)^{n\smile n}],
\end{equation}
is a well-defined group homomorphism.
\end{lemma}
\begin{proof}
The cup product descends to cohomology by \eqref{eq:cochain-leibniz}.
By \eqref{eq:stacking-cup-one-identities} the two mixed terms in the
square of a sum cancel in mod-$2$ cohomology. Thus the cup square is
additive, and composing it with the coefficient map $r\mapsto(-1)^r$
gives $Q$.
\end{proof}

\begin{definition}
\label{def:reference-obstruction}
Choose an even Morita system $(D,M)$ and define
\begin{equation}
\label{eq:reference-obstruction-class}
 o_4(k)=O_{D,M}\operatorname{im}Q
 \in H^4(\BW(k),k^\times)/\operatorname{im}Q,
\end{equation}
which does not depend on the choice. 
\end{definition}

\begin{proposition}
\label{prop:scalar-obstruction}
A reference system exists if and only if $o_4(k) = 1$. 
\end{proposition}
\begin{proof}
Even Morita systems exist by Lemma~\ref{lem:morita-systems-sections},
and Lemmas~\ref{lem:representatives-independence} and~\ref{lem:pentagon-shift}
make the coset $O_{D,M}\operatorname{im}Q$ independent of the choice.
If the coset is trivial, choose a cocycle $n$ with $O_{D,M}Q([n])=1$. Then
$\Pi^nM$ has trivial pentagon class, so a rescaling as in
Lemma~\ref{lem:comparison-change} makes its pentagon scalar $1$.
Conversely, a reference system has scalar $1$ and represents the trivial coset.
\end{proof}

\subsection{Construction and classification of reference categories}
\label{subsec:reference-category-classification}

We connect reference systems to Brauer--Wall reference categories:
a reference system $(D,M,A)$ produces a reference category of graded
right $D_a$-modules, and conversely every reference category yields a
reference system, showing that a reference category exists if and only
if $o_4(k)=1$. When one exists, twisting it by elements of
$\SH_k(\BW(k))$ gives a bijection onto the monoidal superequivalence
classes of Brauer--Wall reference categories over $k$.

\begin{definition}
\label{def:reference-modules}
Let $(D,M,A)$ be a reference system.
We define a monoidal $\Pi$-supercategory $\mathcal R$ as follows. 
The objects are
the families $X=(X_a)_{a\in\BW(k)}$ of
finite-dimensional graded right $D_a$-modules $X_a$ such that $X_a=0$ for all but finitely many $a$.
The morphism spaces are
\begin{equation}
\label{eq:reference-module-homs}
 \mathcal R(X,Y)
 =\bigoplus_{a\in\BW(k)}
   \Hom_{D_a}(X_a,Y_a),
\end{equation}
with componentwise composition.
The tensor product of $X,Y\in\mathcal R$ is the family
\begin{equation}
\label{eq:reference-module-product}
 (X\otimes Y)_c
 =\bigoplus_{a+b=c}(X_a\otimes_kY_b)\otimes_{D_a\gt D_b}M_{a,b}.
\end{equation}
The tensor product of $f=(f_a)_{a}\in\mathcal R(X,X')_d$ and
$g=(g_b)_{b}\in\mathcal R(Y,Y')_{d'}$ is the element
$f\otimes g\in\mathcal R(X\otimes Y,X'\otimes Y')_{d+d'}$ with
\begin{equation}
\label{eq:reference-module-product-maps}
 (f\otimes g)\bigl((x\otimes y)\otimes_{D_a\gt D_b}m\bigr)
 =(-1)^{d'|x|}\bigl(f_a(x)\otimes g_b(y)\bigr)\otimes_{D_a\gt D_b}m
\end{equation}
for homogeneous $x\in X_a$, $y\in Y_b$, and $m\in M_{a,b}$.
The unit object $\one$ is the family supported at $0$ with $\one_0=k$; the parity object $\pi$ is the family
supported at $0$ with $\pi_0=\Pi k$; and $\zeta\in\mathcal R(\pi,\one)_1$
is the odd isomorphism with $\zeta_0(\Pi\lambda)=\lambda$ for
$\lambda\in k$.
For $X,Y,Z$ supported at $a,b,c\in\BW(k)$, respectively, put
$s=a+b+c$, $E=(D_a\gt D_b)\gt D_c$, and
$V=(X_a\otimes_k Y_b)\otimes_k Z_c$ with its natural graded right
$E$-action. The component at $s$ of the associator
$\alpha_{X,Y,Z}:(X\otimes Y)\otimes Z\to X\otimes(Y\otimes Z)$
is the composite
\begin{equation}
\label{eq:reference-module-associator}
\begin{aligned}
 \bigl((X\otimes Y)\otimes Z\bigr)_s
 &\xrightarrow{\sim} V\otimes_E L_{a,b,c}\\
 &\xrightarrow{\ \id_V\otimes_E A_{a,b,c}\ } V\otimes_E U_{a,b,c}\\
 &\xrightarrow{\sim} \bigl(X\otimes(Y\otimes Z)\bigr)_s .
\end{aligned}
\end{equation}
The outer arrows are the canonical identifications constructed as in
\eqref{eq:iterated-identifications}, using graded interchange to move
$M_{a,b}$ past $Z_c$ in the first arrow.
For arbitrary $X,Y,Z$, the component at $s$ is the direct sum of these
maps over all $a+b+c=s$.
The unitors $l_X:\one\otimes X\to X$ and $r_X:X\otimes\one\to X$ are the
even isomorphisms with
\begin{equation}
\label{eq:reference-module-unitors}
 l_X\bigl((\lambda\otimes x)\otimes_{k\gt D_a}u\bigr)=\lambda xu,\qquad
 r_X\bigl((x\otimes\lambda)\otimes_{D_a\gt k}u\bigr)=\lambda xu
\end{equation}
for $\lambda\in k$, $x\in X_a$, and $u\in D_a=M_{0,a}=M_{a,0}$.
For $a\in\BW(k)$, let $S_a\in\mathcal R$ be the family supported at $a$
with $(S_a)_a=D_a$, the regular right $D_a$-module.
\end{definition}

\begin{proposition}
  \label{prop:construction-of-reference-category}
  Let $(D,M,A)$ be a reference system.
Then $\mathcal R$ is a
Brauer--Wall reference category.
\end{proposition}
\begin{proof}
In degrees $a,b$, tensoring is the external tensor product followed by
the Morita equivalence $M_{a,b}$. Under the balanced tensor
identifications the associator is $\id\otimes A_{a,b,c}$, hence even
and supernatural. The pentagon follows from $P_L=P_R$ and the triangle
from unit normalization. These identities extend by direct sums.

Graded freeness over $D_a$ gives semisimplicity. In degree $a$ every
simple is evenly isomorphic to $S_a$ or $\Pi S_a$, and left multiplication
identifies $\End(S_a)$ with $D_a$. Thus there is one homogeneous
isomorphism class of simples of each type, and $\End(\one)=k$.
For $X,X'$ in degree $a$ and $Y,Y'$ in degree $b$, graded freeness
and the Morita equivalence $M_{a,b}$ give
\begin{equation}
\label{eq:reference-module-tensor-homs}
 \mathcal R(X,X')\otimes_k\mathcal R(Y,Y')
 \xrightarrow{\sim}\mathcal R(X\otimes Y,X'\otimes Y').
\end{equation}

For simple $X$ in degree $a$, choose $Y=S_{-a}$ or $\Pi S_{-a}$ with
a nonzero even $\ev:Y\otimes X\to\one$. Since every nonzero homogeneous
endomorphism of $X$ is invertible, $g\mapsto\ev\circ(\id_Y\otimes g)$
embeds $\End(X)$ into $\mathcal R(Y\otimes X,\one)$. Both spaces have
dimension $\dim_kD_a$ by \eqref{eq:reference-module-tensor-homs}, since
$Y\otimes X$ has degree $0$ and $D_{-a}\cong D_a^{\mathrm{sop}}$.
Thus evaluation induces
\[
 \mathcal R(T,X\otimes V)\xrightarrow{\sim}\mathcal R(Y\otimes T,V)
\]
first for $T=X\otimes V'$ by \eqref{eq:reference-module-tensor-homs},
then for all $T$ since such products generate each component under sums
and even summands. Hence $Y$ is a left dual of $X$. Right duals follow
symmetrically, and duals extend over finite sums. The nonzero product
$X\otimes Y$ in degree $0$ is a sum of copies of $\one$ and $\pi$, so
$X$ is weakly invertible.
\end{proof}

\begin{proposition}
\label{thm:reference-existence}
Let $\mathcal C$ be a Brauer--Wall reference category.
We define a reference system $(D,M,A)$ as follows. 
Choose simple objects $S_a\in\mathcal C_a$ with $S_0=\one$ and superalgebra
isomorphisms $D_a\cong\End(S_a)$. Put
\begin{equation}
\label{eq:reference-extracted-bimodules}
 M_{a,b}=\mathcal C(S_{a+b},S_a\otimes S_b)\quad(a,b\ne0),\qquad
 M_{0,a}=M_{a,0}=D_a,
\end{equation}
with the left action of $D_a\gt D_b\cong\End(S_a\otimes S_b)$ on
$\mathcal C(S_{a+b},S_a\otimes S_b)$ by postcomposition and the right
action of $D_{a+b}$ by precomposition; here $D_a$ replaces the bimodules
$\mathcal C(S_a,\one\otimes S_a)$ and $\mathcal C(S_a,S_a\otimes\one)$
along the even isomorphisms
\begin{equation}
\label{eq:reference-unit-replacements}
\begin{aligned}
 \mathcal C(S_a,\one\otimes S_a)&\longrightarrow\End(S_a)\cong D_a,&
 f&\longmapsto l_{S_a}\circ f,\\
 \mathcal C(S_a,S_a\otimes\one)&\longrightarrow\End(S_a)\cong D_a,&
 f&\longmapsto r_{S_a}\circ f.
\end{aligned}
\end{equation}
Then $(D,M)$ is a Morita system; the associator of $\mathcal C$, transported
along \eqref{eq:reference-unit-replacements} at the positions where an
argument is $0$, induces comparisons $A_{a,b,c}$, and $(D,M,A)$ is a
reference system.
\end{proposition}
\begin{proof}
By \eqref{eq:component-module-equivalence} and
Proposition~\ref{prop:simple-tensor}, the left action identifies
$D_a\gt D_b$ with $\End_{D_{a+b}}(M_{a,b})$. Since $M_{a,b}$ is
nonzero and graded free over $D_{a+b}$, it is a Morita equivalence.
For a word $w$, let $S_w$ denote the corresponding bracketed tensor
product. By Lemma~\ref{lem:component-tensor-homs} and
\eqref{eq:component-module-equivalence}, composition gives
\[
 M_w\cong\mathcal C(S_{|w|},S_w),
\]
using \eqref{eq:reference-unit-replacements} at the unit.
Under these identifications, postcomposition with $\alpha_{S_a,S_b,S_c}$
defines $A_{a,b,c}$. The triangle gives the prescribed unit comparisons,
and the pentagon in $\mathcal C$ gives $P_L=P_R$.
\end{proof}

\begin{theorem}
  \label{thm:existence-of-reference}
  The following are equivalent:

  (1) A Brauer--Wall reference category over $k$ exists,

  (2) A reference system over $k$ exists,

  (3) $o_4(k)=1 \in H^4(\BW(k),k^\times)/\operatorname{im}Q$.
\end{theorem}
\begin{proof}
Propositions~\ref{prop:construction-of-reference-category}
and~\ref{thm:reference-existence} give $(1)\Leftrightarrow(2)$.
Proposition~\ref{prop:scalar-obstruction} gives $(2)\Leftrightarrow(3)$.
\end{proof}

\begin{proposition}
\label{prop:reference-torsor}
Let $\mathcal C$ be a Brauer--Wall reference category over $k$.
For normalized $n\in Z^2(\BW(k),\Z/2\Z)$ and $\theta\in C^3(\BW(k),k^\times)$
with $\delta\theta=(-1)^{n\smile n}$, the pullback
$\mathcal P_{\mathcal C}(\BW(k),\id,n,\theta)$ is a Brauer--Wall reference
category, and
\begin{equation}
\label{eq:reference-torsor}
 [(n,\theta)]\longmapsto[\mathcal P_{\mathcal C}(\BW(k),\id,n,\theta)]
\end{equation}
is a bijection from $\SH_k(\BW(k))$ onto the set of Brauer--Wall
reference categories up to monoidal superequivalence. 
Here, for infinite groups $G$, we use the same cochain formulas for
twists, gauges, and $\SH_k(G)$, and form
the pullback using finitely supported families.
\end{proposition}
\begin{proof}
The underlying supercategory is
$\bigoplus_{a\in\BW(k)}\mathcal C_a$, with split unit and one homogeneous
simple class of each type. The pentagon and triangle checks in the
proof of Proposition~\ref{prop:parity-pentagon} extend from homogeneous
summands by finite direct sums. The left dual constructed in the proof
of Proposition~\ref{prop:pullback-laws} sends a simple of degree $a$ to
a simple of degree $-a$. These two degrees contain at most four even
isomorphism classes of simples. Left duality maps this finite set
injectively to itself, hence bijectively. Thus every simple has a right
dual, and duals extend over finite sums.
A simple tensored with its left dual is nonzero of degree $0$, hence a
sum of copies of the unit and its parity shift. This gives weak
invertibility and proves that the pullback is a reference category.

Index any other reference category $\mathcal C'$ by type using
Proposition~\ref{prop:intrinsic-group}. Choose component equivalences
$\mathcal C'_a\to\mathcal C_a$ and normalized binary comparisons as in
the proof of Proposition~\ref{prop:reconstruction}. Their degrees
$n(a,b)$ and the triple comparison scalars $\theta(a,b,c)$ satisfy
$\delta n=0$ and $\delta\theta=(-1)^{n\smile n}$ by the triple and
pentagon identities. The component equivalences assemble on finite
direct sums as in the proof of Proposition~\ref{prop:reconstruction-functor},
giving $\mathcal C'\simeqm\mathcal P_{\mathcal C}(\BW(k),\id,n,\theta)$.

The functors and tensorators in the proof of
Proposition~\ref{prop:gauge-functors} are determined on homogeneous
summands and extend over finite direct sums, so two such pullbacks are
monoidally superequivalent exactly when their twists are gauge equivalent.
A monoidal superequivalence preserves types, so $\phi=\id_{\BW(k)}$. The
remaining gauge equations are precisely those defining $\SH_k(\BW(k))$.
\end{proof}

\section{Reference categories over \texorpdfstring{$\C$, $\R$, and $\F_5$}{C, R, and F5}}
\label{sec:fields}

For $k = \mathbb C, \mathbb R$, the target spaces of $o_4$ are trivial, and reference categories exist. On the other hand, for $k = \mathbb F_5$, no such reference category exists.

\begin{lemma}
\label{lem:brauer-wall-fields}
For $u\in k^\times$, let $\Cl_1(u)$ denote the Clifford superalgebra
$k[\xi]/(\xi^2-u)$ with $\xi$ odd.
Let $\epsilon\colon\BW(k)\to\Z/2\Z$ be the parity homomorphism of
Definition~\ref{def:brauer-wall-group},  sending $[\Cl_1(u)]$ to $1$
for every $u\in k^\times$.

(1) $\BW(\C)\cong\Z/2\Z$, generated by $[\Cl_1(1)]$, and
$\epsilon$ is the identity.

(2) $\BW(\R)\cong\Z/8\Z$, generated by $[\Cl_1(1)]$, and
$\epsilon(a)=a\bmod2$.

(3) $\BW(\mathbb F_5) \cong \mathbb Z/2\mathbb Z \times \mathbb Z/2\mathbb Z$, generated by $[\Cl_1(1)]$ and $[\Cl_1(2)]$, and $\epsilon(a,b) = a+b$.
\end{lemma}
\begin{proof}
The ordinary center of $\Cl_1(u)$ is $k\oplus k\xi$, so
\eqref{eq:bw-parity} gives $\epsilon([\Cl_1(u)])=1$.
Over $\C$ and $\F_5$ the Brauer group is trivial and $-1$ is a square.
Wall's multiplication law \cite[Theorem~3]{Wall} therefore gives
\[
 \BW(k)\cong\Z/2\Z\times k^\times/(k^\times)^2,\qquad
 [\Cl_1(u)]\longmapsto(1,\bar u).
\]
The square classes are represented by $1$ over $\C$ and by $1,2$
over $\F_5$. Over $\R$, Wall's calculation gives a cyclic group of
order eight, in which every class of parity $1$ is a generator.
The formulas for $\epsilon$ follow by additivity.
\end{proof}

\begin{lemma}
\label{lem:cyclic-cohomology}
Let $M$ be an abelian group written multiplicatively with trivial
$\Z/n\Z$-action, and put $M^n=\{t^n:t\in M\}$ and $M[n]=\{t\in M:t^n=1\}$.
Then there are isomorphisms, natural in $M$,
\begin{equation}
\label{eq:cyclic-cohomology}
 H^2(\Z/n\Z,M)\cong M/M^n,\qquad
 H^3(\Z/n\Z,M)\cong M[n],\qquad
 H^4(\Z/n\Z,M)\cong M/M^n.
\end{equation}
Moreover, for a generator $u\in H^2(\Z/n\Z,\Z)$, the cup product
\begin{equation}
\label{eq:cyclic-periodicity}
 -\smile u:H^2(\Z/n\Z,M)\xrightarrow{\ \sim\ }H^4(\Z/n\Z,M)
\end{equation}
with respect to the pairing $M\times\Z\to M$, $(t,r)\mapsto t^r$, is an
isomorphism.
\end{lemma}
\begin{proof}
For $\Z/n\Z$ with generator $g$, the standard periodic resolution has
differentials $g-1$ and $N=1+g+\cdots+g^{n-1}$
\cite[Theorem~6.2.2]{Weibel}. Since the action on $M$ is trivial, the
induced cochain maps are $t\mapsto 1$ and $t\mapsto t^n$, giving the
natural isomorphisms \eqref{eq:cyclic-cohomology}. Cupping with a
generator $u\in H^2(\Z/n\Z,\Z)\cong\Z/n\Z$ gives the periodicity
isomorphism \eqref{eq:cyclic-periodicity}; see
\cite[Chapter~VI, \S9]{Brown}.
\end{proof}

\begin{lemma}
\label{lem:fields-complex-real-pentagons}
For $k=\C$ and for $k=\R$, the set $\Sec(D)$ of classes of even
Morita systems has two elements.

(1) Over $\C$, both classes have $O_{D,M}=1$.

(2) Over $\R$, exactly one class has $O_{D,M}=1$.
\end{lemma}
\begin{proof}
By Lemmas~\ref{lem:morita-systems-sections}
and~\ref{lem:cyclic-cohomology}, $\Sec(D)$ is a torsor under
$H^2(\BW(k),\Z/2\Z)\cong\Z/2\Z$ in both cases.
Over $\C$, $H^4(\BW(\C),\C^\times)\cong\C^\times/(\C^\times)^2=1$, giving (1).

Over $\R$, we have $\BW(\R)\cong\Z/8\Z$. Let $u\in H^2(\BW(\R),\Z)$ be a generator,
and let $p:\Z\to\Z/2\Z$ and $\psi:\Z\to\R^\times$ be $p(r)=r\bmod 2$ and
$\psi(r)=(-1)^r$. By the naturality of \eqref{eq:cyclic-cohomology},
$\bar u=p_*(u)$ is the unique nonzero element of $H^2(\BW(\R),\Z/2\Z)\cong\Z/2\Z$,
and $\psi_*(u)$ is the unique nontrivial element of
$H^2(\BW(\R),\R^\times)\cong\R^\times/(\R^\times)^8=\R^\times/\R_{>0}$.
For a cocycle $c\in Z^2(\BW(\R),\Z)$ representing $u$, we have
$(-1)^{(p_*c)\smile(p_*c)}=(\psi_*c)\smile c$ in $Z^4(\BW(\R),\R^\times)$ with
respect to the pairing $(t,r)\mapsto t^r$, so
\[
 Q(\bar u)=\psi_*(u)\smile u\in H^4(\BW(\R),\R^\times).
\]
By \eqref{eq:cyclic-periodicity}, $Q(\bar u)$ is the unique nontrivial
element of $H^4(\BW(\R),\R^\times)\cong\R^\times/\R_{>0}$. By
Lemma~\ref{lem:pentagon-shift}, the two elements of $\Sec(D)$ therefore
have distinct pentagon classes in this two-element group, proving (2).
\end{proof}

\begin{lemma}
\label{lem:fields-f5-pentagon}
$o_4(\F_5)\ne 1$.
\end{lemma}
\begin{proof}
Put $k=\F_5$, $D=\Cl_1(1)$, $a=[D]$, and restrict to $H=\{0,a\}$.
Use the regular bimodules and unit comparisons at $0$.
Take $M_{a,a}=N=D$ as a graded space,
with scalar right action and with $\xi\gt1$ and $1\gt\xi$ acting by
\begin{equation}
\label{eq:f5-N-action}
 Xr=\xi r,\qquad Yr=2(-1)^{|r|}r\xi .
\end{equation}
These operators square to $1$ and anticommute. Their product $XY$ has
eigenvalues $2,-2$ on $1,\xi$, and $X$ exchanges the two eigenlines.
Thus $D\gt D\cong\End_k(N)$, with $N$ its column module.
At $(a,a,a)$ the source $N\gt D$ is generated by $u=1\gt1$ with
$\xi_2u=2\xi_1u$ and $\xi_3u=u\xi$, where $\xi_i$ denotes $\xi$ in
the $i$th tensor factor. The vector $1\gt1-2\xi\gt\xi$ in $D\gt N$
satisfies these relations, giving the even bimodule map
\begin{equation}
\label{eq:f5-A0}
 A:N\gt D\longrightarrow D\gt N,\qquad
 A(r\gt s)=rs\gt1-2(-1)^{|s|}r\xi s\gt\xi
\end{equation}
on homogeneous tensors. It is invertible by Proposition~\ref{prop:morita-equivalences}.

In the fourfold pentagon, write $v_L$ for the tensor of units in the
source. Identify the target with $N\otimes_k N$ by
$(1\gt(1\gt r))\otimes s\mapsto r\otimes s$, and put $v=1\otimes1$.
The four odd generators $\xi_1,\ldots,\xi_4$ act there as
$\id\gt X,\id\gt Y,X\gt\id,Y\gt\id$, respectively.
The short path uses only unit comparisons. On the long path,
\eqref{eq:f5-A0} and $Y(1)=2X(1)$ give the factors
$1-\xi_1\xi_3$ and $1-\xi_2\xi_4$ on the two overlapping triples. Hence
\begin{equation}
\label{eq:f5-O2}
 P_R(v_L)=v,\qquad P_L(v_L)=(1-\xi_1\xi_3)(1-\xi_2\xi_4)v=2v .
\end{equation}
Indeed, for $w=\xi\otimes\xi$ one has $\xi_2\xi_4v=w$,
$\xi_1\xi_3v=-w$ and $\xi_1\xi_3w=v$. By \eqref{eq:pentagon-scalar}
this gives $O_{D,M,A}(a,a,a,a)=2$.

Every other choice of $M_{a,a}$ is evenly isomorphic to $N$ or $\Pi N$.
The shift contributes $-1$ at $(a,a,a,a)$ by Lemma~\ref{lem:pentagon-shift},
and replacing $A$ by $\lambda A$ contributes $\lambda^2$ because the
long path uses $A$ twice. Here $\lambda\in k^\times$.
Thus the possible scalars are $\pm2\lambda^2\in\{2,3\}$, since $-1=2^2$ in $k$.
A reference system would restrict to one with scalar $1$, so
Proposition~\ref{prop:scalar-obstruction} gives $o_4(k)\ne1$.
\end{proof}

\begin{theorem}
\label{thm:three-fields}
Up to monoidal superequivalence, the Brauer--Wall reference categories
over the three fields are as follows.

(1) Over $\C$ there are four classes of reference categories.

(2) Over $\R$ there are two classes of reference categories.

(3) Over $\F_5$ no reference category exists.
\end{theorem}
\begin{proof}
Lemmas~\ref{lem:fields-complex-real-pentagons}
and~\ref{lem:fields-f5-pentagon}, together with
Theorem~\ref{thm:existence-of-reference}, give existence over $\C$ and
$\R$ and nonexistence over $\F_5$.
Fix a reference category $\mathcal C$ over $k=\C$ or $\R$.
By Proposition~\ref{prop:reference-torsor} and
\eqref{eq:supercohomology-extension}, the number of classes is
\[
 |\SH_k(\BW(k))|=|H^3(\BW(k),k^\times)|\cdot|\ker Q|.
\]
Over $\C$, $H^4(\BW(\C),\C^\times)=1$ implies
\[
 \ker Q=H^2(\BW(\C),\Z/2\Z)\cong\Z/2\Z,
\]
while $H^3(\BW(\C),\C^\times)\cong\C^\times[2]\cong\Z/2\Z$ by
Lemma~\ref{lem:cyclic-cohomology}, giving $|\SH_\C(\BW(\C))|=4$.
Over $\R$, the proof of Lemma~\ref{lem:fields-complex-real-pentagons}
shows that $Q$ is injective on $H^2(\BW(\R),\Z/2\Z)$, so $\ker Q=0$, while
$H^3(\BW(\R),\R^\times)\cong\R^\times[8]=\{1,-1\}$, giving $|\SH_\R(\BW(\R))|=2$.
\end{proof}

\begin{remark}
  \label{rem:explicit-construction}
  When $k = \mathbb R$, there are two classes of even Morita systems; however, the above argument does not specify which one serves as the reference system.
  In the Lean formalization \cite{RealRModLean}, we explicitly constructed a reference category over $\mathbb R$.
  Although an explicit construction is unnecessary for the argument of this paper, we carried it out partially because we did not notice the quotient by $\operatorname{im} Q$ in \eqref{eq:reference-obstruction-class} at first, and partially because AI's brutal computing power enabled it.
  The reason we studied the case $k = \mathbb F_5$ is that we wanted to know whether there exists a general mechanism making the obstruction class trivial without direct computation.
\end{remark}

\begin{remark}
  \label{rem:homotopy-theory}
  As in \cite{EtingofNikshychOstrik}, there is also a topological description in which the algebraic computations are hidden in the machinery of homotopy theory; understanding the cases $k = \mathbb C, \mathbb R$ from this viewpoint is left to the author's future work.
\end{remark}

\begin{remark}
  \label{rem:operator-algebra}
  Another short proof of the existence of reference categories for $k = \mathbb C,\mathbb R$ using operator-algebraic techniques is given in Appendix~\ref{app:short-existence-proof}.
\end{remark}

\begin{remark}
  \label{rem:other-systems}
  The other class of even systems over $\mathbb R$, such that $O_{D,M} \neq 1$, can also be used as in Section~\ref{sec:reference-category} to construct weakly pointed fusion $\Pi$-supercategories, by modifying the twist definition by $\delta \omega = (-1)^{\nu\smile \nu}\cdot \rho^* O_{D,M,A}^{-1}$.
  More generally, over an even system $(D,M,A)$ over $k$, we can construct a weakly pointed fusion $\Pi$-supercategory from a twist $(G,\rho,\nu,\omega)$ satisfying $\delta \omega = (-1)^{\nu\smile \nu}\cdot \rho^* O_{D,M,A}^{-1}$. 
  For a given $(G,\rho,\nu)$, such an $\omega$ exists
if and only if $[\rho^*O_{D,M,A}]=[(-1)^{\nu\smile\nu}]\in H^4(G,k^\times)$. 
  The results of this section imply that, over $k = \mathbb C, \mathbb R$, we can take an even system such that $O_{D,M,A} = 1$ as a system over which weakly pointed fusion $\Pi$-supercategories are constructed. 
\end{remark}

\section{Stacking over \texorpdfstring{$\C$ and $\R$}{C and R}}
\label{sec:complex-real}

Since the stacking formula in \eqref{eq:twist-product} depends on the exchange data, we cannot compute it without explicitly constructing the reference categories. 
However, for $k = \mathbb C, \mathbb R$, we show that we can normalize the exchange data, making the formula \eqref{eq:twist-product} explicit.

\begin{definition}
\label{def:cyclic-cochains}
For $a\in\Z/n\Z$ let $\tilde a\in\{0,\ldots,n-1\}$ be its representative,
and let $\gamma_n\in Z^2(\Z/n\Z,\Z)$ be the carry cocycle,
\begin{equation}
\label{eq:carry-cocycle}
 \tilde a+\tilde b=\widetilde{a+b}+n\,\gamma_n(a,b),\qquad
 \gamma_n(a,b)\in\{0,1\}.
\end{equation}
We identify $\BW(\R)$ with $\Z/8\Z$ by sending
$[\Cl_1(-1)]=-[\Cl_1(1)]$ to $1$. Write $\gamma=\gamma_8$ and define
$q:\Z/8\Z\to\Z/2\Z$ by
\begin{equation}
\label{eq:bw-cochains}
 q(a)=\frac{\tilde a(\tilde a+1)}2\bmod2.
\end{equation}
For a finite group $G$ and $\rho\in Z^1(G,\Z/8\Z)$, define
$x_\rho,x'_\rho,x''_\rho:G\to\Z/2\Z$ and $\gamma_\rho\in Z^2(G,\Z/2\Z)$ by
\begin{equation}
\label{eq:real-bits}
 \widetilde{\rho(g)}=x_\rho(g)+2x'_\rho(g)+4x''_\rho(g),\qquad
 \gamma_\rho=\rho^*\gamma\bmod2.
\end{equation}
\end{definition}

\begin{lemma}
\label{lem:real-cochain-identities}
For $a,b\in\Z/8\Z$ and $\rho\in Z^1(G,\Z/8\Z)$,
\begin{equation}
\label{eq:bw-quadratic-identity}
 q(a+b)=q(a)+q(b)+\epsilon(a)\epsilon(b),\qquad
 q(-a)=q(a)+\epsilon(a),\qquad
 q(a)=\epsilon(a)+\lfloor\tilde a/2\rfloor\bmod2,
\end{equation}
\begin{equation}
\label{eq:real-second-bit}
 x_\rho=\epsilon\circ\rho,\qquad
 \delta x'_\rho=x_\rho\smile x_\rho,\qquad
 q\circ\rho=x_\rho+x'_\rho.
\end{equation}
\end{lemma}
\begin{proof}
The function $a(a+1)/2\bmod2$ is $8$-periodic, and
\[
 \frac{(a+b)(a+b+1)}2
 =\frac{a(a+1)}2+\frac{b(b+1)}2+ab
\]
gives the first identity. Set $b=-a$ for the second and separate
even and odd $a$ for the third. Binary addition modulo $4$ gives
$x(g)+x(h)=x(gh)+2x(g)x(h)$ as integers, hence
$\delta x'=x\smile x$. The remaining formulas follow from the bits.
\end{proof}

\begin{definition}
\label{def:standard-hexagon-scalars}
Identify $\BW(\C)=\Z/2\Z$ with $\{0,1\}\subset\Z$, let $\BW(\R)=\Z/8\Z$
with $\epsilon(a)=a\bmod2$ (Lemma~\ref{lem:brauer-wall-fields}), and let
$q:\Z/8\Z\to\Z/2\Z$ be the function \eqref{eq:bw-cochains}.
The \emph{standard hexagon scalars} over $\C$ and over $\R$ are the
pairs $\kappa^{\C}=(\kappa^{\C}_1,\kappa^{\C}_2)$ of maps
$\BW(\C)^3\to\C^\times$ and $\kappa^{\R}=(\kappa^{\R}_1,\kappa^{\R}_2)$
of maps $\BW(\R)^3\to\R^\times$ given by
\begin{equation}
\label{eq:standard-hexagon-scalars}
\begin{alignedat}{2}
 \kappa^{\C}_1(a;b,c)&=1,&\qquad \kappa^{\C}_2(a,b;c)&=i^{abc},\\
 \kappa^{\R}_1(a;b,c)&=1,&\qquad
 \kappa^{\R}_2(a,b;c)&=(-1)^{\epsilon(a)\epsilon(b)q(c)}.
\end{alignedat}
\end{equation}
\end{definition}

\begin{lemma}
\label{lem:real-theta-twist}
Let $k=\R$, let $(R,\beta_1,\beta_2)$ be an exchange data of
a Brauer--Wall reference category $\mathcal C$, and put
$\tau=(\Z/8\Z,\id,0,\omega_-)$, where
\begin{equation}
\label{eq:real-sign-cocycle}
 \omega_-(a,b,c)=(-1)^{\tilde a\gamma(b,c)}.
\end{equation}
Then Definition~\ref{def:pullback-construction} defines a Brauer--Wall reference category
$\mathcal C^\tau=\mathcal P_{\mathcal C}(\tau)$ with components
$\mathcal C^\tau_a$ the families supported at $a$, and
\begin{equation}
\label{eq:theta-twisted-exchange}
 R^\tau_{X,Y}=r_{Y\otimes X}^{-1}R_{X,Y}r_{X\otimes Y}
 :(X\otimes Y)\otimes\one\longrightarrow(Y\otimes X)\otimes\one,
\end{equation}
on families supported at single elements, is part of an exchange data
$(R^\tau,\beta^-_1,\beta^-_2)$ of $\mathcal C^\tau$ such that
\begin{equation}
\label{eq:real-theta-twist}
 \beta^-_1(a;b,c)=\beta_1(a;b,c)\,(-1)^{\tilde a\gamma(b,c)},\qquad
 \beta^-_2(a,b;c)=\beta_2(a,b;c)\,(-1)^{\gamma(a,b)\tilde c}.
\end{equation}
\end{lemma}
\begin{proof}
The exponent of $\omega_-$ modulo $2$ is the cup product
$\epsilon\smile(\gamma\bmod2)$, so $\omega_-$ is a normalized $3$-cocycle.
The reference-category assertion follows from
Proposition~\ref{prop:reference-torsor}. Delete the inserted units;
then $R^\tau$ becomes $R$ and the associator becomes
$\theta(a,b,c)\alpha$, with $\theta=\omega_-$.
Substitution in the two hexagons gives
\begin{align*}
 \beta^\theta_1(a;b,c)
 &=\beta_1(a;b,c)\frac{\theta(a,b,c)\theta(b,c,a)}{\theta(b,a,c)},\\
 \beta^\theta_2(a,b;c)
 &=\beta_2(a,b;c)\frac{\theta(a,c,b)}{\theta(c,a,b)\theta(a,b,c)}.
\end{align*}
The formula for $\omega_-$ and symmetry of $\gamma$ reduce these
ratios to the two signs in \eqref{eq:real-theta-twist}.
\end{proof}

\begin{proposition}
\label{prop:standard-exchange-data}
(1) Let $k=\C$. Every Brauer--Wall reference category $\mathcal C$ over
$\C$ has exchange degree $b_{\mathcal C}(a,b)=ab$, and there is a
Brauer--Wall reference category over $\C$ with an exchange data
$(R,\beta_1,\beta_2)$ such that $(\beta_1,\beta_2)=\kappa^{\C}$.

(2) Let $k=\R$. Every Brauer--Wall reference category $\mathcal C$ over
$\R$ has exchange degree $b_{\mathcal C}(a,b)=\epsilon(a)\epsilon(b)$,
and there is a Brauer--Wall reference category over $\R$ with an
exchange data $(R,\beta_1,\beta_2)$ such that
$(\beta_1,\beta_2)=\kappa^{\R}$.
\end{proposition}
\begin{proof}
The exchange degrees follow from cyclicity and
\eqref{eq:exchange-degree-laws}.

(1) Put $A=\beta_1(1;1,1)$ and $B=\beta_2(1,1;1)$.
The mixed hexagon identity gives $A^2=-B^2$.
Rescale $R_{1,1}$ by a square root of $A$ to obtain
$(\beta_1,\beta_2)=(1,\varsigma^{abc})$, with $\varsigma=\pm i$.
If necessary, complex conjugation of the reference category changes
$-i$ to $i$ and preserves its types.

(2) Divide the hexagon scalars by $\kappa^\R$, which satisfies
\eqref{eq:exchange-datum-identities} by
Lemma~\ref{lem:real-cochain-identities}. The quotient satisfies
the same identities without the mixed sign.
For a normalized $2$-cocycle $f$ on $\Z/8\Z$ with values in
$\R^\times$ or $\{\pm1\}$, put $t=\prod_{j=0}^7 f(j,1)$ and
$\eta(b)=\prod_{j=0}^{\tilde b-1}f(j,1)^{-1}$.
The cocycle identity gives $f(b,c)=t^{\gamma(b,c)}(\delta\eta)(b,c)$.
For each $a$, this lets us rescale the first component to
$u(a)^{\gamma(b,c)}$.
The mixed identity makes the sign of $u$ a character; removing
its positive part leaves $u(a)=(-1)^{e_1\tilde a}$.
The second component is then a character in $c$ with a
$\{\pm1\}$-valued $2$-cocycle as coefficient. A further rescaling
therefore gives the normal form
\[
 \kappa^e_1(a;b,c)=(-1)^{e_1\tilde a\gamma(b,c)},\qquad
 \kappa^e_2(a,b;c)
 =(-1)^{\epsilon(a)\epsilon(b)q(c)+e_2\gamma(a,b)\tilde c}.
\]
To show $e_1=e_2$, use commutativity of stacking on
$G=\Z/8\Z\times\Z/2\Z$. Take $\rho$ the first projection,
$\rho'=4y$ for the second projection $y$, and $x=\epsilon\rho$.
The products $p_\rho\star p_{\rho'}$ and $p_{\rho'}\star p_\rho$
have zero parity cochain and scalar cochains
\[
 (-1)^{e_2y\smile y\smile x},\qquad
 (-1)^{e_1x\smile y\smile y}.
\]
For the generators $g,h$ of the two factors, the functional
$\Psi(f)=f(g,h,h)+f(h,g,h)+f(h,h,g)$ vanishes on mod-$2$
coboundaries and takes value $1$ on both displayed exponents
without $e_1,e_2$. Gauge equivalence, followed by taking signs,
therefore forces $e_1=e_2$.
The case $(0,0)$ is standard; in the case $(1,1)$, twisting by
$\omega_-$ removes both carry factors by
Lemma~\ref{lem:real-theta-twist}.
\end{proof}

\begin{definition}
\label{def:standard-reference-categories}
By Proposition~\ref{prop:standard-exchange-data}, we may fix a
Brauer--Wall reference category $\mathcal C_{\C}$ over $\C$ with an
exchange data $(R,\beta_1,\beta_2)$ such that
$(\beta_1,\beta_2)=\kappa^{\C}$, and a Brauer--Wall reference category
$\mathcal C_{\R}$ over $\R$ with an exchange data $(R,\beta_1,\beta_2)$
such that $(\beta_1,\beta_2)=\kappa^{\R}$. We call them the
\emph{standard reference categories} and the fixed exchange data their
\emph{standard exchange data}.
In the rest of this section, for $k=\C$ and for $k=\R$, the product
$\star$ of Definition~\ref{def:cochain-operations} and the group
structure on $\Tcal_k(G)$ of Theorem~\ref{thm:stacking-group-categorical}
are those defined by the standard reference category over $k$ with its
standard exchange data.
\end{definition}

\begin{proposition}
\label{prop:standard-products}
Let $G$ be a finite group and let $s,t\in\Tw_k(G)$.

(1) Let $k=\C$, and let $x=\rho_s$ and $y=\rho_t$, regarded as elements
of $Z^1(G,\Z/2\Z)$. Then
\begin{equation}
\label{eq:complex-twist-product}
\begin{aligned}
 s\star t&=\bigl(G,\ \rho_s+\rho_t,\ \nu_s+\nu_t+y\smile x,\
   \omega_s\omega_t\,i^{\,y\smile y\smile x}\,(-1)^{\Xi^{\C}_{s,t}}\bigr),\\
 \Xi^{\C}_{s,t}&=\nu_s\smile_1\nu_t+(y\smile x)\smile_1(\nu_s+\nu_t)
   +y\smile(yx)\smile x,
\end{aligned}
\end{equation}
and
\begin{equation}
\label{eq:complex-twist-inverse}
\begin{aligned}
 s^\vee&=\bigl(G,\ -\rho_s,\ \nu_s+x\smile x,\
   \omega_s^{-1}\,i^{\,x\smile x\smile x}\,(-1)^{\Xi^{\C\vee}_{s}}\bigr),\\
 \Xi^{\C\vee}_{s}&=(\nu_s+x\smile x)\smile_1\nu_s,
\end{aligned}
\end{equation}
where $(yx)(h)=y(h)x(h)$, $i^{\,y\smile y\smile x}(g,h,\ell)=i^{y(g)y(h)x(\ell)}$, and
$i^{\,x\smile x\smile x}(g,h,\ell)=i^{x(g)x(h)x(\ell)}$.

(2) Let $k=\R$, and let $x=\epsilon\circ\rho_s$ and
$y=\epsilon\circ\rho_t$. Then
\begin{equation}
\label{eq:real-twist-product}
\begin{aligned}
 s\star t&=\bigl(G,\ \rho_s+\rho_t,\ \nu_s+\nu_t+y\smile x,\
   \omega_s\omega_t(-1)^{\Xi^{\R}_{s,t}}\bigr),\\
 \Xi^{\R}_{s,t}(g,h,\ell)
 &=(\nu_s\smile_1\nu_t)(g,h,\ell)
   +((y\smile x)\smile_1(\nu_s+\nu_t))(g,h,\ell)\\
 &\qquad+y(g)y(h)x(h)x(\ell)+y(g)y(h)\,q(\rho_s(\ell)),
\end{aligned}
\end{equation}
and
\begin{equation}
\label{eq:real-twist-inverse}
\begin{aligned}
 s^\vee&=\bigl(G,\ -\rho_s,\ \nu_s+x\smile x,\
   \omega_s^{-1}(-1)^{\Xi^{\R\vee}_{s}}\bigr),\\
 \Xi^{\R\vee}_{s}(g,h,\ell)
 &=((\nu_s+x\smile x)\smile_1\nu_s)(g,h,\ell)+x(g)x(h)\,q(\rho_s(\ell)).
\end{aligned}
\end{equation}
\end{proposition}
\begin{proof}
Substitute the standard hexagon scalars into
\eqref{eq:stacking-factor-cochains}. In both cases $c_{s,t}=y\smile x$,
while
\[
 K_{s,t}(g,h,\ell)=
 \begin{cases}
 i^{y(g)y(h)x(\ell)}&k=\C,\\
 (-1)^{y(g)y(h)q(\rho_s(\ell))}&k=\R.
 \end{cases}
\]
This gives the product formulas. For the inverses use
\begin{equation}
\label{eq:cup-one-vanishing}
 (x\smile x)\smile_1(x\smile x)=0,
\end{equation}
which follows directly from additivity of $x$.
Over $\C$, the remaining cubic sign combines with $i^{-x^3}$
to give $i^{x^3}$, where $x^3=x\smile x\smile x$.
Over $\R$, $q(-a)=q(a)+\epsilon(a)$ cancels the same cubic sign.
\end{proof}

\begin{theorem}
\label{thm:complex-stacking-group}
Let $k=\C$ and let $G$ be a finite group.
If some $\rho\in Z^1(G,\Z/2\Z)$ is not the reduction modulo $2$ of a
homomorphism $G\to\Z/4\Z$, then the exact sequence
\begin{equation}
\label{eq:complex-degree-sequence}
 0\longrightarrow\SH_\C(G)\longrightarrow\Tcal_\C(G)
 \xrightarrow{\ d\ } Z^1(G,\Z/2\Z)\longrightarrow0
\end{equation}
of Theorem~\ref{thm:stacking-group} does not split.
For example, if $G=\Z/2\Z$ and $\rho=\id$, then 
\begin{equation}
\label{eq:complex-stacking-group}
 \SH_\C(\Z/2\Z)\cong\Z/4\Z,\qquad
 \Tcal_\C(\Z/2\Z)\cong\Z/8\Z\ \text{generated by }[p_\rho],
\end{equation}
and \eqref{eq:complex-degree-sequence} does not split.
\end{theorem}
\begin{proof}
The product formula gives
\[
 p_\rho\star p_\rho
 =\bigl(G,0,\rho\smile\rho,(-i)^{\rho\smile\rho\smile\rho}\bigr).
\]
If the degree sequence had a section $\sigma$, then
$f(\rho)=[p_\rho]\sigma(\rho)^{-1}\in\SH_\C(G)$ would satisfy
$[p_\rho]^2=f(\rho)^2$. Mapping to $H^2(G,\Z/2\Z)$ gives
$[\rho\smile\rho]=0$.
Choose $\lambda$ with $\delta\lambda=\rho\smile\rho$ and integer
lifts taking values in $\{0,1\}$. Then
\[
 \delta(\tilde\rho-2\tilde\lambda)
 =2(\tilde\rho\smile\tilde\rho-\delta\tilde\lambda)
 \equiv0\pmod4,
\]
so $\rho$ lifts to a homomorphism into $\Z/4\Z$.
Now take $G=\Z/2\Z$ and $\rho=\id$, and put
$\omega=(-i)^{\rho\smile\rho\smile\rho}$.
By \eqref{eq:degree-zero-product} and \eqref{eq:cup-one-vanishing},
$[p_\rho]^4$ is the image of $[\omega^2]$ under the injection in
\eqref{eq:supercohomology-extension}. The cocycle $\omega^2$ takes
the value $-1$ at $(1,1,1)$, where every normalized $3$-coboundary
takes the value $1$. Since $(\omega^2)^2=1$, $[p_\rho]$ has order eight.
By Lemma~\ref{lem:cyclic-cohomology}, both $H^3(G,\C^\times)$ and
$H^2(G,\Z/2\Z)$ have order two. Thus \eqref{eq:supercohomology-extension}
gives $|\SH_\C(G)|\le4$, so $[p_\rho]^2$ generates $\SH_\C(G)$.
Exactness of \eqref{eq:complex-degree-sequence} now gives
$\Tcal_\C(G)=\langle[p_\rho]\rangle\cong\Z/8\Z$.
\end{proof}

\begin{theorem}
\label{thm:real-stacking-group}
Let $k=\R$ and let $G$ be a finite group. The exact sequence
\begin{equation}
\label{eq:real-degree-sequence}
 0\longrightarrow\SH_\R(G)\longrightarrow\Tcal_\R(G)
 \xrightarrow{\ d\ } Z^1(G,\Z/8\Z)\longrightarrow0
\end{equation}
of Theorem~\ref{thm:stacking-group} splits.
\end{theorem}
\begin{proof}
Write $x=x_\rho$, $x'=x'_\rho$. The real product formula and
$q\rho=x+x'$ give
\[
 p_\rho\star p_\rho
 =\bigl(G,2\rho,x\smile x,(-1)^{x\smile x\smile x'}\bigr).
\]
Since $\delta x'=x\smile x$, the gauge
$(\id_G,x',(-1)^{x'\smile x'})$ transforms this into $p_{2\rho}$:
the scalar cancellation is
$\delta(x'\smile x')=x\smile x\smile x'+x'\smile x\smile x$.
Hence $[p_\rho]^2=[p_{2\rho}]$.
Iteration shows $[p_\rho]^{2^j}=[p_{2^j\rho}]$; since every $\rho$
has order dividing eight and $d[p_\rho]=\rho$, their orders agree.
Choose cyclic generators $e_i$ of $Z^1(G,\Z/8\Z)$ and define
\[
 h\Bigl(\sum_i a_i e_i\Bigr)=\prod_i[p_{e_i}]^{a_i}.
\]
The order equality and commutativity make $h$ a well-defined
homomorphism, and $d\circ h=\id$.
\end{proof}
\appendix
\section{Cochains}
\label{app:cochains}

\begin{definition}
\label{def:cochains}
Let $G$ be a group and $M$ an abelian group.
For $n\geq1$, a \emph{normalized inhomogeneous $n$-cochain} is a map
$c:G^n\to M$ such that $c(g_1,\ldots,g_n)=0$ if some $g_i=1$.
We write $C^n(G,M)$ for the abelian group of these cochains, and we
put $C^0(G,M)=M$.
When $M$ is written multiplicatively, normalization means value $1$ and
the formulas below are read multiplicatively; all coefficient groups
carry the trivial $G$-action.
The
\emph{coboundary map} $\delta:C^n(G,M)\to C^{n+1}(G,M)$ is 
\begin{equation}
\label{eq:cochain-coboundary}
\begin{aligned}
 (\delta c)(g_1,\ldots,g_{n+1})
  ={}&c(g_2,\ldots,g_{n+1})\\
  &+\sum_{i=1}^n(-1)^i
    c(g_1,\ldots,g_i g_{i+1},\ldots,g_{n+1})\\
  &+(-1)^{n+1}c(g_1,\ldots,g_n).
\end{aligned}
\end{equation}
On $C^0(G,M)=M$ the coboundary map is $0$, the constant $1$ in
multiplicative notation, in accordance with the trivial action; thus
$B^1(G,M)=0$.
A cochain is a \emph{cocycle} if its coboundary is $0$. We write
$Z^n(G,M)$ for the group of cocycles, $B^n(G,M)=\delta C^{n-1}(G,M)$ for
the coboundaries, and $H^n(G,M)=Z^n(G,M)/B^n(G,M)$.
A group homomorphism $\phi:G\to H$ induces the \emph{pullback}
$\phi^*:C^n(H,M)\to C^n(G,M)$ given by
\begin{equation}
\label{eq:cochain-pullback}
 (\phi^*c)(g_1,\ldots,g_n)=c(\phi(g_1),\ldots,\phi(g_n)).
\end{equation}
A homomorphism of abelian groups $\psi:M\to M'$ induces
$\psi_*:C^n(G,M)\to C^n(G,M')$, $c\mapsto\psi\circ c$; for
$c\in C^n(G,\Z/2\Z)$ we write $(-1)^c\in C^n(G,k^\times)$ for
$\psi_*c$ with $\psi:\Z/2\Z\to k^\times$, $r\mapsto(-1)^r$.
If $M$ is a ring, the \emph{cup product} of $a\in C^p(G,M)$ and
$b\in C^q(G,M)$ is the cochain $a\smile b\in C^{p+q}(G,M)$ given by
\begin{equation}
\label{eq:cochain-cup-product}
 (a\smile b)(g_1,\ldots,g_{p+q})
  =a(g_1,\ldots,g_p)b(g_{p+1},\ldots,g_{p+q}).
\end{equation}
More generally, any biadditive pairing $\langle-,-\rangle:M_1\times M_2\to M_3$
of abelian groups defines a cup product
$C^p(G,M_1)\times C^q(G,M_2)\to C^{p+q}(G,M_3)$ by
$(a\smile b)(g_1,\ldots,g_{p+q})
=\langle a(g_1,\ldots,g_p),b(g_{p+1},\ldots,g_{p+q})\rangle$,
which induces a cup product on cohomology.
For $v,w\in C^2(G,\Z/2\Z)$, put
\begin{equation}
\label{eq:cup-one}
 (v\smile_1w)(g,h,\ell)=v(gh,\ell)w(g,h)+v(g,h\ell)w(h,\ell).
\end{equation}
\end{definition}

\begin{lemma}
\label{lem:coboundary-identities}
Let $G$ and $H$ be groups, let $M$ and $M'$ be abelian groups, let
$\phi:G\to H$ be a group homomorphism, and let $\psi:M\to M'$ be a
homomorphism of abelian groups.

(1) $\delta\circ\delta=0$, and $H^n(G,M)$ is the cohomology of $G$
with coefficients in $M$ regarded as a trivial $G$-module
\cite[Chapter~I, \S5 and Chapter~III, \S1]{Brown}.

(2) $\delta(\phi^*c)=\phi^*(\delta c)$ for $c\in C^n(H,M)$, and
$\delta(\psi_*c)=\psi_*(\delta c)$ for $c\in C^n(G,M)$. In particular,
for $c\in C^n(G,\Z/2\Z)$,
\begin{equation}
\label{eq:sign-coboundary}
 \delta\bigl((-1)^{c}\bigr)=(-1)^{\delta c}.
\end{equation}

(3) If $M$ is a commutative ring, then for $a\in C^p(G,M)$ and
$b\in C^q(G,M)$,
\begin{equation}
\label{eq:cochain-leibniz}
 \delta(a\smile b)=\delta a\smile b+(-1)^p a\smile\delta b .
\end{equation}
\end{lemma}
\begin{proof}
(1) Applying \eqref{eq:cochain-coboundary} twice produces terms
that occur in pairs with opposite signs, so $\delta^2=0$.
The normalized cochains form the cochain complex of the normalized
bar resolution, which computes group cohomology
\cite[Chapter~I, \S5 and Chapter~III, \S1]{Brown}.
(2) Pullback and coefficient homomorphisms commute termwise with
the coboundary formula.
(3) Expand $\delta a\smile b+(-1)^p a\smile\delta b$ using
\eqref{eq:cochain-coboundary}. The last term of $\delta a\smile b$
and the first term of $(-1)^p a\smile\delta b$ cancel; the remaining
terms are exactly $\delta(a\smile b)$, giving the Leibniz rule.
\end{proof}

\begin{lemma}
\label{lem:cup-identities}
Let $v,w,\nu\in Z^2(G,\Z/2\Z)$ and $\lambda\in C^1(G,\Z/2\Z)$.

(1) 
\begin{equation}
\label{eq:stacking-cup-one-identities}
 \delta(v\smile_1w)=v\smile w+w\smile v.
\end{equation}

(2) 
\begin{equation}
\label{eq:gauge-obstruction-coboundary}
\begin{aligned}
 (\nu+\delta\lambda)\smile(\nu+\delta\lambda)+\nu\smile\nu
 &=\delta\bigl((\nu+\delta\lambda)\smile\lambda+\lambda\smile\nu\bigr)\\
 &=\delta(\lambda\smile\nu+\nu\smile\lambda+\lambda\smile\delta\lambda).
\end{aligned}
\end{equation}
\end{lemma}
\begin{proof}
For (1), expand \eqref{eq:cup-one} and use the cocycle equations
for $v,w$; the surviving terms are
$v(g,h)w(\ell,m)+w(g,h)v(\ell,m)$.
For (2), put $\bar\nu=\nu+\delta\lambda$ and apply the Leibniz rule:
\[
 \delta(\bar\nu\smile\lambda+\lambda\smile\nu)
 =\bar\nu\smile(\bar\nu+\nu)+(\bar\nu+\nu)\smile\nu
 =\bar\nu\smile\bar\nu+\nu\smile\nu.
\]
Expanding
$\delta(\lambda\smile\nu+\nu\smile\lambda+\lambda\smile\delta\lambda)$
by the same rule gives the same result.
\end{proof}

\section{Another proof of the existence over \texorpdfstring{$\C$ and $\R$}{C and R}}
\label{app:short-existence-proof}

The following is another short proof of the existence of Brauer--Wall
reference categories over $k=\C$ and $k=\R$, obtained by consulting with
AI and Y.~Tachikawa.

Throughout this appendix, $k$ is $\R$ or $\C$, and
\begin{equation}
 n=4\ \text{ if }k=\R,\qquad n=1\ \text{ if }k=\C .
\end{equation}
We use the notions of Definitions~\ref{def:weakly-pointed}
and~\ref{def:brauer-wall-reference} and the description of $\BW(k)$ in
Lemma~\ref{lem:brauer-wall-fields}.
For $k=\C$ the central division superalgebras are $D_0=\C$ and
$D_1=\C[\xi]/(\xi^2-1)$ with $\xi$ odd.
For $k=\R$, let $\Cl_{p,q}$ be the real Clifford algebra with $p$
positive and $q$ negative squares. In the coordinates of
Definition~\ref{def:cyclic-cochains}, $D_r$ is its division representative
when $r\equiv q-p\pmod 8$.
In both cases the Brauer--Wall group is $\Z/2n\Z$.

\begin{theorem}
\label{alt:thm:main}
Let $k$ be $\R$ or $\C$.
There exists a Brauer--Wall reference category over $k$.
\end{theorem}

Explicitly, we construct a monoidal supercategory $\mathcal R$ over $k$
with an object $\pi$ and an odd isomorphism $\zeta:\pi\to\one$ such that
the underlying category $\underline{\mathcal R}$ is semisimple with
finitely many simple objects and rigid; $\End(\one)=k$; for every simple
object $X$ the object $X\otimes X^*$ is a direct sum of copies of $\one$
and $\pi$ in $\underline{\mathcal R}$; the endomorphism superalgebra of
every simple object is isomorphic to some $D_r$; every $r\in\Z/2n\Z$
occurs; and simple objects with the same $r$ are isomorphic in
$\mathcal R$.

The category is realized by endomorphisms of an algebra $L'$ which is
obtained from the Clifford algebra $\Cl_{\infty,\infty}$, on generators
$e_m$ and $f_m$ ($m\ge1$) with $e_m^2=1$ and $f_m^2=-1$, by adjoining two
elements.
The generating object $\sigma_+$ shifts the $e_m$ and its dual $\sigma_-$
shifts the $f_m$. The first adjoined element is an isometry $s$ which
implements the shift $\rho=\sigma_+\sigma_-$ of all generators; it gives
$\sigma_+\sigma_-\cong\one\oplus\pi$ and expresses that $\Cl_{1,1}$ is a
matrix algebra. The second is an invertible element $U$ which gives
$\sigma_+^n\cong\sigma_-^n$, making the category finite, and expresses an
isomorphism $\Cl_{0,n}\cong\Cl_{n,0}$ of superalgebras over $k$. The proof
is uniform in $k$ and uses no facts about these algebras beyond what is
proved here.

\subsection{Endomorphisms of a superalgebra}

Let $A$ be a unital superalgebra with grading automorphism $\gamma$.

\begin{definition}
\label{alt:def:endcat}
The supercategory $\EndCat(A)$ has as objects the unital algebra
endomorphisms of $A$ commuting with $\gamma$. For $p\in\Z/2\Z$,
\begin{equation}
 \Hom(\rho,\tau)_p=\{x\in A_p:x\rho(a)=\gamma^p(\tau(a))\,x
 \ \text{ for all } a\in A\} .
\end{equation}
The composition of morphisms is the product of $A$. We set
$\rho\otimes\rho'=\rho\circ\rho'$ and $x\otimes x'=x\,\rho(x')$ for
$x\in\Hom(\rho,\tau)$, $x'\in\Hom(\rho',\tau')$.
\end{definition}

\begin{lemma}
\label{alt:lem:strict}
$\EndCat(A)$ is a strict monoidal supercategory with unit object $\id$.
\end{lemma}
\begin{proof}
For $x\in\Hom(\rho,\tau)_p$ and $x'\in\Hom(\rho',\tau')_{p'}$ we have
\begin{equation}
 x\rho(x')\,\rho\rho'(a)=x\,\rho\bigl(\gamma^{p'}(\tau'(a))\bigr)\rho(x')
 =\gamma^{p+p'}(\tau\tau'(a))\,x\rho(x'),
\end{equation}
so $x\otimes x'\in\Hom(\rho\rho',\tau\tau')_{p+p'}$. For
$y\in\Hom(\tau,\upsilon)$ and $y'\in\Hom(\tau',\upsilon')_{q'}$,
\begin{equation}
 (y\otimes y')(x\otimes x')=y\,\tau(y')\,x\,\rho(x')
 =(-1)^{q'p}\,yx\,\rho(y'x'),
\end{equation}
which is the interchange law. Associativity and unitality are strict.
\end{proof}

In particular an element of $\Hom(\rho,\tau)$ also lies in
$\Hom(\rho\delta,\tau\delta)$ for every object $\delta$, and $\delta$ maps
$\Hom(\rho,\tau)$ into $\Hom(\delta\rho,\delta\tau)$, bijectively if
$\delta$ is an automorphism.

\subsection{The Clifford algebra and the algebra \texorpdfstring{$L$}{L}}

Let $C=\Cl_{\infty,\infty}$ be the $k$-algebra with generators $e_m,f_m$
($m\ge1$) which pairwise anticommute and satisfy $e_m^2=1$, $f_m^2=-1$. It
is a superalgebra with the generators odd. We call the $e_m$ and $f_m$
\emph{the generators}, write $\Cl(S)$ for the subalgebra generated by a
set $S$ of generators, and put $C_m=\Cl(e_1,\dots,e_m,f_1,\dots,f_m)$. Then
$C_m\cong M_{2^m}(k)$, the monomials in the generators form a basis of
$C$, and $C=\bigcup_mC_m$ is simple. If $x\in C_m$ has parity $p$ and $y$
is a generator $e_j$ or $f_j$ with $j>m$, then
\begin{equation}
\label{alt:eq:far}
 xy=(-1)^p\,yx .
\end{equation}
Let $\rho$ be the endomorphism of $C$ with $\rho(e_m)=e_{m+1}$ and
$\rho(f_m)=f_{m+1}$. It is injective and commutes with $\gamma$, and by
\eqref{alt:eq:far} $x\rho(a)=\rho(\gamma^p(a))\,x$ for $x\in C_1$ of parity
$p$ and $a\in C$. Finally let $P=\frac12(1+f_1e_1)\in C_1$. Since
$(f_1e_1)^2=1$, it is an even idempotent, and
\begin{equation}
\label{alt:eq:P}
 Pe_1=\tfrac12(e_1+f_1),\qquad e_1P=\tfrac12(e_1-f_1),\qquad
 e_1Pe_1=1-P,\qquad Pe_1P=0 .
\end{equation}

\begin{definition}
\label{alt:def:L}
$L$ is the algebra generated by $C$ and two elements $s,t$ with relations
\begin{equation}
\label{alt:eq:presentation}
 ts=1,\qquad st=P,\qquad sa=\rho(a)\,s\quad(a\in C).
\end{equation}
It is a superalgebra with $s,t$ even, and it has a $\Z$-grading in which
$C$ has degree $0$, $s$ has degree $1$ and $t$ has degree $-1$; $L_d$ is
the part of degree $d$. We put
\begin{equation}
 s_0=s,\quad s_1=e_1s,\quad t_0=t,\quad t_1=te_1,\qquad P_d=s^dt^d .
\end{equation}
\end{definition}

Thus $L$ is the Clifford algebra with one isometry adjoined which
implements the shift $\rho$. It acts on the Fock space of a half-infinite
chain, where $s$ inserts a new first site in its vacuum state; see the
proof of (4) below.

\begin{lemma}
\label{alt:lem:L}
(1) $at=t\rho(a)$ for $a\in C$. Moreover $s=Ps$ and $t=tP$, and the
element $P_d$ equals $P\rho(P)\cdots\rho^{d-1}(P)$, which lies in $C_d$.

(2) $t_is_j=\delta_{ij}$ and $s_0t_0+s_1t_1=1$, and the four elements
$s_it_j$ form a basis of $C_1$.

(3) $L_d=Cs^d$ and $L_{-d}=t^dC$ for $d\ge0$.

(4) The map $C\to L$ is injective. So $L_0=C$.

(5) The formula $\rho(x)=s_0xt_0+s_1\gamma(x)t_1$ defines an injective
object of $\EndCat(L)$ which extends $\rho$, and
\begin{equation}
\label{alt:eq:rhos}
 \rho(x)s_i=s_i\gamma^i(x),\qquad t_i\rho(x)=\gamma^i(x)t_i\qquad(x\in L),
\end{equation}
so that $s_i\in\Hom(\id,\rho)_i$ and $t_i\in\Hom(\rho,\id)_i$.
\end{lemma}
\begin{proof}
(1) $s=sts=Ps$ and $t=tst=tP$. As $P$ commutes with $\rho(a)$, we get
$at=ts\,at=t\rho(a)st=t\rho(a)P=t\rho(a)$, and
$P_{d+1}=sP_dt=\rho(P_d)\,st=\rho(P_d)P$.

(2) By \eqref{alt:eq:P}, $te_1s=tPe_1Ps=0$ and
$s_0t_0+s_1t_1=P+e_1Pe_1=1$, and
$s_0t_0=P$, $s_1t_1=1-P$, $s_0t_1=Pe_1$, $s_1t_0=e_1P$ span $C_1$.

(3) Every $c\in C$ is a sum of products $x\rho(a)$ with $x\in C_1$, and
$tx\rho(a)s=(txs)\,a$, where $txs=tPxPs$ is a scalar because $PC_1P=kP$.
Hence $tCs\subset C$. Together with $sa=\rho(a)s$ and $at=t\rho(a)$ this
shows that the span of the elements $a\,s^it^j\,b$ with $a,b\in C$ is
closed under multiplication, so it is $L$. Moreover $s^it^j=s^{i-j}P_j$
for $i\ge j$ and $s^it^j=P_it^{j-i}$ for $i\le j$, so
$a\,s^it^j\,b$ lies in $Cs^{i-j}$, respectively in $t^{j-i}C$.

(4) Let $V$ be the vector space with basis the sequences
$|i_1i_2i_3\cdots\rangle$ with $i_m\in\{0,1\}$ and $i_m=0$ for large $m$.
Let $X_m$ flip the $m$-th entry and let $Z_m$ multiply by $(-1)^{i_m}$.
The operators $e_m=Z_1\cdots Z_{m-1}X_m$ and $f_m=Z_1\cdots Z_mX_m$ satisfy
the relations of $C$. Let $s$ insert a new first entry $0$, and let $t$
remove the first entry if it is $0$ and give zero otherwise. Then $ts=1$,
$st=\frac12(1+Z_1)=P$, $se_m=e_{m+1}s$ and $sf_m=f_{m+1}s$. So $L$ acts on
$V$; in particular $L\ne0$, and $C\to L$ is injective because $C$ is
simple.

(5) By (2) the formula defines a unital homomorphism commuting with
$\gamma$ and satisfying \eqref{alt:eq:rhos}, and it is injective because
$t_0\rho(x)s_0=x$. For $c\in C$ it gives
$\rho(c)P+e_1\rho(\gamma(c))Pe_1=\rho(c)(P+e_1Pe_1)=\rho(c)$.
\end{proof}

\begin{remark}
\label{alt:rem:leavitt}
The elements $s_0,s_1,t_0,t_1$ generate $L$, and by
Lemma~\ref{alt:lem:L}(2) they satisfy the defining relations of the
Leavitt algebra \cite{Leavitt}, the algebraic part of the Cuntz algebra
$O_2$; one can show that these are defining relations of $L$. We do not
use this.
\end{remark}

Endomorphisms of $L$ are constructed as follows.

\begin{lemma}
\label{alt:lem:objects}
Let $\tau_0$ be an object of $\EndCat(C)$ and let $W\in C$ be even and
invertible, with
\begin{equation}
\label{alt:eq:E}
 \tau_0(x)=WxW^{-1}\quad(x\in C_1),\qquad
 \tau_0(\rho(a))=W\rho(\tau_0(a))W^{-1}\quad(a\in C).
\end{equation}
Then there is a unique object $\tau$ of $\EndCat(L)$ which extends
$\tau_0$ and satisfies $\tau(s)=Ws$, $\tau(t)=tW^{-1}$. We write
$\tau=[\tau_0,W]$. Moreover:

(1) $\tau(s_i)=Ws_i$, $\tau(t_i)=t_iW^{-1}$, and
$\tau(\rho(x))=W\rho(\tau(x))W^{-1}$ for $x\in L$.

(2) If $x\in\Hom(\tau,\tau)$, then $xW=W\rho(x)$.

(3) $[\tau_0,W]\circ[\tau'_0,W']=[\tau_0\tau'_0,\ \tau_0(W')W]$.

\noindent
It suffices to check \eqref{alt:eq:E} for $x=e_1,f_1$ and for the
generators $a$.
\end{lemma}
\begin{proof}
The relations \eqref{alt:eq:presentation} are preserved: $tW^{-1}Ws=1$,
$Ws\,tW^{-1}=WPW^{-1}=\tau_0(P)$, and
$Ws\,\tau_0(a)=W\rho(\tau_0(a))s=\tau_0(\rho(a))Ws$. (1) follows from
$\tau(e_1)=We_1W^{-1}$ and from the formula for $\rho$. (2) For $x$ of
parity $p$, $xWs_i=(-1)^{pi}Ws_ix=W\rho(x)s_i$ by \eqref{alt:eq:rhos};
multiply by $t_i$ and sum over $i$. (3) The composition sends $s$ to
$\tau_0(W')Ws$, and one checks \eqref{alt:eq:E} directly.
\end{proof}

The next two lemmas compute all the morphism spaces we need.

\begin{lemma}
\label{alt:lem:degree}
Let $x\in L$ and $a,b\in\Z$, and suppose that $xe_j=\pm e_{j+a}x$ and
$xf_j=\pm f_{j+b}x$ for all sufficiently large $j$, with signs which may
depend on $j$. If $a\ne b$ then $x=0$, and if $a=b$ then $x\in L_a$.
\end{lemma}
\begin{proof}
The components of $x$ of each $\Z$-degree and parity satisfy the same
relations, so let $x\in L_d$ have a definite parity. If $d\ge0$, then
$x=cs^d$ with $c\in C_i$ of definite parity, for some $i$. For large $j$
we have $xe_j=c\,e_{j+d}\,s^d$ and $e_{j+a}x=\pm c\,e_{j+a}s^d$ by
\eqref{alt:eq:far}.
Multiplying by $t^d$ on the right, and using that $P_d$ commutes with
$e_j$ for $j>d$, we get $cP_d\,(e_{j+d}\pm e_{j+a})=0$. If $a\ne d$, the
element $e_{j+d}\pm e_{j+a}$ is invertible, so $cP_d=0$ and
$x=cP_ds^d=0$. In the same way $x=0$ if $b\ne d$. If $d=-m<0$, then
$x=t^mc$, and $e_{j+a}x=t^me_{j+a+m}c$ while $xe_j=\pm t^me_jc$;
multiplying by $s^m$ on the left gives $(e_{j+a+m}\pm e_j)P_mc=0$, so
$x=t^mP_mc=0$ unless $a=d$, and likewise unless $b=d$.
\end{proof}

\begin{lemma}
\label{alt:lem:commutant}
Let $S$ be a finite set of generators and $p\in\Z/2\Z$. If $x\in C$
satisfies $xy=(-1)^pyx$ for every generator $y\notin S$, then
$x\in\Cl(S)$ and $x$ has parity $p$. In particular, for $S=\emptyset$,
$x$ is a scalar, and $x=0$ if $p=1$.
\end{lemma}
\begin{proof}
Let $x\in C_m$ with $\Cl(S)\subset C_m$. The parity components of $x$
satisfy the same relations, and by \eqref{alt:eq:far} the component of
parity $q$ satisfies $x_qe_{m+1}=(-1)^qe_{m+1}x_q$; hence $x$ has parity
$p$. Expand $x$ in monomials in the generators of $C_m$. Conjugation by a
generator $y$ multiplies a monomial $M$ of parity $p$ by $(-1)^p$ if $y$
does not occur in $M$ and by $-(-1)^p$ if it does. So only monomials in
the elements of $S$ occur.
\end{proof}

We say that an object $\tau$ of $\EndCat(L)$ \emph{has shift $(a,b)$} if
$\tau(e_j)=\pm e_{j+a}$ and $\tau(f_j)=\pm f_{j+b}$ for all sufficiently
large $j$. Shifts add under composition, and $\rho$ has shift $(1,1)$. If
$\tau$ and $\tau'$ have shifts $(a,b)$ and $(a',b')$, then by
Lemma~\ref{alt:lem:degree}
\begin{equation}
\label{alt:eq:index}
 \Hom(\tau,\tau')=0\ \text{ unless } a'-a=b'-b,\quad\text{and then }
 \Hom(\tau,\tau')\subset L_{a'-a}.
\end{equation}

\subsection{The shifts \texorpdfstring{$\sigma_+$ and $\sigma_-$}{sigma+ and sigma-}}

Let $\sigma_+$ and $\sigma_-$ be the endomorphisms of $C$ with
\begin{equation}
\label{alt:eq:sigmaC}
 \sigma_+(e_m)=-e_{m+1},\quad\sigma_+(f_m)=f_m,\qquad
 \sigma_-(e_m)=-e_m,\quad\sigma_-(f_m)=f_{m+1}\qquad(m\ge1),
\end{equation}
and let
\begin{equation}
 v_+=\tfrac1{\sqrt2}(e_1-e_2)e_1=\tfrac1{\sqrt2}(1+e_1e_2),\qquad
 v_-=\tfrac1{\sqrt2}(f_1-f_2)e_1 .
\end{equation}
These are even invertible elements of $C_2$ with $v_-^2=1$. Conjugation by
$v_+$ sends $e_1\mapsto-e_2$, $e_2\mapsto e_1$, conjugation by $v_-$ sends
$f_1\mapsto f_2$, $f_2\mapsto f_1$, $e_1\mapsto-e_1$, and both fix the
other generators. Comparing with \eqref{alt:eq:sigmaC}, the conditions
\eqref{alt:eq:E} hold for $(\sigma_\pm,v_\pm)$, so we have objects of
$\EndCat(L)$
\begin{equation}
 \sigma_+=[\sigma_+,v_+],\qquad\sigma_-=[\sigma_-,v_-],
\end{equation}
of shifts $(1,0)$ and $(0,1)$.

\begin{proposition}
\label{alt:prop:sigma}
(1) $\sigma_+\circ\sigma_-=\sigma_-\circ\sigma_+=\rho$.

(2) $e_1\in\Hom(\sigma_+,\sigma_+)_1$ and
$f_1\in\Hom(\sigma_-,\sigma_-)_1$. For $m\ge0$,
\begin{equation}
 \Hom(\sigma_+^m,\sigma_+^m)=\Cl(e_1,\dots,e_m),\quad
 \Hom(\sigma_-^m,\sigma_-^m)=\Cl(f_1,\dots,f_m),\quad
 \Hom(\rho,\rho)=C_1 .
\end{equation}

(3) $t_0\,\sigma_\pm(s_0)=\sigma_\pm(t_0)\,s_0=1/\sqrt2$.
\end{proposition}
\begin{proof}
(1) The three endomorphisms agree on $C$. Since
$\sigma_+(x)=v_+xv_+^{-1}$ for $x\in\Cl(e_1,f_1,f_2)$ and
$\sigma_-(v_+)=v_+$, both compositions send $s$ to $Ws$ and $t$ to
$tW^{-1}$ with $W=v_+v_-=\frac12(f_1-f_2)(e_1-e_2)$, by
Lemma~\ref{alt:lem:objects}(3). From $f_1e_1s=s$ we get $f_1s_0=-s_1$ and
$f_1s_1=s_0$, and $e_2s_i=(-1)^is_ie_1$, $f_2s_i=(-1)^is_if_1$ by
\eqref{alt:eq:rhos}. Hence
\begin{equation}
 2Ws_0=(f_1-f_2)(s_1-s_0e_1)=s_0(1+f_1e_1)+s_1(e_1+f_1)
 =2(s_0s_0t_0+s_1s_0t_1)=2\rho(s_0).
\end{equation}
Since $We_1W^{-1}=e_2=\rho(e_1)$, also
$Ws_1=We_1s_0=\rho(e_1)\rho(s_0)=\rho(s_1)$.
Therefore $W=\sum_iWs_it_i=\sum_i\rho(s_i)t_i$, whose inverse is
$\sum_is_i\rho(t_i)$, and $t_0W^{-1}=\rho(t_0)$.

(2) Since $e_1v_+=v_+e_2=v_+\rho(e_1)$, we get from \eqref{alt:eq:rhos}
\begin{equation}
 e_1\sigma_+(s)=\sigma_+(s)e_1,\qquad e_1\sigma_+(t)=\sigma_+(t)e_1 ,
\end{equation}
and $e_1$ anticommutes with the odd elements $\sigma_+(e_m)$ and
$\sigma_+(f_m)$; likewise for $f_1$ and $\sigma_-$, using
$f_1v_-=v_-\rho(f_1)$.
Hence the element $e_i=\pm1_{\sigma_+^{i-1}}\otimes e_1$ lies in
$\Hom(\sigma_+^i,\sigma_+^i)$, which is contained in
$\Hom(\sigma_+^m,\sigma_+^m)$, for $i\le m$. Conversely, an element $x$ of
$\Hom(\sigma_+^m,\sigma_+^m)_p$ lies in $L_0=C$ by \eqref{alt:eq:index},
and $xy=(-1)^pyx$ for the generators $y=\pm\sigma_+^m(e_j)=e_{j+m}$ and
$y=\sigma_+^m(f_j)=f_j$; so $x\in\Cl(e_1,\dots,e_m)$ by
Lemma~\ref{alt:lem:commutant}. The same argument applies to $\sigma_-^m$
and to $\rho$, where $C_1$ is spanned by the elements
$s_it_j\in\Hom(\rho,\rho)$.

(3) We have $e_2s_0=s_0e_1$, $f_2s_1=-s_1f_1$ and $t_0s_1=0$, so
$t_0v_+^{\pm1}s_0=\frac1{\sqrt2}(1\pm t_0e_1s_0e_1)=\frac1{\sqrt2}$ and
$t_0v_-s_0=\frac1{\sqrt2}t_0(f_1-f_2)s_1=\frac1{\sqrt2}t_0(s_0+s_1f_1)
=\frac1{\sqrt2}$, and $v_-^{-1}=v_-$.
\end{proof}

\subsection{The automorphism \texorpdfstring{$\beta$}{beta}}

The following is the only place where $k$ matters. We fix odd, pairwise
anticommuting elements
\begin{equation}
\label{alt:eq:conversion}
 g_1,\dots,g_n\in\Cl(e_1,\dots,e_n)\quad\text{with } g_i^2=-1,\ \text{ generating }
 \Cl(e_1,\dots,e_n).
\end{equation}
For $k=\C$ and $n=1$ we take $g_1=\sqrt{-1}\,e_1$. For $k=\R$ and $n=4$ we
take $g_i=e_i\,\omega$ with $\omega=e_1e_2e_3e_4$: since $\omega^2=1$ and
$\omega e_i=-e_i\omega$, these are odd, anticommute and square to $-1$,
and $e_i=g_i\,g_1g_2g_3g_4$. In both cases the even central elements of
$\Cl(e_1,\dots,e_n)$ are the scalars. (The $g_i$ exhibit an isomorphism
$\Cl_{0,n}\cong\Cl_{n,0}$ of superalgebras over $k$; over $\R$ this
requires $n$ to be a multiple of $4$.)

Let $\beta_0$ be the automorphism of $C$ with
\begin{equation}
\label{alt:eq:beta0}
 \beta_0(e_m)=e_{m+n},\qquad\beta_0(f_i)=g_i\quad(i\le n),\qquad
 \beta_0(f_{m+n})=f_m\qquad(m\ge1).
\end{equation}
It is well defined because the images are odd, anticommute and have the
right squares, and bijective because they generate $C$. By
\eqref{alt:eq:sigmaC},
\begin{equation}
\label{alt:eq:beta0sigma}
 \beta_0\circ\sigma_-^n=\sigma_+^n\quad\text{on } C .
\end{equation}
By Lemma~\ref{alt:lem:objects}(3), $\sigma_\pm^n=[\sigma_\pm^n,W_\pm]$ with
$W_\pm=\sigma_\pm^{n-1}(v_\pm)\cdots\sigma_\pm(v_\pm)\,v_\pm$.

\begin{lemma}
\label{alt:lem:beta}
Let $B=\beta_0(W_-)^{-1}W_+$. Then the conditions \eqref{alt:eq:E} hold for
$(\beta_0,B)$, and $\beta=[\beta_0,B]$ is an automorphism of $L$ with
\begin{equation}
 \beta\circ\sigma_-^n=\sigma_+^n .
\end{equation}
It has shift $(n,-n)$.
\end{lemma}
\begin{proof}
Lemma~\ref{alt:lem:objects}(1) for $\sigma_\pm^n$ and
\eqref{alt:eq:beta0sigma} give, for $x\in C_1$ and $a\in C$,
\begin{align}
 BxB^{-1}&=\beta_0(W_-)^{-1}\,\sigma_+^n(x)\,\beta_0(W_-)
 =\beta_0\bigl(W_-^{-1}\sigma_-^n(x)W_-\bigr)=\beta_0(x),\\
 B\,\rho(\sigma_+^n(a))\,B^{-1}&=\beta_0(W_-)^{-1}\,\sigma_+^n(\rho(a))\,\beta_0(W_-)
 =\beta_0\bigl(W_-^{-1}\sigma_-^n(\rho(a))W_-\bigr)=\beta_0\bigl(\rho(\sigma_-^n(a))\bigr).
 \nonumber
\end{align}
The first line is the first condition of \eqref{alt:eq:E}. The second
line says that
\begin{equation}
\label{alt:eq:Brho}
 B\,\rho(\beta_0(y))\,B^{-1}=\beta_0(\rho(y))
\end{equation}
for $y\in\sigma_-^n(C)$. It also holds for $y=f_i$ with $i\le n$: by
Proposition~\ref{alt:prop:sigma}(2) and Lemma~\ref{alt:lem:objects}(2) we
have $W_+\rho(g_i)W_+^{-1}=g_i$ and $W_-^{-1}f_iW_-=\rho(f_i)$, so
$B\rho(g_i)B^{-1}=\beta_0(W_-^{-1}f_iW_-)=\beta_0(\rho(f_i))$. As $C$ is
generated by $\sigma_-^n(C)$ and $f_1,\dots,f_n$, \eqref{alt:eq:Brho}
holds for all $y\in C$, which is the second condition of
\eqref{alt:eq:E}. By Lemma~\ref{alt:lem:objects}(3),
$\beta\sigma_-^n=[\beta_0\sigma_-^n,\beta_0(W_-)B]=[\sigma_+^n,W_+]=\sigma_+^n$.

The image of $\beta$ contains $C$, hence also $s=B^{-1}\beta(s)$ and $t$,
so $\beta$ is surjective. It preserves the $\Z$-degree. If
$\beta(cs^d)=0$ with $c\in C$, then $\beta_0(cP_d)=\beta(cs^dt^d)=0$, so
$cP_d=0$ and $cs^d=cP_ds^d=0$; similarly for $t^dc$. So $\beta$ is
injective by Lemma~\ref{alt:lem:L}(3).
\end{proof}

\subsection{Adjoining \texorpdfstring{$U$}{U}}

Let $\varepsilon\in\Z/2\Z$; it will be fixed in Lemma~\ref{alt:lem:t}.

\begin{definition}
\label{alt:def:Lprime}
$L'$ is the superalgebra $\bigoplus_{j\in\Z}L\,U^j$, where $U$ is
invertible of parity $\varepsilon$ and
$Ua=\beta(\gamma^\varepsilon(a))\,U$ for $a\in L$.
\end{definition}

\begin{lemma}
\label{alt:lem:crossed}
(1) Let $\tau$ be an object of $\EndCat(L)$ and let $c\in L$ be even and
invertible with $\tau\circ\beta=\Ad(c)\circ\beta\circ\tau$. Then $\tau$
extends to an object of $\EndCat(L')$ with $\tau(U)=cU$.

(2) Let $\tau,\tau'$ be objects of $\EndCat(L)$ with $\tau'$ injective,
and let $c,c',d\in L$ be even and invertible with
$\tau'\circ\beta=\Ad(c')\circ\beta\circ\tau'$,
$\tau'\tau\circ\beta=\Ad(d)\circ\beta\circ\tau'\tau$ and
$\tau'(c)\,c'=d$. Then $\tau\circ\beta=\Ad(c)\circ\beta\circ\tau$.

(3) Let $\tau,\tau'$ be objects of $\EndCat(L')$ which preserve $L$, and
let $x=\sum_jx_jU^j$ with $x_j\in L$ be an element of $\Hom(\tau,\tau')$.
Then $x_j\in\Hom(\beta^j\circ\tau,\tau')$ in $\EndCat(L)$.
\end{lemma}
\begin{proof}
(1) The element $cU$ is invertible of parity $\varepsilon$, and since
$\tau$ commutes with $\gamma$,
\begin{equation}
 cU\tau(a)(cU)^{-1}=c\,\beta\gamma^\varepsilon\tau(a)\,c^{-1}
 =\tau(\beta\gamma^\varepsilon(a)) .
\end{equation}

(2) We have
$\tau'\tau\beta=\Ad(d)\beta\tau'\tau=\Ad(dc'^{-1})\tau'\beta\tau
=\tau'\circ\Ad(c)\beta\tau$, and $\tau'$ is injective.

(3) Let $x$ have parity $p$, so that $x_j$ has parity $p+j\varepsilon$.
Comparing the coefficients of $U^j$ in
$x\tau(a)=\gamma^p(\tau'(a))x$ gives
\begin{equation}
 x_j\,\beta^j\gamma^{j\varepsilon}\tau(a)=\gamma^p(\tau'(a))\,x_j\qquad(a\in L);
\end{equation}
now replace $a$ by $\gamma^{j\varepsilon}(a)$.
\end{proof}

Let $p_j=\frac1{\sqrt2}(1+e_je_{j+1})$, so that conjugation by $p_j$ sends
$e_j\mapsto-e_{j+1}$, $e_{j+1}\mapsto e_j$ and fixes the other generators.
For $\kappa\in k^\times$, to be fixed in Lemma~\ref{alt:lem:t}, put
\begin{equation}
 c_+=\kappa\,p_1p_2\cdots p_n\in\Cl(e_1,\dots,e_{n+1}),\qquad
 c_-=c_+^{-1}B^{-1} .
\end{equation}
Then $\Ad c_+$ sends $e_i\mapsto-e_{i+1}$ for $i\le n$ and
$e_{n+1}\mapsto e_1$, and fixes the other generators. In particular
\begin{equation}
\label{alt:eq:cplus}
 c_+xc_+^{-1}=\sigma_+(x)\qquad\text{for } x\in\Cl(e_1,\dots,e_n,f_1,f_2,\dots).
\end{equation}

\begin{proposition}
\label{alt:prop:U}
The objects $\sigma_+$, $\sigma_-$ and $\rho$ of $\EndCat(L)$ extend to
objects of $\EndCat(L')$ with
\begin{equation}
 \sigma_+(U)=c_+U,\qquad\sigma_-(U)=c_-U,\qquad
 \rho(x)=s_0xt_0+s_1\gamma(x)t_1\quad(x\in L'),
\end{equation}
and these satisfy:

(1) $\sigma_+\circ\sigma_-=\sigma_-\circ\sigma_+=\rho$ on $L'$;

(2) $e_1\in\Hom(\sigma_+,\sigma_+)_1$ and
$f_1\in\Hom(\sigma_-,\sigma_-)_1$ in $\EndCat(L')$.
\end{proposition}
\begin{proof}
\emph{Extension of $\sigma_+$.} We claim
$\sigma_+\beta=\Ad(c_+)\beta\sigma_+$. First,
$\sigma_+\beta\sigma_-^n=\sigma_+^{n+1}$, while
$\Ad(c_+)\beta\sigma_+\sigma_-^n=\Ad(c_+)\beta\sigma_-^n\sigma_+
=\Ad(c_+)\sigma_+^{n+1}=\sigma_+^{n+1}$, because $c_+$ is an even element
of $\Hom(\sigma_+^{n+1},\sigma_+^{n+1})$ by
Proposition~\ref{alt:prop:sigma}(2). Second, for $i\le n$, by
\eqref{alt:eq:cplus},
$\sigma_+\beta(f_i)=\sigma_+(g_i)=c_+g_ic_+^{-1}=c_+\,\beta\sigma_+(f_i)\,c_+^{-1}$.
So the two sides agree on $\sigma_-^n(L)$ and on $f_1,\dots,f_n$. These
generate $L$: indeed $\sigma_-^n(L)$ contains $e_m=\pm\sigma_-^n(e_m)$ and
$f_{m+n}=\sigma_-^n(f_m)$, and $s=W_-^{-1}\sigma_-^n(s)$,
$t=\sigma_-^n(t)W_-$. Now Lemma~\ref{alt:lem:crossed}(1) applies.

\emph{The object $\rho$.} The formula defines an object of $\EndCat(L')$
extending $\rho$, with
\begin{equation}
 \rho(U)=\sum_is_i\gamma^i(U)t_i=\sum_is_i\beta(t_i)\,U=B^{-1}U .
\end{equation}

\emph{Extension of $\sigma_-$, and (1).} By
Lemma~\ref{alt:lem:objects}(1), $\rho\beta=\Ad(B^{-1})\beta\rho$ on $L$.
The conditions \eqref{alt:eq:E} for $(\beta_0,B)$ give
$Be_1B^{-1}=e_{n+1}$ and $Bf_1B^{-1}=g_1$, and, applied to $a=e_m$ and
$a=f_{m+n}$, they show that $B$ commutes with $e_j$ for $j\ge n+2$ and
with $f_j$ for $j\ge2$. Since $c_+$ commutes with these $e_j$ and with
all $f_j$, and $c_+e_{n+1}c_+^{-1}=e_1$, the even element
$c_-=c_+^{-1}B^{-1}$ commutes with $e_j$ for $j\ge n+1$ and with $f_j$ for
$j\ge2$. By Lemma~\ref{alt:lem:commutant},
$c_-\in\Cl(e_1,\dots,e_n,f_1)$, and then \eqref{alt:eq:cplus} gives
\begin{equation}
 \sigma_+(c_-)\,c_+=c_+c_-=B^{-1}.
\end{equation}
By Lemma~\ref{alt:lem:crossed}(2), applied to $\tau'=\sigma_+$, which is
injective because $\sigma_-\sigma_+=\rho$ is, and to $\tau=\sigma_-$, we
get $\sigma_-\beta=\Ad(c_-)\beta\sigma_-$, so that $\sigma_-$ extends.
Moreover $\sigma_+\sigma_-(U)=\sigma_+(c_-)c_+U=B^{-1}U$ and
$\sigma_-\sigma_+(U)=\sigma_-(c_+)c_-U=c_+c_-U=B^{-1}U$, since $\sigma_-$
fixes the even elements of $\Cl(e_1,e_2,\dots)$. With
Proposition~\ref{alt:prop:sigma}(1) this proves (1).

(2) By Proposition~\ref{alt:prop:sigma}(2) it remains to check that
$e_1\sigma_+(U)=(-1)^\varepsilon\sigma_+(U)e_1$ and
$f_1\sigma_-(U)=(-1)^\varepsilon\sigma_-(U)f_1$. As
$(-1)^\varepsilon Uy=\beta(y)U$ for odd $y\in L$, this means
$e_1c_+=c_+\beta(e_1)$ and $f_1c_-=c_-\beta(f_1)$, and indeed
\begin{equation}
 c_+\beta(e_1)c_+^{-1}=c_+e_{n+1}c_+^{-1}=e_1,\qquad
 c_-\beta(f_1)c_-^{-1}=c_+^{-1}B^{-1}g_1Bc_+=c_+^{-1}f_1c_+=f_1 .
 \qedhere
\end{equation}
\end{proof}

By Proposition~\ref{alt:prop:U}, $\sigma_\pm^n(U)=C_\pm U$ with
$C_\pm=\sigma_\pm^{n-1}(c_\pm)\cdots\sigma_\pm(c_\pm)\,c_\pm\in L$. These
elements do not depend on $\varepsilon$.

\begin{lemma}
\label{alt:lem:t}
There is $t\in k^\times$ with $C_+=t\,\beta(C_-)$. One can choose $\kappa$
so that $t=1$ if $k=\C$, and $t=\pm1$ if $k=\R$. With this choice, and
with $\varepsilon$ defined by $(-1)^\varepsilon=t$, the element $U$ lies
in $\Hom(\sigma_-^n,\sigma_+^n)_\varepsilon$ in $\EndCat(L')$. Thus
$\sigma_+^n$ and $\sigma_-^n$ are isomorphic.
\end{lemma}
\begin{proof}
Applying $\sigma_\pm^n$ to $UaU^{-1}=\beta\gamma^\varepsilon(a)$ gives
$\sigma_\pm^n\beta=\Ad(C_\pm)\beta\sigma_\pm^n$ on $L$. Applying $\beta$ to
the second of these and using $\beta\sigma_-^n=\sigma_+^n$, we get
$\sigma_+^n\beta=\Ad(\beta(C_-))\beta\sigma_+^n$ as well. Hence
$z=\beta(C_-)^{-1}C_+$ is an even element of
$\Hom(\beta\sigma_+^n,\beta\sigma_+^n)=\beta(\Cl(e_1,\dots,e_n))$.

By Proposition~\ref{alt:prop:U}(2), in $\EndCat(L')$ the space
$\Hom(\sigma_+^n,\sigma_+^n)$ contains $\Cl(e_1,\dots,e_n)$ and
$\Hom(\sigma_-^n,\sigma_-^n)$ contains $\Cl(f_1,\dots,f_n)$.
Evaluating the intertwining relations on $U$, as in the proof of
Proposition~\ref{alt:prop:U}(2), this means
\begin{equation}
 x=C_+\,\beta(x)\,C_+^{-1}\quad(x\in\Cl(e_1,\dots,e_n)),\qquad
 y=C_-\,\beta(y)\,C_-^{-1}\quad(y\in\Cl(f_1,\dots,f_n)).
\end{equation}
Applying $\beta$ to the second relation, and using that
$\beta(\Cl(f_1,\dots,f_n))=\Cl(g_1,\dots,g_n)=\Cl(e_1,\dots,e_n)$, we get
$x=\beta(C_-)\,\beta(x)\,\beta(C_-)^{-1}$ for $x\in\Cl(e_1,\dots,e_n)$.
Hence $z$ commutes with $\beta(\Cl(e_1,\dots,e_n))$. So $z$ is the image
under $\beta$ of an even central element of $\Cl(e_1,\dots,e_n)$, that
is, a scalar $t$.

Replacing $\kappa$ by $\kappa\xi$ multiplies $c_+$ by $\xi$, $c_-$ by
$\xi^{-1}$, $C_\pm$ by $\xi^{\pm n}$ and $t$ by $\xi^{2n}$, which gives the
normalization. Finally, for $a\in L$,
$U\sigma_-^n(a)=\beta\gamma^\varepsilon\sigma_-^n(a)\,U
=\gamma^\varepsilon(\sigma_+^n(a))\,U$, and
\begin{equation}
 U\sigma_-^n(U)=UC_-U=\beta(C_-)\,U^2=t^{-1}C_+U^2
 =(-1)^\varepsilon\sigma_+^n(U)\,U=\gamma^\varepsilon(\sigma_+^n(U))\,U .
 \qedhere
\end{equation}
\end{proof}

From now on $\kappa$ and $\varepsilon$ are as in Lemma~\ref{alt:lem:t}.

\begin{proposition}
\label{alt:prop:homs}
In $\EndCat(L')$:

(1) $\Hom(\sigma_+^p\sigma_-^q,\sigma_+^{p'}\sigma_-^{q'})=0$ unless
$p-q\equiv p'-q'$ modulo $2n$;

(2) $\Hom(\sigma_+^m,\sigma_+^m)=\Cl(e_1,\dots,e_m)$,
$\Hom(\sigma_-^m,\sigma_-^m)=\Cl(f_1,\dots,f_m)$ and $\Hom(\rho,\rho)=C_1$,
for $m\ge0$.
\end{proposition}
\begin{proof}
Put $\tau=\sigma_+^p\sigma_-^q$ and $\tau'=\sigma_+^{p'}\sigma_-^{q'}$, and
let $x=\sum_jx_jU^j$ be a morphism from $\tau$ to $\tau'$. By
Lemma~\ref{alt:lem:crossed}(3), the coefficient $x_j$ lies in
$\Hom(\beta^j\tau,\tau')$ in $\EndCat(L)$. The shifts of these two objects
are $(p+jn,q-jn)$ and $(p',q')$, so $x_j=0$ by \eqref{alt:eq:index} unless
$2nj=(p'-q')-(p-q)$. This proves (1), and for $\tau=\tau'$ it shows that
$\Hom(\tau,\tau)$ in $\EndCat(L')$ is contained in the same space in
$\EndCat(L)$, which is given in Proposition~\ref{alt:prop:sigma}(2). The
converse inclusions follow from Proposition~\ref{alt:prop:U}(2), and from
the fact that $s_it_j$ lies in $\Hom(\rho,\rho)$.
\end{proof}

\subsection{The reference category}

Let $\mathcal E$ be the full subsupercategory of $\EndCat(L')$ on the
objects $\sigma_+^p\sigma_-^q$, $p,q\ge0$. By
Proposition~\ref{alt:prop:U}(1) it is closed under composition, so it is
a strict monoidal supercategory.

\begin{definition}
\label{alt:def:R}
$\mathcal R$ is the completion $\overline{\mathcal E}$ of $\mathcal E$
(Definition~\ref{def:completions}), that is, the idempotent completion
of the additive envelope of $\mathcal E$ with respect to even morphisms,
with the monoidal structure of
Proposition~\ref{prop:monoidal-completion}.
For $\tau\in\mathcal E$ and an even idempotent $q$ of $\Hom(\tau,\tau)$,
we write $(\tau,q)$ for the corresponding object of $\mathcal R$; the
morphisms from $(\tau,q)$ to $(\tau',q')$ are the elements of
$q'\Hom(\tau,\tau')\,q$, and the tensor product is chosen such that
$(\tau,q)\otimes(\tau',q')=(\tau\tau',q\,\tau(q'))$. We write $\tau$ for
$(\tau,1)$ and put
\begin{equation}
 \one=\id,\qquad\pi=(\rho,s_1t_1),\qquad
 \zeta=t_1\in\Hom(\pi,\one)_1 ,
\end{equation}
an odd isomorphism with inverse $s_1$.
\end{definition}

\begin{remark}
\label{alt:rem:sums}
The formal completion is only a convenience: direct sums and parity
shifts already exist in $\EndCat(L)$, and in the same way in
$\EndCat(L')$.

(1) \emph{Direct sums.} The elements
\begin{equation}
 a_0=s_0,\qquad a_1=s_1e_1,\qquad b_0=t_0,\qquad b_1=e_1t_1
\end{equation}
are even and satisfy $b_ia_j=\delta_{ij}$ and $a_0b_0+a_1b_1=1$. Hence for
any two objects $\tau_0,\tau_1$ the formula
\begin{equation}
 \tau(x)=a_0\,\tau_0(x)\,b_0+a_1\,\tau_1(x)\,b_1
\end{equation}
defines an object, and the even morphisms $a_i\in\Hom(\tau_i,\tau)$ and
$b_i\in\Hom(\tau,\tau_i)$ make it a direct sum of $\tau_0$ and $\tau_1$.

(2) \emph{Parity shift.} Let $\pi'=\Ad(e_1)\circ\gamma$, the automorphism
which on $C$ changes the sign of $e_1$ and fixes the other generators.
Then $e_1\in\Hom(\pi',\id)_1$ is an odd isomorphism, and for every object
$\tau$ the objects $\pi'\circ\tau$ and $\tau\circ\pi'$ are oddly
isomorphic to $\tau$, via $e_1$ and $\tau(e_1)$.

(3) With these, $\rho(x)=a_0\,x\,b_0+a_1\,\pi'(x)\,b_1$, that is,
$\rho=\id\oplus\pi'$ in the sense of (1), and $a_1$ is an even
isomorphism from $\pi'$ to the object $\pi=(\rho,s_1t_1)$ of $\mathcal R$.
One can check that the images of even idempotents lying in $C$ exist in
$\EndCat(L)$ as well, so that $\mathcal R$ is equivalent, as a
supercategory, to a full subcategory of $\EndCat(L')$. We do not need
this.
\end{remark}

For $k=\R$ the superalgebras $\Cl(e_1,\dots,e_m)\cong\Cl_{m,0}$ and
$\Cl(f_1,\dots,f_m)\cong\Cl_{0,m}$ have Brauer--Wall classes $-m$ and $m$,
respectively, in these coordinates. We use the structure theory of
central simple superalgebras \cite{Wall}: such an algebra $A$ has a
semisimple even part; if $q$ is a minimal even idempotent, $qAq$ is the
division superalgebra in the class of $A$; and for two minimal even
idempotents $q,q'$ the space $q'Aq$ contains a homogeneous element $u$
with inverse in $qAq'$.

\begin{theorem}
\label{alt:thm:R}
$(\mathcal R,\pi,\zeta)$ is a Brauer--Wall reference category over $k$.
Its simple objects are, up to isomorphism in $\mathcal R$, the object
$\one$ of type $0$, the summands of $\sigma_+^m$ for $1\le m\le n$, of
type $-m$, and the summands of $\sigma_-^m$ for $1\le m\le n-1$, of type
$m$.
\end{theorem}
\begin{proof}
Let $T_m=\sigma_+^m$ for $1\le m\le n$ and $T_{-m}=\sigma_-^m$ for
$1\le m\le n-1$, and let $A_j=\Hom(T_j,T_j)$. Call two objects
\emph{oddly isomorphic} if there is an odd isomorphism between them in
$\mathcal R$.

\emph{Step 1: every object of $\mathcal E$ is isomorphic in
$\underline{\mathcal R}$ to a direct sum of copies of $\one$, $\pi$ and
the $T_j$.} Let us say that such a direct sum is \emph{of standard form}.
An object which is oddly isomorphic to an object of standard form is
itself isomorphic in $\underline{\mathcal R}$ to one: $\one$ and $\pi$
are oddly isomorphic, and an odd isomorphism to $T_j$ can be composed
with the odd automorphism $e_1$ or $f_1$ of $T_j$. Now let
$\tau\in\mathcal E$. The even idempotents $\tau(s_0t_0)$ and
$\tau(s_1t_1)$ of $\Hom(\tau\rho,\tau\rho)$ add up to $1$, and
$\tau(s_0)$ and $\tau(s_1)$ are an even and an odd isomorphism from
$\tau$ to the corresponding summands of $\tau\rho$. So if $\tau$ is
isomorphic to an object of standard form, then so is $\tau\rho$. We prove
the claim for $\sigma_+^p\sigma_-^q$ by induction on $p+q$, treating the
case $p,q\ge1$ first: then the object is $\tau\rho$ with
$\tau=\sigma_+^{p-1}\sigma_-^{q-1}$.
If $q=0$ and $p\le n$, or $p=0$ and $q\le n-1$, it is $\one$ or some
$T_j$. Otherwise we use Lemma~\ref{alt:lem:t}: the invertible element
$\sigma_+^{p-n}(U)$ is an isomorphism, even or odd, from
$\sigma_+^{p-n}\sigma_-^n$ to $\sigma_+^p$, and $\sigma_-^{q-n}(U)$ is one
from $\sigma_-^q$ to $\sigma_-^{q-n}\sigma_+^n$. The latter is $T_n$ if
$q=n$, and in the remaining cases we are back in the first case, with the
same $p+q$.

\emph{Step 2: semisimplicity and types.} By
Proposition~\ref{alt:prop:homs} there are no non-zero morphisms between
$T_i$ and $T_j$ for $i\ne j$, nor between $T_j$ and $\rho$, and hence
none between $T_j$ and $\one$ or $\pi$. Moreover $A_j$ is a Clifford
algebra, $\Hom(\one,\one)=k$ and $\Hom(\pi,\pi)=s_1t_1C_1s_1t_1=k$; and
$\Hom(\id,\rho)=\Hom(\rho,\rho)\,s_0=ks_0\oplus ks_1$, so that
$\Hom(\one,\pi)_0=s_1t_1\cdot ks_0=0$. Since the even part of each $A_j$
is semisimple, the full subcategory of $\underline{\mathcal R}$ on the
direct sums of $\one$, $\pi$ and the objects $(T_j,q)$ is semisimple,
with simple objects $\one$, $\pi$ and $(T_j,q)$ for $q$ a minimal even
idempotent of $A_j$. It is closed under direct summands, so by Step~1 it
is all of $\underline{\mathcal R}$. By the structure theory, the
endomorphism superalgebra $qA_jq$ of $(T_j,q)$ is the division
superalgebra of the class of $A_j$, which is $-j$, and $(T_j,q)$ and
$(T_j,q')$ are isomorphic in $\mathcal R$. (For $k=\C$ there is only
$T_1=\sigma_+$, with $A_1=\Cl(e_1)\cong D_1$ and $q=1$.)
The classes of $\one$ and of the $T_j$ are
$0,-1,\dots,-n,1,\dots,n-1$, which are all the elements of $\Z/2n\Z$,
each once. So every simple object has a unique type, every type occurs,
and simple objects of the same type are isomorphic in $\mathcal R$.

\emph{Step 3: rigidity.} By Proposition~\ref{alt:prop:U}(1), $s_0$ is an
even morphism from $\one$ to $\sigma_+\sigma_-=\sigma_-\sigma_+$ and $t_0$
is one in the other direction, and by
Proposition~\ref{alt:prop:sigma}(3) the four zigzags
$t_0\,\sigma_\pm(s_0)$ and $\sigma_\pm(t_0)\,s_0$ are equal to
$1/\sqrt2$. So $\sigma_-$ is a left and a right dual of $\sigma_+$. Hence
all objects of $\mathcal E$ are rigid, and so are their direct sums and
direct summands. If $X$ is a simple summand of $\sigma_\pm^m$, then $X^*$
is a summand of $\sigma_\mp^m$ and $X\otimes X^*$ is a summand of
$\rho^m$, which is a direct sum of copies of $\one$ and $\pi$ by Step~1.

Thus $\underline{\mathcal R}$ is a semisimple rigid monoidal category
with finitely many simple objects and $\End(\one)=k$, every simple object
is weakly invertible, and the conditions on the types hold.
\end{proof}

\begin{remark}
\label{alt:rem:scalar}
(1) The scalar $t$ of Lemma~\ref{alt:lem:t} is not needed explicitly, but
it can be computed. For $k=\C$ and $g_1=\pm\sqrt{-1}\,e_1$ one finds
$t=\mp\sqrt{-1}\,\kappa^2$. For $k=\R$ and $g_i=e_i\omega$ one finds
$t=\kappa^8$, so that $\kappa=\pm1$ and $U$ is even; if $g_1$ is replaced
by $-g_1$, which represents the other class of isomorphisms
$\Cl_{0,4}\cong\Cl_{4,0}$ up to even inner automorphisms, then
$t=-\kappa^8$ and $U$ is odd.

(2) On the chain $\cdots f_3\,f_2\,f_1\,|\,e_1\,e_2\,e_3\cdots$ the
endomorphism $\sigma_+$ translates the right half and $\sigma_-$ the left
half, each leaving one generator behind, which generates its odd
endomorphisms; $\rho$ translates both halves, and $s$ fills the two sites
left behind with the vacuum of the pair $e_1,f_1$. The automorphism
$\beta$ translates the whole chain by $n$ sites to the right, and the $n$
generators which cross the wall are converted by
\eqref{alt:eq:conversion}. The difference of the two components of the
shift is additive, it vanishes on $L$ by Lemma~\ref{alt:lem:degree}, and
it is $2n$ for $\beta$; this is the grading of $\mathcal R$ by
$\Z/2n\Z$.
\end{remark}

\section*{Acknowledgements}
The author is deeply grateful to Yuji Tachikawa for his generous support
and encouragement throughout this work.
The author also thanks Sean Sanford for insightful comments.
This research was supported by FoPM, WINGS Program, the University of Tokyo.
This work was supported by JST SPRING, Japan Grant Number JPMJSP2108.

\section*{Use of AI tools}
AI tools were used to draft proofs from the statements, to carry out
computations, especially in Sections~\ref{sec:stacking}
and~\ref{sec:complex-real}, and to formalize the results for $k=\R$ in Lean.
The author takes full responsibility for the paper.

\bigskip
\noindent
\begin{minipage}{\textwidth}
\small
\textsc{Kavli Institute for the Physics and Mathematics of the Universe (WPI),
The University of Tokyo, Kashiwa, Chiba 277-8583, Japan}\\
\textit{E-mail address}:
\href{mailto:kotaro.kawasumi@gmail.com}{\texttt{kotaro.kawasumi@gmail.com}}
\end{minipage}

\end{document}